\documentclass[letterpaper]{amsart} 
 
\usepackage{amssymb,latexsym} 
 
\numberwithin{equation}{section} 
 
\newtheorem{theorem}{Theorem}[section] 
\newtheorem{proposition}[theorem]{Proposition} 
 
\newtheorem{corollary}[theorem]{Corollary} 
\newtheorem{lemma}[theorem]{Lemma} 
 
\theoremstyle{definition} 
 
\newtheorem{remark}[theorem]{Remark} 
\newtheorem{example}[theorem]{Example} 
\newtheorem{definition}[theorem]{Definition}

\newcommand{\credit}[1]{\smallskip\noindent {\textbf{#1.\ }}} 
 
\def\proof{\smallskip\noindent {\bf Proof.\ }} 
\def\endproof{\hfill$\square$\medskip} 
 
\def\AA{\mathcal{A}} 
 
\def\ZZ{\mathbb{Z}} 
\def\CC{\mathbb{C}} 
\def\PP{\mathbb{P}} 
\def\RR{\mathbb{R}} 
\def\QQ{\mathbb{Q}} 
 
\def\xx{\mathbf{x}} 
\def\pp{\mathbf{p}} 
\def\FFcal{\mathcal{F}} 
 
\def\subplus{{\,\uplus\,} } 
\def\verylongleftrightarrow{\longleftarrow\!\!\!\longrightarrow} 
\def\height{{\rm ht}} 
\def\Poly{\mathbf{P}} 
\def\Trop{{\rm Trop}}

\newcommand{\bmat}[4]{\left[\!\!\begin{array}{cc} 
#1 & #2 \\ #3 & #4 \\ 
\end{array}\!\!\right]} 
 
\begin{document} 
 
\title[Cluster algebras II] 
{Cluster algebras II: \\ 
Finite type classification 
}

\author{Sergey Fomin} 
\address{Department of Mathematics, University of Michigan, 
Ann Arbor, MI 48109, USA} \email{fomin@umich.edu}

\author{Andrei Zelevinsky} 
\address{\noindent Department of Mathematics, Northeastern University, 
 Boston, MA 02115} 
\email{andrei@neu.edu} 
 

\date{Received: September 7, 2002. Revised version: February~6, 2003} 
 
 \thanks{Research 
supported 
by NSF (DMS) grants \# 0070685 (S.F.), 9971362, and 0200299 (A.Z.).} 
 
 
 
\maketitle 
 
\tableofcontents 
 
\vspace{-.3in} 
 
\pagebreak[3] 
 
 
\section{Introduction and main results} 
\label{sec:intro-main-results} 
 
\subsection{Introduction} 
 
The origins of cluster algebras, 
first introduced 
in~\cite{fz-clust1}, 
lie in the desire to understand, 
in concrete algebraic and combinatorial terms, the structure of 
``dual canonical bases'' in (homogeneous) 
coordinate rings of various algebraic varieties 
related to semisimple groups. 
Several classes of such varieties---among them Grassmann~and Schubert 
varieties, base affine spaces, 
and double Bruhat cells---are expected (and in many cases proved) to 
carry a cluster algebra structure. 
This structure includes the description of the ring in question 
as a commutative ring generated inside its ambient field by a 
distinguished family of generators called cluster variables. 
Even though most of the rings of interest to us are finitely 
generated, their set of cluster variables may well be infinite. 
A cluster algebra has \emph{finite type} if it has 
a finite number of cluster variables. 
 
The main result of this paper (Theorem~\ref{th:fin-type-class}) 
provides a complete classification of the cluster algebras of finite 
type. 
This classification turns out to be identical 
to the Cartan-Killing 
classification of semisimple Lie algebras and finite root systems. 
%
%
This result is particularly intriguing since 
in most cases, the symmetry exhibited by the Cartan-Killing type 
of a cluster algebra is not apparent at all from its geometric 
realization. 
For instance, the coordinate ring of the base affine space of the 
group $SL_5$ turns out to be a cluster algebra of the 
Cartan-Killing type~$D_6\,$. 
Other examples of similar nature can be found in 
Section~\ref{sec:geometric-realization}, 
in which we show how cluster algebras of types $ABCD$ arise as 
coordinate rings of some classical 
algebraic varieties. 
 
In order to understand a cluster algebra of finite type, 
one needs to study the combinatorial structure behind it, 
which is captured by its \emph{cluster complex}. 
Roughly speaking, it is defined as follows. 
The cluster variables for a given cluster algebra 
are not given from the outset 
but are obtained from some 
initial ``seed'' by an explicit process of ``mutations''; 
each mutation exchanges a cluster variable in the current seed by 
a new cluster variable according to a particular set of rules. 
In a cluster algebra of finite type, this process ``folds'' to 
produce a finite set of seeds, each containing the same number $n$ of 
cluster variables (along with some extra information needed to perform 
mutations). These $n$-element collections of cluster variables, 
called \emph{clusters}, are the maximal faces of the (simplicial) 
cluster complex.

In Theorem~\ref{th:finite-type-complex}, we identify this complex 
as the dual simplicial complex of a \emph{generalized associahedron} 
associated with the corresponding root system. 
These complexes (indeed, convex polytopes~\cite{cfz}) 
were introduced in \cite{yga} in relation to 
our proof of Zamolodchikov's periodicity conjecture for algebraic 
\hbox{$Y$-systems.} 
A generalized associahedron of type $A$ is the usual 
associahedron, or the Stasheff polytope \cite{stasheff}; 
in types $B$ and~$C$, it is the 
cyclohedron, or the Bott-Taubes polytope~\cite{bott-taubes,simion-B}. 
 
One of the crucial steps in our proof of the classification theorem 
is a new combinatorial characterization of Dynkin diagrams. 
In Section~\ref{sec:diagram-mutations}, 
we introduce an equivalence relation, 
called \emph{mutation equivalence}, 
on finite directed graphs with weighted edges. 
We then prove 
that a connected graph 
$\Gamma$ is mutation equivalent to an orientation of a Dynkin diagram 
if and only if every graph that is mutation equivalent to  $\Gamma$ 
has all edge weights $\leq 3$. 
We do not see a direct way to relate this description to any 
previously known 
characterization of the Dynkin diagrams. 
 
\pagebreak[2]

We already mentioned that the initial motivation for the study of 
cluster algebras came from representation theory; 
see \cite{camblec} for a more detailed discussion of the 
representation-theoretic context. 
Another source of inspiration was Lusztig's theory of 
total positivity in semisimple Lie groups, which 
was further developed in a series of papers of the 
present authors and their collaborators (see, e.g.,~\cite{Lu2, FZ1} and 
references therein). 
The mutation mechanism used for creating new cluster variables from the initial 
``seed" was designed to ensure that in concrete geometric realizations, these 
variables become regular functions taking positive values on all 
totally positive elements of a variety in question, 
a property that the elements of the dual canonical basis are known to
possess~\cite{Lu1}.  
 
Following the foundational paper \cite{fz-clust1}, several 
unexpected connections and appearances of cluster algebras have been 
discovered and explored. 
They included: Laurent phenomena in number theory and 
combinatorics~\cite{fz-laurent}, 
$Y$-systems and thermodynamic Bethe ansatz~\cite{yga}, 
quiver representations \cite{mrz}, 
and Poisson geometry~\cite{GSV}.

While this text belongs to an ongoing series of papers devoted to cluster 
algebras and related topics, it is designed to be read independently of 
any other publications on the subject. 
Thus, all definitions and results from \cite{fz-clust1,yga,cfz} 
that we need are presented ``from scratch'', in the form most suitable for our 
current purposes. 
In particular, the core concept of a normalized cluster algebra 
\cite{fz-clust1} is defined anew in 
Section~\ref{sec:def-of-normalized-ca}, while 
Section~\ref{sec:root-systems-gen-associahedra} provides 
the relevant background on generalized associahedra~\cite{yga,cfz}.

The main new results 
(Theorems~\ref{th:fin-type-class}-\ref{th:finite-type-complex}) 
are stated in 
Sections~\ref{sec:finite-type-classification}--\ref{sec:cluster-complexes}. 
The organization of the rest of the paper is outlined in 
Section~\ref{sec:organization-ca2}. 
 
 
\pagebreak[2] 
 
\subsection{Basic definitions} 
\label{sec:def-of-normalized-ca} 
 
We start with 
the definition of a (normalized) cluster 
algebra $\AA$ (cf.~\cite[Sections~2 and~5]{fz-clust1}). 
This is a commutative ring embedded in an ambient field~$\FFcal$ 
defined as follows. 
Let $(\PP,\oplus, \cdot)$ be a \emph{semifield}, i.e., 
an abelian multiplicative group supplied with an \emph{auxiliary addition} 
$\oplus$ which is commutative, associative, and 
distributive with respect to the multiplication in~$\PP$. 
The following example (see \cite[Example~5.6]{fz-clust1}) will be of 
particular importance to us: 
let $\PP$ be a free abelian group, written multiplicatively, 
with a finite set of generators $p_j \, (j \in J)$, and with 
auxiliary addition $\oplus$ given by 
\begin{equation} 
\label{eq:tropical addition} 
\prod_j p_j^{a_j} \oplus \prod_j p_j^{b_j}  = 
\prod_j p_j^{\min (a_j, b_j)}  \ . 
\end{equation} 
We denote this semifield by $\Trop (p_j: j \in J)$. 
The multiplicative group of any semifield $\PP$ is torsion-free 
\cite[Section~5]{fz-clust1}, hence its group ring $\ZZ\PP$ is a 
domain. 
As an \emph{ambient field} for~$\AA$, we take a field $\FFcal$ 
isomorphic to the field of rational functions in $n$ independent 
variables 
(here $n$ is the \emph{rank} of~$\AA$), 
with coefficients in~$\ZZ \PP$. 
 
\pagebreak[2] 
 
A \emph{seed} in $\FFcal$ is a triple $\Sigma = (\xx, \pp, B)$, where 
 
\begin{itemize} 
 
\item $\xx$ is an $n$-element subset of $\FFcal$ which is a 
transcendence basis over the field of fractions of $\ZZ \PP$. 
 
\item $\pp = (p_x^\pm)_{x \in \xx}$ is a $2n$-tuple of elements of $\PP$ satisfying the 
\emph{normalization condition} 
$p_x^+ \oplus p_x^-  = 1$ for all $x \in \xx$. 
 
\item $B = (b_{xy})_{x,y \in \xx}$ is an $n\!\times\! n$ integer matrix 
with rows and columns indexed by~$\xx$, which is 
\emph{sign-skew-symmetric} \cite[Definition~4.1]{fz-clust1}; that~is, 
\begin{equation} 
\label{eq:sss} 
\text{for every $x,y \in \xx$, either $b_{xy} = b_{yx} = 0$, or 
$b_{xy} b_{yx} < 0$.} 
\end{equation} 
\end{itemize} 
 
We will need to recall the notion of \emph{matrix mutation} 
\cite[Definition~4.2]{fz-clust1}. 
Let 
$B = (b_{ij})
$ 
and $B' = (b'_{ij})
$ 
be 
real square matrices of the same size. 
We say that $B'$ is obtained from $B$ by a \emph{matrix mutation} in direction~$k$ 
if 
\begin{equation} 
\label{eq:mutation} 
b'_{ij} = 
\begin{cases} 
-b_{ij} & \text{if $i=k$ or $j=k$;} \\[.05in] 
b_{ij} + \displaystyle\frac{|b_{ik}| b_{kj} + 
b_{ik} |b_{kj}|}{2} & \text{otherwise.} 
\end{cases} 
\end{equation} 

\begin{definition} 
\label{def:seed-mutation} 
\emph{(Seed mutations)} 
Let $\Sigma = (\xx, \pp, B)$ be a seed in $\FFcal$, as above, 
and let $z \in \xx$. 
Define the triple 
$\overline \Sigma = (\overline \xx, \overline \pp, \overline B)$ as 
follows: 
\begin{itemize} 
 
\item 
$\overline \xx = \xx - \{z\} \cup \{\overline z\}$, where $\overline z \in \FFcal$ is determined 
by the \emph{exchange relation} 
\begin{equation} 
\label{eq:exchange-rel-xx} 
z\,\overline z 
=p_z^+ \, \prod_{\substack{x\in\xx \\ b_{xz}>0}} x^{b_{xz}} 
+p_z^- \, \prod_{\substack{x\in\xx \\ b_{xz}<0}} x^{-b_{xz}} 
\end{equation} 
 
\item 
the $2n$-tuple $\overline \pp=(\overline p ^\pm_x)_{x\in\overline 
  \xx}$ 
is uniquely determined by the  normalization conditions 
$\overline p_x^+ \oplus \overline p_x^-  = 1$ together with 
 \begin{equation} 
\label{eq:p-mutation} 
{\overline p ^+_x}/{\overline p ^-_x} = 
\begin{cases} 
{p^-_z}/{p^+_z} & \text{if $x=\overline z$};\\[.05in] 
(p_z^+)^{b_{zx}}{p ^+_x}/{p ^-_x}   & \text{if $b_{zx}\geq 0$};\\[.05in] 
(p_z^-)^{b_{zx}}{p ^+_x}/{p ^-_x}   & \text{if $b_{zx}\leq 0$}. 
\end{cases} 
\end{equation}

\item 
the matrix $\overline B$ is obtained from $B$ by applying the matrix 
mutation in direction~$z$ and then relabeling one row and one column 
by replacing $z$ by~$\overline z$. 
 
\end{itemize} 
If the triple $\overline \Sigma$  is again a seed in $\FFcal$ 
(i.e., if the matrix $\overline B$ is sign-skew-symmetric), then we say that 
$\Sigma$ admits a \emph{mutation} in the direction $z$ that results 
in~$\overline \Sigma$. 
\end{definition} 
 
We note that the exchange relation (\ref{eq:exchange-rel-xx}) is a 
reformulation of \cite[(2.2),(4.2)]{fz-clust1}, 
while the rule (\ref{eq:p-mutation}) is a rewrite of \cite[(5.4), 
(5.5)]{fz-clust1}. 
The elements $\overline p^\pm_x$ are determined by 
(\ref{eq:p-mutation}) uniquely since $p \oplus q = 1$ and 
$p/q = u$ imply $p = u/(1 \oplus u)$ and $q = 1/(1 \oplus u)$. 
In particular, the first case in (\ref{eq:p-mutation}) yields 
$\overline p_{\overline z}^\pm = p_z^\mp$. 
 
It is easy to check that the mutation of $\overline\Sigma$ in direction $\overline z$ recovers~$\Sigma$. 
 
\begin{definition} 
\label{def:cluster-algebra} 
\emph{(Normalized cluster algebra)} 
Let $\mathcal{S}$ be a set of seeds in $\mathcal{F}$ with 
the following properties: 
\begin{itemize} 
\item 
every seed $\Sigma\in\mathcal{S}$ admits mutations in all $n$ 
conceivable directions, and the results of all these mutations belong 
to~$\mathcal{S}$; 
\item 
any two seeds in $\mathcal{S}$ are obtained from each other by a 
sequence of mutations. 
\end{itemize} 
The sets $\xx$, for $\Sigma=(\xx, \pp, B)\in\mathcal{S}$, 
are called~\emph{clusters;} their elements are the~\emph{cluster 
  variables}; 
the set of all cluster variables is denoted by~$\mathcal{X}$. 
The set of all elements $p^\pm_x\!\in\! \pp$, for all seeds 
$\Sigma=(\xx, \pp, B)\!\in\!\mathcal{S}$, is denoted 
by~$\mathcal{P}$. 
The \emph{ground ring} $\ZZ[\mathcal{P}]$ is the subring of $\FFcal$ 
generated by~$\mathcal{P}$. 
The (normalized) \emph{cluster algebra} 
$\AA=\AA(\mathcal{S})$ 
is the $\ZZ[\mathcal{P}]$-subalgebra of $\FFcal$ generated 
by~$\mathcal{X}$. 
The \emph{exchange graph} of $\AA(\mathcal{S})$ 
is the $n$-regular graph whose vertices are labeled by the seeds 
in~$\mathcal{S}$, and whose edges correspond to mutations. (This 
is easily seen to be equivalent 
to~\cite[Definition~7.4]{fz-clust1}.) 
\end{definition} 
 
Definition~\ref{def:cluster-algebra} is a bit more restrictive than the one given 
in~\cite{fz-clust1}, where we allowed to use any subring with unit in 
$\ZZ\PP$ containing $\mathcal{P}$ as a ground~ring. 
Some concrete examples of cluster algebras will be given in 
Section~\ref{sec:ca-type-A} below.

\begin{remark} 
\label{rem:sign-symmetry} 
There is an involution $\Sigma \mapsto \Sigma^\vee$ on the set of 
seeds in $\FFcal$ acting by $(\xx,\pp,B) \mapsto (\xx, \pp^\vee,-B)$, 
where $(p^\vee_x)^\pm = p_x^\mp$. 
An easy check show that this involution commutes with seed 
mutations. 
Therefore, if a collection of seeds $\mathcal{S}$ 
satisfies the conditions in Definition~\ref{def:cluster-algebra}, 
then so does the collection $\mathcal{S}^\vee$. 
The corresponding cluster algebras $\AA(\mathcal{S})$ and 
$\AA(\mathcal{S}^\vee)$ are canonically identified with each other. 
(This is a reformulation of ~\cite[(2.8)]{fz-clust1}.) 
\end{remark}

Two cluster algebras $\AA(\mathcal{S})\subset\FFcal$ and 
$\AA(\mathcal{S}')\subset\FFcal'$ over the same semifield $\PP$ 
are called \emph{strongly isomorphic} 
if there exists a 
$\ZZ \PP$-algebra 
isomorphism $\FFcal\to\FFcal'$ 
that transports some (equivalently, any) seed in $\mathcal{S}$ into a 
seed  in~$\mathcal{S}'$, thus inducing a bijection 
$\mathcal{S}\to\mathcal{S}'$ 
and an algebra isomorphism $\AA(\mathcal{S})\to\AA(\mathcal{S}')$. 

The set of seeds $\mathcal{S}$ for a cluster algebra 
$\AA=\AA(\mathcal{S})$ (hence the algebra itself) 
is uniquely determined by any single seed 
$\Sigma= (\xx, \pp, B)\in\mathcal{S}$. 
Thus, $\AA$ is determined by $B$ and $\pp$ up to a strong 
isomorphism, justifying the notation $\AA=\AA(B,\pp)$. 
In general, an $n\times n$ matrix $B$ and a $2n$-tuple $\pp$ 
satisfying the normalization conditions define a cluster algebra 
$\AA(B,\pp)$ if and only if any matrix obtained from $B$ by a sequence 
of mutations is sign-skew-symmetric. 
This condition is in particular satisfied whenever $B$ is 
\emph{skew-symmetrizable} \cite[Definition~4.4]{fz-clust1}, i.e., 
there exists a diagonal matrix $D$ 
with positive diagonal entries such that $DB$ is skew-symmetric. 
Indeed, matrix mutations preserve skew-symmetrizability 
\cite[Proposition~4.5]{fz-clust1}, 
and any skew-symmetrizable matrix is sign-skew-symmetric.

Every cluster algebra over a fixed semifield~$\PP$ belongs to a 
\emph{series} $\AA(B,-)$ consisting of all 
cluster algebras of the form $\AA(B,\pp)$, where $B$ is fixed, and 
$\pp$ is allowed to vary. 
Two series $\AA(B,-)$ and $\AA(B',-)$ are 
\emph{strongly isomorphic} if there is a bijection sending each 
cluster algebra $\AA(B,\pp)$ to a strongly isomorphic cluster 
algebra $\AA(B',\pp')$. (This amounts to requiring that $B$ and 
$B'$ can be obtained from each other by a sequence of matrix 
mutations, modulo simultaneous relabeling of rows and columns.) 
 
\subsection{Finite type classification} 
\label{sec:finite-type-classification} 
A cluster algebra $\AA(\mathcal{S})$ is said to be of \emph{finite 
type} if the set of seeds $\mathcal{S}$ 
is finite.

Let  $B=(b_{ij})$  be an integer square matrix. 
Its \emph{Cartan counterpart} is a generalized Cartan matrix $A = 
A(B)=(a_{ij})$ of the same size defined by 
\begin{equation} 
\label{eq:assoc-cartan} a_{ij} = 
\begin{cases} 
2 & \text{if $i=j$;} \\ 
- |b_{ij}| & \text{if $i\neq j$.} 
\end{cases} 
\end{equation} 
 
The following classification theorem is our main result. 
 
\begin{theorem} 
\label{th:fin-type-class} 
All cluster algebras in any series $\AA(B,-)$ 
are simultaneously of finite or infinite type. 
There is a canonical bijection between the Cartan matrices of finite 
type and the strong isomorphism classes of 
series of cluster algebras of finite type. 
Under this bijection, a Cartan matrix $A$ of finite type 
corresponds to the series $\AA(B,-)$, where $B$ is an arbitrary 
sign-skew-symmetric matrix with $A(B)=A$. 
\end{theorem} 
 
We note that in the last claim of Theorem~\ref{th:fin-type-class}, 
the series $\AA(B,-)$ is well defined since 
$A$ is symmetrizable and therefore $B$ must be skew-symmetrizable.

By Theorem~\ref{th:fin-type-class}, each cluster 
algebra of finite type has a well-defined \emph{type} (e.g., 
$A_n,B_n,\dots$), mirroring the Cartan-Killing classification. 
 
 
We prove Theorem~\ref{th:fin-type-class} by splitting it into the 
following three statements 
(Theorems~\ref{th:finite-type-class-new}--\ref{th:series-CK} below). 
 
\begin{theorem} 
\label{th:finite-type-class-new} 
Suppose that 
\begin{align} 
\label{eq:data-A} &\text{$A$~is a Cartan matrix of finite type;} 
\\ 
\label{eq:data-B} &\text{$B_\circ=(b_{ij})$~is a sign-skew-symmetric matrix 
such that $A=A(B_\circ)$ and} 
\\ 
\nonumber &\text{$b_{ij} b_{ik} \geq 0$ for all $i,j,k$;} 
\\ \label{eq:data-pp} &\text{$\pp_\circ$~is  a 
$2n$-tuple of elements in $\PP$ satisfying the normalization 
conditions.} 
\end{align} 
Then 
$\AA(B_\circ,\pp_\circ)$ is a 
cluster algebra of finite type. 
\end{theorem} 
 
It is easy to see that for any Cartan matrix $A$ of finite type, 
there is a matrix $B_\circ$ satisfying (\ref{eq:data-B}). 
Indeed, the sign-skew-symmetric matrices $B$ with $A(B) = A$ are 
in a bijection with orientations of the Coxeter graph of $A$ 
(recall that this graph has $I$ as the set of vertices, with 
$i$ and $j$ joined by an edge whenever $a_{ij} \neq 0$): 
under this bijection, $b_{ij} > 0$ if and only if the edge 
$\{i,j\}$ is oriented from $i$ to $j$. 
Condition (\ref{eq:data-B}) means that $B_\circ$ corresponds to an 
orientation such that every vertex is a source or a sink; since 
the Coxeter graph is a tree, hence a bipartite graph, such an 
orientation exists. 
 
\begin{theorem} 
\label{th:finite-type-class-exhaust} Any cluster algebra $\AA$ of 
finite type is strongly isomorphic to a cluster algebra 
$\AA(B_\circ,\pp_\circ)$ for some data of the 
form~{\rm(\ref{eq:data-A})--(\ref{eq:data-pp})}. 
\end{theorem} 
 
\begin{theorem} 
\label{th:series-CK} 
Let $B$ and $B'$ be sign-skew-symmetric matrices such that 
$A(B)$ and $A(B')$ are Cartan matrices of finite type. 
Then the series $\AA(B,-)$ and $\AA(B',-)$ 
are strongly isomorphic if and only if $A(B)$ and $A(B')$ 
are of the same Cartan-Killing 
type. 
\end{theorem} 
 
In the process of proving these theorems, 
we obtain the following characterizations of the 
cluster algebras of finite type. 
 
\pagebreak[3] 
 
\begin{theorem} 
\label{th:fin-type-characterizations} 
For a cluster algebra 
$\mathcal{A}$, the following are equivalent: 
\begin{itemize} 
\item[{\rm(i)}] 
$\mathcal{A}$ is of finite type; 
\item[{\rm(ii)}] 
the set $\mathcal{X}$ of all cluster variables is finite; 
\item[{\rm(iii)}] 
for every seed  $(\xx, \pp, B)$ in~$\mathcal{A}$, the entries of 
the matrix $B=(b_{xy})$ satisfy the inequalities 
$|b_{xy}b_{yx}|\leq 3$, for all $x,y\in\xx$. 
\item[{\rm(iv)}] 
$\AA = \AA(B_\circ,\pp_\circ)$ for some data of the 
form~{\rm(\ref{eq:data-A})--(\ref{eq:data-pp}).} 
\end{itemize} 
\end{theorem} 
 
The equivalence ${\rm(i)}\!\Longleftrightarrow\!{\rm(iv)}$ 
in Theorem~\ref{th:fin-type-characterizations} is tantamount to 
Theorems~\ref{th:finite-type-class-new}--\ref{th:finite-type-class-exhaust}. 
 
\subsection{Cluster variables in the finite type} 
\label{sec:laurent-finite-type} 
 
The techniques  in our proof of 
Theorem~\ref{th:finite-type-class-new} allow us to 
enunciate the basic properties of cluster algebras of finite type. 
We begin by providing an explicit description of the set of 
cluster variables in terms of the corresponding root system. 
 
For the remainder of Section~\ref{sec:intro-main-results}, 
$A=(a_{ij})_{i,j\in I}$ is a Cartan matrix of finite type 
and $\AA=\AA(B_\circ,\pp_\circ)$ a cluster algebra (of finite type) 
related to $A$ as in Theorem~\ref{th:finite-type-class-new}. 
Let $\Phi$ be the root system associated with~$A$, 
with the set of simple roots $\Pi = \{\alpha_i : i \in I\}$ 
and the set of positive roots $\Phi_{>0}\,$. 
(Our convention on the correspondence 
between $A$ and $\Phi$ is that the simple reflections $s_i$ act on simple roots by 
$s_i (\alpha_j) = \alpha_j - a_{ij} \alpha_i\,$.) 
Let $\xx_\circ = \{x_i : i \in I\}$ be the cluster for the initial 
seed $(\xx_\circ, \pp_\circ, B_\circ)$. 
(By an abuse of notation, we label the rows and columns of $B_\circ$ by the 
elements of $I$ rather than by the variables $x_i\,$, for $i\in I$.) 
We will use the shorthand 
$x^\alpha = \prod_{i \in I} x_i^{a_i}$ 
for any vector $\alpha = \sum_{i \in I} a_i \alpha_i$ in the root lattice. 
 
The following result shows that the cluster variables of $\AA$ are naturally parameterized 
by the set $\Phi_{\geq -1} = \Phi_{> 0} \cup (- \Pi)$ of 
\emph{almost positive roots}. 
 
\begin{theorem} 
\label{th:cluster-variable-denominators-1} 
There is a unique bijection 
$\alpha\mapsto x[\alpha]$ 
between the almost positive roots in~$\Phi$ 
and the cluster variables in 
$\AA$ 
such that, for any $\alpha\in\Phi_{\geq -1}\,$, 
the cluster variable $x[\alpha]$ is expressed in terms of the initial 
cluster $\xx_\circ = \{x_i : i \in I\}$ as 
\begin{equation} 
\label{eq:laurent-finite-type} 
x[\alpha]=\frac{P_\alpha (\xx_\circ)}{x^\alpha} 
\,, 
\end{equation} 
where $P_\alpha$ is a polynomial over $\ZZ\PP$ with nonzero constant term. 
Under this bijection, $x[-\alpha_i]=x_i\,$. 
\end{theorem} 
 
Formula (\ref{eq:laurent-finite-type}) is an example of 
the \emph{Laurent phenomenon} established in~\cite{fz-clust1} for 
arbitrary cluster algebras: 
every cluster variable can be written as a Laurent polynomial in the 
variables of an arbitrary fixed cluster and the elements of~$\PP$. 
In~\cite{fz-clust1}, we conjectured that the coefficients of these Laurent 
polynomials are always nonnegative. 
Our next result establishes this 
conjecture (indeed, strengthens it) in the special case of the 
distinguished cluster $\xx_\circ$ in a cluster algebra of finite type. 
 
\begin{theorem} 
\label{th:cluster-variable-denominators-2} 
Every coefficient of each polynomial $P_\alpha$ 
(see {\rm (\ref{eq:laurent-finite-type})}) 
can be written as a polynomial in the elements of $\mathcal{P}$ 
(see Definition~\ref{def:cluster-algebra}) with positive integer 
coefficients. 
\end{theorem}

 
\subsection{Cluster complexes} 
\label{sec:cluster-complexes} 
We next focus on the combinatorics of clusters. 
As before, $\AA$ is a cluster algebra of finite type 
associated with a root system~$\Phi$. 
 
\begin{theorem} 
\label{th:exch-rel-dep-on-vars-only} 
The exact form of each exchange relation 
{\rm(\ref{eq:exchange-rel-xx})} in $\AA$ 
(that is, the cluster variables, exponents, and coefficients appearing in 
the right-hand side) 
\linebreak[3] 
depends only on the ordered pair $(z,\overline z)$ of cluster 
variables, and not on the particular choice of clusters (or seeds) containing them. 
\end{theorem} 
 
In fact, we do more: we describe  in concrete root-theoretic 
terms all pairs $(\beta, \beta')$ of almost positive 
roots such that the product $x[\beta] x[\beta']$ appears as a 
left-hand side of an exchange relation, and 
for every such pair, we describe the exponents 
appearing on the right. 
See 
Definition~\ref{def:b-matrices-via-roots} and formula (\ref{eq:exch-rel-fin-type}). 
 
\begin{theorem} 
\label{th:pseudomanifold-clusters} 
Every seed $(\xx,\pp,B)$ in $\AA$ 
is uniquely determined by its cluster~$\xx$. 
For any cluster $\xx$ and any 
$x\in\xx$, 
there is a unique cluster $\xx'$ with $\xx\cap\xx'=\xx-\{x\}$. 
\end{theorem} 
 
We conjecture that the requirement of finite type 
in Theorems~\ref{th:exch-rel-dep-on-vars-only} 
and~ \ref{th:pseudomanifold-clusters} 
can be dropped; that is, \emph{any} cluster algebra conjecturally has 
these properties. 
 
We define the \emph{cluster complex} $\Delta(\AA)$  as the simplicial 
complex whose ground set is $\mathcal{X}$ 
(the set of all cluster variables) and whose maximal 
simplices are the clusters. 
By Theorem~\ref{th:pseudomanifold-clusters}, 
the cluster complex encodes the combinatorics of seed mutations. 
Thus, the dual graph of~$\Delta(\AA)$ is precisely the exchange 
graph of~$\AA$. 
 
Our next result identifies the cluster complex 
$\Delta(\AA)$ with the dual complex $\Delta(\Phi)$ of the 
\emph{generalized associahedron} of the corresponding type. 
The simplicial complexes $\Delta(\Phi)$ were introduced  and 
studied in \cite{yga}; see also \cite{cfz} and 
Section~\ref{sec:root-systems-gen-associahedra}  below.

\begin{theorem} 
\label{th:finite-type-complex} 
Under the bijection 
\hbox{$\Phi_{\geq -1}\to\mathcal{X}$} 
of Theorem~\ref{th:cluster-variable-denominators-1}, 
the cluster complex $\Delta(\AA)$ is  identified 
with the 
simplicial complex~$\Delta(\Phi)$. 
In particular, the cluster complex does not depend on the coefficient 
semifield $\PP$, 
or on the choice of the coefficients $\pp_\circ$ in the initial seed. 
\end{theorem}

\subsection{Organization of the paper} 
\label{sec:organization-ca2} 
 
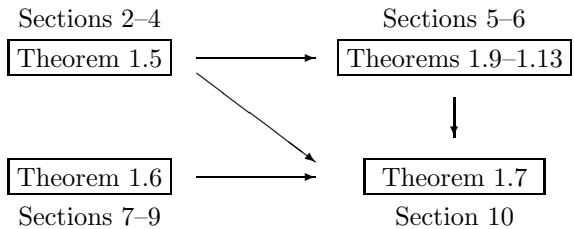
\begin{figure}[ht] 
\begin{center} 
\setlength{\unitlength}{1.5pt} 
\begin{picture}(0,50)(48,2) 
 
\put(3,10){\makebox(0,0){$\boxed{\text{Theorem~\ref{th:finite-type-class-exhaust}}}$}} 
\put(3,40){\makebox(0,0){$\boxed{\text{Theorem~\ref{th:finite-type-class-new}}}$}} 
 
\put(3,0){\makebox(0,0){${\text{Sections~\ref{sec:2-finite}--\ref{sec:diagram-CK-proof}}}$}} 
 
\put(3,50){\makebox(0,0){${\text 
{Sections~\ref{sec:cluster-algebras-via-pseudomanifolds}--\ref{sec:construction-of-finite-type}}}$}}

\put(30,10){\vector(1,0){30}} 
 
\put(30,40){\vector(1,0){30}} 
 
\put(30,36.25){\vector(4,-3){30}} 
 
\put(95,10){\makebox(0,0){$\boxed{\,\,\text{Theorem 
\ref{th:series-CK}}\,\,}$}} 
 
\put(95,40){\makebox(0,0){$\boxed{\text{Theorems~\ref{th:cluster-variable-denominators-1}--\ref{th:finite-type-complex}}}$}} 
 
\put(95,0){\makebox(0,0){$\text{Section~\ref{sec:proof-CK-exactly}}$}} 
 
\put(95,50){\makebox(0,0){$\text{Sections~\ref{sec:proofs-laurent-etc}--\ref{sec:initial-cluster-positivity}}$}} 
 
\put(95,30){\vector(0,-1){10}} 
 
\end{picture} 
\end{center} 
\caption{\hbox{Logical dependences among the proofs of 
Theorems~\ref{th:finite-type-class-new}--\ref{th:finite-type-complex}}} 
\label{fig:flowchart} 
\end{figure}

The bulk of the paper is devoted to the proofs of 
Theorems~\ref{th:fin-type-class}--\ref{th:fin-type-characterizations}. 
We already noted that Theorem~\ref{th:fin-type-class} 
follows from 
Theorems~\ref{th:finite-type-class-new}--\ref{th:series-CK}, 
and that the implications ${\rm(iv)}\!\Longrightarrow\!{\rm(i)}$ 
and ${\rm(i)}\!\Longrightarrow\!{\rm(iv)}$ 
in Theorem~\ref{th:fin-type-characterizations} are~essentially 
Theorems~\ref{th:finite-type-class-new} 
and~\ref{th:finite-type-class-exhaust}, respectively. 
Furthermore, 
${\rm(i)}\!\Longrightarrow\!{\rm(ii)}$ is trivial, 
while ${\rm(ii)}\!\Longrightarrow\!{\rm(iii)}$ follows from 
\cite[Theorem~6.1]{fz-clust1}. 
Thus, we need to prove the following: 
\begin{itemize} 
\item 
Theorem~\ref{th:finite-type-class-new}; 
\item 
Theorem~\ref{th:finite-type-class-exhaust} via the implication 
${\rm(iii)}\!\Longrightarrow\!{\rm(iv)}$ 
in Theorem~\ref{th:fin-type-characterizations}; 
\item 
Theorem~\ref{th:series-CK}. 
\end{itemize} 
Figure~\ref{fig:flowchart} shows the logical dependences between these 
proofs, and the sections containing them. 
Theorems~\ref{th:cluster-variable-denominators-1}--\ref{th:finite-type-complex}, 
which only rely on Theorem~\ref{th:finite-type-class-new}, 
are proved in Sections~\ref{sec:proofs-laurent-etc}-- 
\ref{sec:initial-cluster-positivity}, 
following the completion of the proof of 
Theorem~\ref{th:finite-type-class-new}.

The concluding Section~\ref{sec:geometric-realization} 
provides explicit geometric realizations for some special 
cluster algebras of the classical types~$ABCD$. 
These examples are based on a general criterion given 
in Section~\ref{sec:geometric-type-ca} 
for a cluster algebra to be isomorphic to a $\ZZ$-form 
of the coordinate ring of some algebraic variety. 
In particular, we show that a $\ZZ$-form of the homogeneous coordinate 
ring of the Grassmannian ${\rm Gr}_{2,m}$ ($m\geq 5$) 
in its Pl\"ucker embedding carries two different cluster algebra 
structures of types $A_{m-3}$ and $B_{m-2}\,$, respectively. 
 
\pagebreak[2] 
 
\section{Cluster algebras via pseudomanifolds} 
\label{sec:cluster-algebras-via-pseudomanifolds} 
 
\subsection{Pseudomanifolds and geodesic loops} 
\label{sec:pseudomanifolds} 
 
This section begins our proof of 
Theorem~\ref{th:finite-type-class-new}. 
Its main result (Proposition~\ref{prop:finite-type-criterion}) 
provides sufficient conditions 
ensuring that a cluster algebra that 
arises from a particular combinatorial construction 
is of finite type. 
 
The first ingredient of this construction is an $(n-1)$-dimensional 
pure simplicial complex $\Delta$ (finite or infinite) on the ground set~$\Psi$. 
Thus, every maximal simplex in $\Delta$ is an $n$-element subset of~$\Psi$ 
(a~``cluster''). 
A simplex of codimension $1$ (i.e., an $(n-1)$-element subset 
of~$\Psi$) is called a \emph{wall}. 
The vertices of the \emph{dual graph} $\Gamma$ are the clusters 
in~$\Delta$; 
two clusters are connected by an edge in $\Gamma$ if they share a wall. 
 
We assume that $\Delta$ is a \emph{pseudomanifold}, 
i.e., 
\begin{align} 
\label{eq:pseudomanifold-1} 
&\text{every wall is contained in precisely two 
maximal simplices (clusters);} \hspace{.3in} 
\\ 
\label{eq:pseudomanifold-2} 
&\text{the dual graph $\Gamma$ is connected.} 
\end{align} 
In view of (\ref{eq:pseudomanifold-1}), the graph $\Gamma$ is 
$n$-regular, i.e., there are precisely $n$ 
edges incident to every vertex~$C\in\Gamma$. 
 
\begin{example} 
Let $n=1$. 
Then (\ref{eq:pseudomanifold-1}) is saying that 
the empty simplex is contained in precisely two $0$-dimensional 
simplices (points in~$\Psi$). 
Thus a $0$-dimensional pseudomanifold must be a disjoint union of two points 
(a $0$-dimensional sphere). 
The dual graph $\Gamma$ has these points as vertices, 
with an edge connecting them. 
 
For $n=2$, a $1$-dimensional pseudomanifold is nothing but a 
$2$-regular connected graph---thus, either an infinite chain or a 
cycle. Ditto for its dual graph. 
\end{example} 
 
For a non-maximal simplex $D \in \Delta$, we denote by $\Delta_D$ the 
\emph{link} of~$D$. 
This is the simplicial complex on the ground set 
$\Psi_D = \{\alpha \in \Psi - D \, : \, D \cup \{\alpha\} \in \Delta\}$ 
such that $D'$ is a simplex in $\Delta_D$  if and only if 
$D \cup D'$ is a simplex in~$\Delta$. 
The link $\Delta_D$ is a pure simplicial complex of dimension 
$n-|D|-1$ satisfying property (\ref{eq:pseudomanifold-1}) 
of a pseudomanifold. 
 
We will assume that $\Delta$ satisfies the following additional condition: 
\begin{align} 
\label{eq:links-connected} 
&\text{the link of every non-maximal simplex $D$ in $\Delta$ is a 
pseudomanifold.} \hspace{.3in} 
\end{align} 
Equivalently, the dual graph $\Gamma_D$ of $\Delta_D$ is 
connected. 
 
Conditions (\ref{eq:pseudomanifold-1})--(\ref{eq:links-connected}) 
can be restated as saying that $\Delta$ is a (possibly infinite, simplicial) 
\emph{abstract polytope} in the sense of \cite{adler} or 
\cite{mcmullen-schulte}; 
another terminology is that $\Delta$ is a \emph{thin, residually connected} 
complex (see, e.g., \cite{aschbacher}). 
 
We identify 
the graph $\Gamma_D$ with an induced subgraph in $\Gamma$ whose 
vertices are the maximal simplices in $\Delta$ that contain~$D$. 
In particular, for $|D| = n-2$, the pseudomanifold $\Delta_D$ is $1$-dimensional, so 
$\Gamma_D$ is either an infinite chain or a finite cycle in~$\Gamma$. 
In the latter case, we call $\Gamma_D$  a \emph{geodesic loop}. 
(This is a \emph{geodesic} in~$\Gamma$ with respect to the canonical 
\emph{connection} on~$\Gamma$, in the sense of~\cite{bgh, guillemin-zara}.) 

We assume that 
\begin{align} 
\label{eq:geodesic-loops-generate} 
&\text{the fundamental group of $\Gamma$ is generated by the 
geodesic loops pinned\quad}\\ 
\nonumber 
&\text{down to a fixed basepoint.} 
\end{align} 
More precisely, by (\ref{eq:geodesic-loops-generate}) we mean that the 
fundamental group of $\Gamma$ is generated by all the 
loops of the form $PL\bar P$, where $L$ is a geodesic loop, 
$P$ is a path originating at the basepoint, 
and $\bar P$ is the inverse path to~$P$. 
 
 
\begin{lemma} 
\label{lem:triangulated-polytope} 
Let $\Delta$ be the boundary complex of an $n$-dimensional simplicial 
convex polytope. 
Then conditions 
{\rm (\ref{eq:pseudomanifold-1})--(\ref{eq:geodesic-loops-generate})} are 
satisfied. 
\end{lemma} 
 
Equivalently, conditions 
(\ref{eq:pseudomanifold-1})--(\ref{eq:geodesic-loops-generate}) hold 
if $\Delta$ is (the simplicial complex of) 
the normal fan of a simple $n$-dimensional convex polytope~$\Delta^*$. 
 
\begin{proof} 
Statements (\ref{eq:pseudomanifold-1})--(\ref{eq:pseudomanifold-2}) 
are trivial. 
The link $\Delta_D$ of each non-maximal face~$D$ in $\Delta$ is 
again a simplicial polytope, 
implying~(\ref{eq:links-connected}). 
Specifically, $\Delta_D$ is canonically identified (see, e.g., 
\cite[Problem~VI.1.4.4]{barvinok-book}) with the dual polytope for the dual 
face $D^*$ in the dual simple polytope~$\Delta^*$. 
Thus, the graph $\Gamma_D$ is the $1$-skeleton of~$D^*$. 
 
It remains to check (\ref{eq:geodesic-loops-generate}). 
The case $n=2$ is trivial, so let us assume that $n\geq 3$. 
Each geodesic loop $\Gamma_D$ 
(for an $(n-2)$-dimensional face~$D$) 
is identified with the $1$-skeleton (i.e, the boundary) 
of the dual $2$-dimensional face~$D^*$ in~$\Delta^*$. 
The boundary cell complex of~$\Delta^*$ is spherical, hence simply 
connected. 
On the other hand, the fundamental group of this complex 
can be obtained as a quotient of the fundamental group of its 
$1$-skeleton~$\Gamma$ by the normal subgroup generated by the 
boundaries of its $2$-dimensional cells 
pinned down to a basepoint 
(see, e.g., \cite[Theorems~VII.2.1, VII.4.1]{massey}), 
or, equivalently, by the subgroup generated by all pinned-down 
geodesic loops. This proves (\ref{eq:geodesic-loops-generate}). 
\end{proof} 
 
\subsection{Sufficient conditions for finite type} 
 
We next describe the second ingredient of our construction. 
Suppose that we have a family of integer matrices 
$B(C)=(b_{\alpha\beta}(C))_{\alpha, \beta\in C}$, for all vertices $C$ 
in~$\Gamma$, 
satisfying the following conditions: 
\begin{align} 
\label{eq:some-matrix-is-ss} 
&\text{all the matrices $B(C)$ are sign-skew-symmetric.} \hspace{0.7in} 
\\ 
\label{eq:B-mutation-gamma} 
&\text{for every edge $(C,\overline C )$ in $\Gamma$, 
with $\overline C =C-\{\gamma\}\cup\{\overline\gamma\}$, the matrix 
$B(\overline C)$ 
} 
\\ 
\nonumber 
&\text{is obtained from $B(C)$ by a matrix mutation in direction $\gamma$ 
followed by} 
\\ 
\nonumber 
&\text{relabeling one row and one column 
by replacing $\gamma$ by~$\overline \gamma$.} 
%
\end{align} 
 
We 
need one more assumption concerning the matrices $B(C)$, which 
will require a little preparation. Fix a geodesic $\Gamma_D$ and 
associate to its every vertex $C = D \cup \{\alpha, \beta\}$ the 
integer $b_{\alpha\beta}(C) b_{\beta\alpha}(C)$. It is trivial to 
check, using (\ref{eq:B-mutation-gamma}), that this integer 
depends only on $D$, not on the particular choice of $\alpha$ and 
$\beta$. We say that $\Gamma_D$ is of \emph{finite type} if 
$b_{\alpha\beta}(C) b_{\beta\alpha}(C) \in \{0,-1,-2,-3\}$ for 
some (equivalently, any) vertex $C = D \cup \{\alpha, \beta\}$ 
on~$\Gamma_D\,$. 
If this is the case, then we associate to $\Gamma_D$ the \emph{Coxeter number} 
$h\in\ZZ_{>0}$ defined by 
\[ 
2 \cos (\pi/h) = 
\sqrt{|b_{\alpha\beta}(C) b_{\beta\alpha}(C)|} 
\,, 
\] 
or, equivalently, by the table 
\begin{center} 
\begin{tabular}{|c|c|c|c|c|} 
\hline 
$|b_{\alpha\beta}(C) b_{\beta\alpha}(C)|$ & 0 & 1 & 2 & 3 \\ 
\hline 
$h(C,x,y) $ & 2 & 3 & 4 & 6 \\ 
\hline 
\end{tabular} \,. 
\end{center} 
 
Our last condition is: 
\begin{align} 
\label{eq:geodesic-loops} 
&\text{every geodesic loop in $\Gamma$ is of finite type, 
and has length 
$h+ 2$, where $h$ is} 
\\ 
\nonumber 
&\text{the corresponding Coxeter number.} 
\end{align} 
 
 
\begin{proposition} 
\label{prop:finite-type-criterion} 
Assume that a simplicial complex $\Delta$ and a family of matrices 
$(B(C))$  satisfy the assumptions 
{\rm(\ref{eq:pseudomanifold-1})--(\ref{eq:geodesic-loops})} above. 
Let $B=B(C)$ for some vertex~$C$, and let $\AA = \AA(B,\pp)$ be the cluster algebra associated 
with $B$ 
and some coefficient tuple $\pp$. There exists a surjection from 
the set 
of vertices of $\Gamma$ onto the set of all seeds for $\AA$. In 
particular, if $\Delta$ is finite, then $\AA$ is of finite type. 
\end{proposition} 
 
We will prove this proposition by showing that, 
 whether $\Delta$ is finite or infinite, 
its dual graph $\Gamma$ is always a covering graph for the exchange 
 graph of~$\AA(B, \pp)$. 
To~formulate this more precisely, we need some preparation. 
 
Let $C$ be a vertex of $\Gamma$. 
A \emph{seed attachment} at $C$ consists of a choice of a seed 
$\Sigma = (\xx, \pp, B)$ and a bijection 
$\alpha \mapsto x[C,\alpha]$ between 
$C$ and $\xx$ identifying the matrices $B(C)$ and~$B$, so that 
$b_{\alpha\beta}(C) = b_{x[C,\alpha],x[C,\beta]}$. 
The \emph{transport} of a seed attachment along an edge 
$(C,\overline C)$ with $\overline 
C=C-\{\gamma\}\cup\{\overline\gamma\}$ is defined as follows: 
the seed $\overline \Sigma = (\overline\xx,\overline \pp,\overline B)$ 
attached to $\overline C$ 
is obtained from $\Sigma$ by the mutation in direction~$x[C,\gamma]$, 
and the corresponding bijection $\overline C\to \overline\xx$ is 
uniquely determined by 
$x[\overline C,\alpha]=x[C,\alpha]$ for all $\alpha\in 
C\cap\overline C$. 
(The remaining cluster variable $x[\overline C,\overline\gamma]$ 
is obtained from $x[C,\gamma]$ by the corresponding exchange 
relation~(\ref{eq:exchange-rel-xx}).) 
We note that transporting the resulting seed attachment backwards from 
$\overline C$ to $C$ recovers the original seed attachment. 
 
Proposition~\ref{prop:finite-type-criterion} 
is an immediate consequence of the following lemma. 
 
\begin{lemma} 
\label{lem:attachments} 
Let $\Delta$ 
and $(B(C))$ satisfy 
{\rm(\ref{eq:pseudomanifold-1})--(\ref{eq:geodesic-loops})}. 
Suppose we are given a vertex $C_\circ$ in~$\Gamma$ together with a 
seed attachment involving a seed $\Sigma_\circ$ for a cluster algebra~$\AA$. 
 
\begin{enumerate} 
\item 
The given seed attachment at $C_\circ$ extends uniquely to a family of 
seed attachments at all vertices in $\Gamma$ such that, 
for every edge $(C,\overline C)$, the seed attachment at 
$\overline C$ is obtained from that at $C$ by transport along~$(C,\overline C)$. 
 
\item 
Let $\Sigma(C)$ denote the seed attached to a vertex~$C$. 
The map $C\mapsto\Sigma(C)$ is a surjection onto the set of all seeds 
for~$\AA$. 
 
\item 
For every vertex $C$ and every $\alpha\in C$, 
the cluster variable $x[C,\alpha]$ attached to $\alpha$ at $C$ 
depends only on~$\alpha$ (so can be denoted by~$x[\alpha]$). 
\item 
The map $\alpha\mapsto x[\alpha]$ is a surjection from the ground set 
$\Psi$ onto the set of all cluster variables for~$\AA$. 
\end{enumerate} 
\end{lemma} 
 
\begin{proof} 
1. Since $\Gamma$ is connected, we can transport the initial seed 
attachment at $C_\circ$ to an arbitrary vertex $C$ along a path 
from $C_\circ$ to~$C$. 
We need to show that the resulting seed attachment 
at~$C$ is independent of the choice of a path. 
For that, it suffices to prove that transporting a 
seed attachment along a loop in $\Gamma$ brings it back unchanged. 
By (\ref{eq:geodesic-loops-generate}), it is enough to show this for the geodesic 
loops. Then the claim follows from (\ref{eq:geodesic-loops}) 
and \cite[Theorem~7.7]{fz-clust1}. 
 
2. Take an arbitrary seed $\Sigma$ for~$\AA$. 
By Definition~\ref{def:cluster-algebra}, $\Sigma$ can be obtained from 
the initial seed $\Sigma_\circ$ by a sequence of mutations. 
This sequence is uniquely lifted to a path $(C_\circ,\dots,C)$ in $\Gamma$ 
such that transporting the initial seed attachment at $C_\circ$ along 
the edges of this path produces the chosen sequence of mutations. 
Hence $\Sigma(C)=\Sigma$, as desired. 
 
3. Let $\alpha\in\Psi$, and let $C$ and $C'$ be two vertices of 
$\Gamma$ such that $\alpha\in C\cap C'$. 
By (\ref{eq:links-connected}), $C$ and $C'$ can be joined by a path 
$(C_1=C,C_2,\dots,C_\ell=C')$ such that $\alpha\in C_i$ for all~$i$. 
Hence $x[C_1,\alpha]=x[C_2,\alpha]=\cdots=x[C_\ell,\alpha]$, as needed. 
 
4. Follows from Part~2. 
\end{proof} 
 
\begin{remark} 
\label{rem:isom-cluster-complex} 
Parts 1 and 2 in Lemma~\ref{lem:attachments} imply that the map 
$C\mapsto \Sigma(C)$ induces a covering of the exchange graph 
of~$\AA$ by the graph~$\Gamma$. 
If, in addition, the map $\alpha\mapsto x[\alpha]$ in 
Lemma~\ref{lem:attachments} is a 
bijection, then the map $C\mapsto \Sigma(C)$ is also a bijection. 
Thus, the latter map establishes an isomorphism between $\Gamma$ and 
the exchange graph of~$\AA$, 
and between $\Delta$ and the cluster complex of~$\AA$. 
\end{remark} 
 
 
\section{Generalized associahedra} 
\label{sec:root-systems-gen-associahedra} 
 
This section contains an exposition of the results in \cite{yga} 
and~\cite{cfz} that will be used later in our proof of 
Theorem~\ref{th:finite-type-class-new}. 
 
 
Let $A = (a_{ij})_{i,j \in I}$ be an indecomposable Cartan matrix 
of finite type, and $\Phi$ the corresponding irreducible root system 
of rank $n = |I|$. 
We retain the notation introduced in Section~\ref{sec:laurent-finite-type}. 
In particular, $\Phi_{\geq -1}=\Phi_{>0}\cup(-\Pi)$ denotes the set of almost positive 
roots. 

The Coxeter graph associated to $\Phi$ is a tree; 
recall that this graph has the index set $I$ as the set of vertices, with $i$ 
and $j$ joined by an edge whenever $a_{ij} \neq 0$. 
In particular, the Coxeter graph is bipartite; the two parts 
$I_+,I_-\subset I$ are determined uniquely up to renaming. 
The \emph{sign function} $\varepsilon:I\to\{+,-\}$ is defined by 
\begin{equation} 
\label{eq:sign-function} 
\varepsilon(i)=\begin{cases} 
+ & \text{if $i\in I_+\,$;}\\ 
- & \text{if $i\in I_-\,$.} 
\end{cases} 
\end{equation} 
 
Let $Q = \ZZ \Pi$ denote the root lattice, 
and $Q_\RR$ the ambient real vector space. 
For $\gamma \!\in\! Q_\RR$, we denote by $[\gamma : \alpha_i]$ 
the coefficient of $\alpha_i$ in the expansion of $\gamma$ in the 
basis~$\Pi$. 
Let $\tau_+$ and $\tau_-$ denote the piecewise-linear 
automorphisms of $Q_\RR$ given~by 
\begin{equation} 
\label{eq:tau-pm-tropical} 
[\tau_\varepsilon  \gamma: \alpha_i] = 
\begin{cases} 
- [\gamma: \alpha_i] - \sum_{j \neq i}  a_{ij} \max ([\gamma: \alpha_j], 0) 
& 
\text{if $i \in I_\varepsilon$;} \\{} 
[\gamma: \alpha_i] & \text{otherwise.} 
\end{cases} 
\end{equation} 
It is easy to see that each of $\tau_+$ and $\tau_-$ is an involution 
that preserves the set $\Phi_{\geq -1}$. 
More specifically, the action of $\tau_+$ and $\tau_-$ on $\Phi_{\geq -1}$ 
can be described as follows: 
\begin{equation} 
\label{eq:tau-pm-on-roots} 
\tau_\varepsilon(\alpha) = 
\begin{cases} 
\displaystyle 
\alpha & \text{if $\alpha = - \alpha_i, \, i \in I_{- \varepsilon}$;}\\ 
(\prod_{i \in I_\varepsilon} s_i) \,(\alpha) & \text{otherwise.} 
\end{cases} 
\end{equation} 
(The product of simple reflections $\prod_{i \in I_\varepsilon} s_i$ is well-defined since 
its factors commute). 
To illustrate, consider the type $A_2$, with $I_+=\{1\}$ and $I_-=\{2\}$. 
Then 
\begin{equation} 
\label{eq:A2-tau-tropical} 
\begin{array}{ccc} 
-\alpha_1 & \stackrel{\textstyle\tau_+}{\verylongleftrightarrow} 
~\alpha_1~ \stackrel{\textstyle\tau_-}{\verylongleftrightarrow} 
~\alpha_1\!+\!\alpha_2~ 
\stackrel{\textstyle\tau_+}{\verylongleftrightarrow} ~\alpha_2~ 
\stackrel{\textstyle\tau_-}{\verylongleftrightarrow} 
& -\alpha_2\,. \\
\circlearrowright & & \circlearrowright \\ \tau_- & & \tau_+ 
\end{array} 
\end{equation} 
We denote by $\langle\tau_-,\tau_+\rangle$  the group generated by 
$\tau_-$ and~$\tau_+\,$. 
 
The Weyl group of $\Phi$ is denoted by~$W$, its longest element 
by~$w_\circ\,$, and its Coxeter number by~$h$. 

\begin{theorem} 
\label{th:dihedral} 
{\rm \cite[Theorems~1.2,~2.6]{yga}} 
{\ } 
\begin{enumerate} 
\item 
The order of $\tau_-  \tau_+$ is 
equal to $(h+2)/2$ if $w_\circ = -1$, and to $h+2$ otherwise. 
Accordingly, $\langle\tau_-,\tau_+\rangle$ 
is a dihedral group of order 
$(h+2)$ or $2(h+2)$. 
\item 
The correspondence $\Omega \mapsto \Omega \cap (- \Pi)$ 
is a bijection between the $\langle\tau_-,\tau_+\rangle$-orbits in  $\Phi_{\geq - 1}$ 
and the $\langle - w_\circ \rangle$-orbits 
in~$(- \Pi)$. 
\end{enumerate} 
\end{theorem} 
 
We note that Theorem~\ref{th:dihedral} is stronger than 
\cite[Theorem~2.6]{yga}, since in the latter, $\tau_-$ 
and $\tau_+$ are treated as permutations of the set $\Phi_{\geq - 
  1}\,$, rather than as transformations of the entire 
space~$\QQ_\RR$. 
This stronger version follows from \cite[Theorem~1.2]{yga} by 
``tropical specialization'' (see \cite[(1.8)]{yga}). 
 
 
 
According to \cite[Section~3.1]{yga}, there exists a unique function 
$(\alpha,\beta)\mapsto (\alpha \| \beta)$ on $\Phi_{\geq 
  -1}\times\Phi_{\geq -1}$ with nonnegative integer values, 
called the \emph{compatibility degree,} 
such that 
\begin{equation} 
\label{eq:compatibility-1} 
(- \alpha_i \| \alpha) 
= \max \ ([\alpha : \alpha_i], 0) 
\end{equation} 
for any $i \in I$ and $\alpha \in\Phi_{\geq -1}\,$, 
and 
\begin{equation} 
\label{eq:compatibility-2} 
(\tau_\varepsilon \alpha \| \tau_\varepsilon \beta) = (\alpha \| 
\beta) 
\end{equation} 
for any $\alpha, \beta \in \Phi_{\geq -1}$ 
and any sign~$\varepsilon$. 
We say that 
$\alpha$ and $\beta$ are 
\emph{compatible} if 
$(\alpha \| \beta)=0$. 
(This is equivalent to $(\beta \| \alpha) = 0$ by 
\cite[Proposition~3.3, Part~2]{yga}.) 
 
Let~$\Delta (\Phi)$ be the simplicial complex on the ground set 
$\Phi_{\geq -1}$ whose simplices are the subsets of mutually compatible 
roots. 
As in Section~\ref{sec:cluster-algebras-via-pseudomanifolds} above, 
the maximal simplices of $\Delta (\Phi)$ are called \emph{clusters}. 
 
\begin{theorem} 
\label{th:cluster-fan} 
{\rm \cite[Theorems~1.8,~1.10]{yga} \cite[Theorem~1.4]{cfz}} 
{\ } 
\begin{enumerate} 
\item 
Each cluster in $\Delta(\Phi)$ is a $\ZZ$-basis of the root lattice~$Q$; 
in particular, all clusters are of the same size~$n$. 
\item 
The cones spanned by the simplices in $\Delta(\Phi)$ form a complete 
simplicial fan in~$Q_\RR$. 
\item 
This simplicial fan 
is the normal fan of a simple $n$-dimension\-al convex polytope, 
the \emph{generalized associahedron} of the corresponding type. 
\end{enumerate} 
\end{theorem} 
 
\begin{figure}[ht] 
\begin{center} 
\setlength{\unitlength}{1.5pt} 
\begin{picture}(120,78)(-60,-44) 
\thicklines 
 
\put(0,0){\circle{20}} 
 
\thinlines 
 
\put(0,0){\circle*{1}} 
 
\put(0,0){\vector(1,0){50}} 
\put(0,0){\vector(-1,0){50}} 
\put(0,0){\vector(-2,3){25}} 
\put(0,0){\vector(2,3){25}} 
\put(0,0){\vector(2,-3){25}} 
 
\put(58,0){\makebox(0,0){$\alpha_1$}} 
\put(-58,0){\makebox(0,0){$-\alpha_1$}} 
\put(38,37.5){\makebox(0,0){$\alpha_1\!+\!\alpha_2$}} 
\put(-32,37.5){\makebox(0,0){$\alpha_2$}} 
\put(34,-37.5){\makebox(0,0){$-\alpha_2$}} 
 
\thicklines 
 
\put(25,15){\line(0,-1){30}} 
\put(25,15){\line(-3,2){25}} 
\put(-25,15){\line(3,2){25}} 
\put(25,-15){\line(-3,-2){50}} 
\put(-25,15){\line(0,-1){63.33}} 
 
\put(10.5,0){\circle*{2}} 
\put(-10.5,0){\circle*{2}} 
\put(5.75,8.62){\circle*{2}} 
\put(-5.75,8.62){\circle*{2}} 
\put(5.75,-8.62){\circle*{2}} 
 
\put(25,15){\circle*{2}} 
\put(25,-15){\circle*{2}} 
\put(-25,15){\circle*{2}} 
\put(0,31.67){\circle*{2}} 
\put(-25,-48.33){\circle*{2}} 
 
\put(-5,-16){\makebox(0,0){$\Delta(\Phi)$}} 
 
\end{picture} 
\end{center} 
\caption{\hbox{The complex $\Delta(\Phi)$ and the corresponding polytope 
in type $A_2$}} 
\label{fig:5-roots-dual} 
\end{figure}
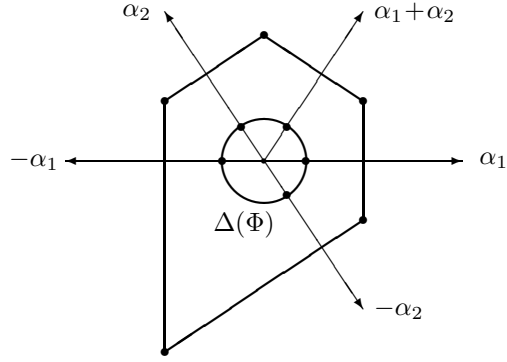 
 
Generalized associahedra of types $ABC$ are: 
in type~$A$, the Stasheff polytope, or ordinary \emph{associahedron} 
(see, e.g., \cite{stasheff,lee} or \cite[Chapter~7]{gkz}); 
in types~$B$ and~$C$, 
the Bott-Taubes polytope, or \emph{cyclohedron} 
(see \cite{bott-taubes, markl, simion-B}). 
Explicit combinatorial descriptions of generalized associahedra of 
types~$ABCD$ in relation to the root system framework 
are discussed in \cite{yga,cfz}; see also Section~\ref{sec:geometric-realization} below.

\begin{proposition} 
\label{prop:inv-cluster-expansion-1} 
{\rm \cite[Theorem~3.11]{yga}} 
Every vector $\gamma \in Q$ has a unique {\rm cluster expansion,} 
that is, $\gamma$ can be expressed uniquely 
as a nonnegative linear combination of mutually compatible roots 
from $\Phi_{\geq -1}$ 
(the {\rm cluster components} of~$\gamma$). 
\end{proposition} 
 
\begin{proposition} 
\label{prop:inv-cluster-expansion-2} 
{\rm \cite[Proposition~1.13]{cfz}} 
Let $[\gamma: \alpha]_{\rm clus}$ denote the coefficient of 
an almost positive root $\alpha$ in the cluster expansion of a 
vector $\gamma \in Q$. 
Then we have 
$[\sigma (\gamma): \sigma (\alpha)]_{\rm clus} = [\gamma: \alpha]_{\rm clus}$ 
for $\sigma\in\langle \tau_+, \tau_- \rangle$. 
\end{proposition} 
 
We call two roots $\beta,\beta'\in\Phi_{\geq -1}$ \emph{exchangeable} 
if $(\beta\|\beta')=(\beta'\|\beta) = 1$. 
The choice of terminology is motivated by the following proposition. 

\begin{proposition} 
\label{pr:two-terms-1} 
{\rm \cite[Lemma~2.2]{cfz}} 
Let $C$ and $C'=C-\{\beta\}\cup\{\beta'\}$ be two adjacent 
clusters. 
Then the roots $\beta$ and $\beta'$ are exchangeable. 
\end{proposition} 
 
\pagebreak[3] 
 
The converse of Proposition~\ref{pr:two-terms-1} is also true: see 
Corollary~\ref{cor:exch-roots-converse} below. 
 
\begin{proposition} 
\label{pr:two-terms-2} 
{\rm \cite[Theorem~1.14]{cfz}} 
If $n>1$ and $\beta,\beta'\in\Phi_{\geq -1}$ are exchangeable, 
then the set 
\[ 
\{\sigma^{-1} (\sigma (\beta) + \sigma (\beta')): \sigma \in 
\langle \tau_+, \tau_- \rangle \} 
\] 
consists of two elements of $Q$, one of 
which is $\beta + \beta'$, and the other will be denoted by 
$\beta \subplus \beta'$. 
In the special case where $\beta'=-\alpha_j\,$, $j\in I$, we have 
\begin{equation} 
\label{eq:subplus-special} 
\begin{array}{rcl} 
(- \alpha_j)\subplus \beta &=& 
\tau_{- \varepsilon (j)}(- \alpha_j+\tau_{- \varepsilon (j)}(\beta)) \\[.1in] 
&=& \beta - \alpha_j + \sum_{i\neq j} a_{ij}\alpha_i\,. 
\end{array} 
\end{equation} 
\end{proposition} 
 
A precise rule for deciding whether an element $\sigma^{-1} (\sigma (\beta) + 
\sigma (\beta'))$ is equal to $\beta+\beta'$ or $\beta\subplus\beta'$ 
is given in Lemma~\ref{lem:cluster-expansion-properties-6} below. 
 
\begin{remark} 
If $n =1$, i.e., $\Phi$ is of type $A_1$ with a unique simple root 
$\alpha_1$, then $\{\beta, \beta'\} = \{- \alpha_1, \alpha_1\}$, 
and the group $\langle \tau_+, \tau_- \rangle$ is just the 
Weyl group $W = \langle s_1 \rangle$. 
Thus, in this case, the set 
in Proposition~\ref{pr:two-terms-2} consists of a single element $\beta + \beta' = 0$. 
It is then natural to set $\beta \subplus \beta' = 0$ as well. 
\end{remark}

\begin{remark} 
\label{rem:allow-decomposable} 
All results in this section extend in an obvious way to the 
case of an arbitrary Cartan matrix of finite type (not 
necessarily indecomposable). 
In that generality, $\Phi$ is a disjoint 
union of irreducible root systems $\Phi^{(1)}, \dots, \Phi^{(m)}$, and 
the clusters for $\Phi$ are the unions 
$C_1 \cup \cdots \cup C_m$, where each $C_k$ is a cluster for $\Phi^{(k)}$. 
\end{remark}

\section{Proof of Theorem~\ref{th:finite-type-class-new}} 
\label{sec:construction-of-finite-type} 
 
In this section, we complete the proof of 
Theorem~\ref{th:finite-type-class-new}. 
The plan is as follows. 
Without loss of generality, we can assume that 
the Cartan matrix~$A$ is indecomposable, so the corresponding 
root system~$\Phi$ is irreducible. 
By Theorem~\ref{th:cluster-fan} 
and Lemma~\ref{lem:triangulated-polytope}, 
conditions 
(\ref{eq:pseudomanifold-1})--(\ref{eq:geodesic-loops-generate}) 
are satisfied. 
By Proposition~\ref{prop:finite-type-criterion}, 
to prove Theorem~\ref{th:finite-type-class-new} it suffices 
to define a family of matrices $B(C)$, 
for each cluster $C$ in~$\Delta(\Phi)$, 
such that (\ref{eq:some-matrix-is-ss})--(\ref{eq:geodesic-loops}) 
hold, together with 
\begin{equation} 
\label{eq:some-is-ss-a} 
\text{for some cluster $C_\circ\,$, the matrix $B_\circ=B(C_\circ)$ is 
as in (\ref{eq:data-B}).} 
\hspace{.5in} 
\end{equation}

Defining the matrices $B(C)$ requires a little preparation. 
Throughout this section, all roots are presumed to belong to the set 
$\Phi_{\geq -1}\,$. 
 
\begin{lemma} 
\label{lem:sign-beta-beta'} 
There exists a sign function $(\beta, \beta') \mapsto 
\varepsilon(\beta, \beta')\in \{\pm 1\}$ on pairs of exchangeable roots, 
uniquely determined by the following properties: 
\begin{align} 
\label{eq:sign-neg-simple:1} 
&\varepsilon(- \alpha_j,\beta') = - \varepsilon (j)\,;\\ 
\label{eq:tau-sign} 
& \varepsilon(\tau \beta,\tau \beta')= - \varepsilon(\beta,\beta') 
\text{\ for $\tau\in\{\tau_+,\tau_-\}$ 
and $\beta,\beta' \notin 
\{-\alpha_j:\tau(-\alpha_j) =-\alpha_j\}$.} 
\end{align} 
Moreover, this function is skew-symmetric: 
\begin{equation} 
\label{eq:sign-symmetry} 
\varepsilon (\beta',\beta) = - \varepsilon(\beta, \beta'). 
\end{equation} 
\end{lemma} 
 
\proof 
The uniqueness of $\varepsilon (\beta, \beta')$ is an easy 
consequence of Theorem~\ref{th:dihedral}, Part~2. 
Let us prove the existence. 
For a root $\beta \in \Phi_{\geq -1}$ and a sign~$\varepsilon$, 
let $k_\varepsilon(\beta)$ denote 
the smallest nonnegative integer $k$ such that 
$
\tau^{(k+1)}_\varepsilon (\beta) = 
\tau^{(k)}_\varepsilon (\beta) \in-\Pi 
$, 
where we abbreviate 
\[\tau^{(k)}_\varepsilon = 
\underbrace{\tau_{\pm} 
\cdots \tau_{-\varepsilon}\tau_\varepsilon}_{k \text{~factors}} 
\] 
(cf.\ Theorem~\ref{th:dihedral}). 
In view of 
\cite[Theorem~3.1]{cfz}, 
we always have 
\begin{equation} 
\label{eq:sum=h+1} 
k_+(\beta)+k_-(\beta)=h+1; 
\end{equation} 
in particular, 
$k_{\varepsilon(j)}(-\alpha_j)=h+1$ and 
$k_{-\varepsilon(j)}(-\alpha_j)=0$. 
It follows from (\ref{eq:sum=h+1}) that if $\beta$ and $\beta'$ are 
incompatible (in particular, exchangeable), then 
$k_\varepsilon (\beta) < k_\varepsilon (\beta')$ for precisely one 
choice of a sign $\varepsilon$. 
Let us define $\varepsilon(\beta,\beta')$ by the condition 
\begin{equation} 
\label{eq:k<k} 
k_{\varepsilon(\beta,\beta')} (\beta) < k_{\varepsilon(\beta,\beta')} (\beta') 
\ . 
\end{equation} 
The properties 
(\ref{eq:sign-neg-simple:1})--(\ref{eq:sign-symmetry}) are 
immediately checked from this definition. 
\endproof 
 
We are now prepared to define the matrices $B(C)= (b_{\alpha\beta}(C))$.

\begin{definition} 
\label{def:b-matrices-via-roots} 
Let $C$ be a cluster in $\Delta(\Phi)$, 
that is, a $\ZZ$-basis of the root lattice $Q$ consisting of $n$ 
mutually compatible roots. 
Let $C'=C-\{\beta\}\cup\{\beta'\}$ be an 
adjacent cluster obtained from $C$ by exchanging 
a root $\beta\in C$ with some other root~$\beta'$. 
The entries 
$b_{\alpha\beta}(C)$, $\alpha\in C$, of the matrix $B(C)$ 
are defined by 
\begin{equation} 
\label{eq:b-matrix-via-roots} 
b_{\alpha\beta}(C) = \varepsilon(\beta,\beta') \cdot 
[\beta + \beta' - (\beta \subplus \beta'): 
\alpha]_C, 
\end{equation} 
where $[\gamma:\alpha]_C$ denotes the coefficient of $\alpha$ in the expansion 
of a vector $\gamma \in Q$ in the basis~$C$. 
\end{definition}

To complete the proof of Theorem~\ref{th:finite-type-class-new}, 
all we need to show is that the matrices $B(C)$ described in 
Definition~\ref{def:b-matrices-via-roots} satisfy 
(\ref{eq:some-is-ss-a}) and 
(\ref{eq:some-matrix-is-ss})--(\ref{eq:geodesic-loops}). 
 
\credit{Proof of {\rm(\ref{eq:some-is-ss-a})}} 
Let $C_\circ=-\Pi$ be the cluster consisting of all the negative simple roots. 
Applying (\ref{eq:b-matrix-via-roots}) and using 
(\ref{eq:sign-neg-simple:1}) and (\ref{eq:subplus-special}), we 
obtain 
\begin{equation} 
\label{eq:B(C_0)} 
b_{-\alpha_i,-\alpha_j}(C_\circ) = - \varepsilon (j) 
\cdot [- \sum_{k\neq j} a_{kj}\alpha_k: - \alpha_i]_{C_\circ} 
=\begin{cases} 
0 & \text{if $i=j$;} \\ 
\varepsilon(j) a_{ij} & \text{if $i\neq j$,} 
\end{cases} 
\end{equation} 
establishing 
(\ref{eq:some-is-ss-a}). 
\endproof 
 
\credit{Proof of {\rm(\ref{eq:some-matrix-is-ss})}} 
We start by summarizing the basic properties of cluster 
expansions of $\beta + \beta'$ and $\beta \subplus \beta'$ for 
an exchangeable pair of roots. 
 
 
\begin{lemma} 
\label{lem:cluster-expansion-properties} 
Let $\beta,\beta'\in\Phi_{\geq -1}$ be exchangeable. 
 
\begin{enumerate} 
 
\item 
No negative simple root can be a cluster component of $\beta + \beta'$. 
 
\item 
The vectors $\beta + \beta'$ and $\beta \subplus \beta'$ 
have no common cluster components. 
That is, no root in $\Phi_{\geq -1}$ can contribute, 
with nonzero coefficient, 
to the cluster expansions of both $\beta + \beta'$ and $\beta \subplus 
\beta'$. 
 
\item 
All cluster components of $\beta + \beta'$ and $\beta \subplus \beta'$ are 
compatible with both $\beta$ and~$\beta'$. 
 
\item 
A root $\alpha \in \Phi_{\geq -1}\,$, $\alpha\notin\{\beta,\beta'\}$, 
is compatible with both $\beta$ and~$\beta'$ if and 
only if it is compatible with all cluster components of 
$\beta + \beta'$ and $\beta \subplus \beta'$. 
 
\item 
If $\alpha \in - \Pi$ is compatible 
with all cluster components of $\beta + \beta'$, 
then it is compatible 
with all cluster components of $\beta \subplus \beta'$. 
 \end{enumerate} 
\end{lemma}

\begin{proof} 
 
 1. Suppose $[\beta+ \beta': - \alpha_i]_{\rm clus} > 0$ 
for some $i \in I$. 
(Here we use the notation from 
Proposition~\ref{prop:inv-cluster-expansion-2}.) 
Since all roots compatible with $- \alpha_i$ and different 
from $- \alpha_i$ do not contain $\alpha_i$ in their simple root 
expansion, it follows that $[\beta+ \beta': \alpha_i] < 0$. 
This can only happen if one of the roots $\beta$ and $\beta'$, say 
$\beta$, is equal to $- \alpha_i$; but since $\beta'$ is 
incompatible with $\beta$, we will still have $[\beta+ \beta': \alpha_i] \geq 0$, 
a contradiction. 
 
 2. Suppose $\alpha$ is a  common cluster component of 
$\beta + \beta'$ and $\beta \subplus \beta'$. 
Applying if necessary a transformation from $\langle \tau_+, 
\tau_- \rangle$, we can assume that $\alpha$ is negative simple 
(see Proposition~\ref{prop:inv-cluster-expansion-2}, 
and Theorem~\ref{th:dihedral}, Part~2). 
But this is impossible by Part~1. 
 
 3. The claim for $\beta + \beta'$ is proved in \cite[Lemma~2.3]{cfz}. 
Since $\beta \subplus \beta' = \sigma^{-1} (\sigma (\beta) + 
\sigma (\beta'))$ for some $\sigma \in \langle \tau_+, \tau_- \rangle$, 
the claim for  $\beta \subplus \beta'$ follows from 
Proposition~\ref{prop:inv-cluster-expansion-2}. 
 
 4. First suppose that $\alpha$ is compatible with both $\beta$ 
and~$\beta'$. 
The fact that $\alpha$ is compatible with all cluster components of 
$\beta + \beta'$ is proved in \cite[Lemma~2.3]{cfz}. 
The fact that $\alpha$ is compatible with all cluster components of 
$\beta \subplus \beta'$ now follows in the same way as in Part 3. 
 
To prove the converse, suppose that $\alpha$ is incompatible with $\beta$. 
As in Part 2 above, we can assume that $\alpha = - \alpha_i$ for some $i$. 
Thus, we have $[\beta : \alpha_i] > 0$. 
Since $\alpha \neq \beta'$, it follows that 
$[\beta + \beta': \alpha_i] > 0$ as well. 
Therefore, $[\gamma : \alpha_i] > 0$ for some cluster component 
$\gamma$ of $\beta + \beta'$, and we are done. 
 
5. This follows from \cite[Theorem~1.17]{cfz}. 
\end{proof} 
 
As a corollary of Lemma~\ref{lem:cluster-expansion-properties}, 
we obtain the converse of Proposition~\ref{pr:two-terms-1}. 
 
\begin{corollary} 
\label{cor:exch-roots-converse} 
Let $\beta$ and $\beta'$ be exchangeable almost positive roots. 
Then there exist two adjacent clusters $C$ and $C'$ 
such that $C'=C-\{\beta\}\cup\{\beta'\}$. 
\end{corollary} 
 
\proof 
By Lemma~\ref{lem:cluster-expansion-properties}, Parts 3 and 4, 
the set consisting of $\beta$ and all cluster components of 
$\beta + \beta'$ and $\beta \subplus \beta'$ is compatible 
(i.e., all its elements are mutually compatible). 
Thus, there exists a cluster $C$ containing this set. 
Again using Lemma~\ref{lem:cluster-expansion-properties}, Part~4, 
we conclude that every element of $C - \{\beta\}$ is compatible 
with $\beta'$. 
Hence  $C-\{\beta\}\cup\{\beta'\}$ is a cluster, as desired. 
\endproof 
 
\begin{corollary} 
\label{cor:cluster-components-lie-in-cluster} 
Let $C$ and $C'=C-\{\beta\}\cup\{\beta'\}$ be adjacent clusters . 
Then all cluster components of $\beta + \beta'$ and $\beta \subplus \beta'$ 
belong to $C\cap C'=C - \{\beta\}$. 
\end{corollary}

\begin{proof} 
By Lemma~\ref{lem:cluster-expansion-properties}, Parts 3 and~4, 
any cluster component of  $\beta + \beta'$ or $\beta \subplus \beta'$ 
is compatible with every element of~$C$, 
hence must belong to~$C$. 
\end{proof} 
 
We next 
derive a useful alternative description of the matrix entries of $B(C)$. 
 
\begin{lemma} 
\label{lem:b-thru-cluster-expansion} 
In the situation of Definition~\ref{def:b-matrices-via-roots}, 
we have 
\begin{equation} 
\label{eq:b-matrix-via-roots-2} 
b_{\alpha\beta}(C) = \varepsilon(\beta,\beta') \cdot 
([\beta + \beta': \alpha]_{\rm clus} - 
[\beta \subplus \beta': \alpha]_{\rm clus}) \ . 
\end{equation} 
In particular, the entry $b_{\alpha\beta}(C)$ is uniquely 
determined by $\alpha$, $\beta$, and~$\beta'$. 
\end{lemma} 
 
\proof 
By Corollary~\ref{cor:cluster-components-lie-in-cluster}, 
the cluster expansions of the 
vectors $\beta + \beta'$ and $\beta \subplus \beta'$ 
are the same as their basis $C$ expansions. 
Thus, (\ref{eq:b-matrix-via-roots}) is equivalent to 
(\ref{eq:b-matrix-via-roots-2}). 
\endproof 
 
We will need the following lemma (cf.\ Proposition~\ref{pr:two-terms-2}). 
 
\begin{lemma} 
\label{lem:cluster-expansion-properties-6} 
For a pair of exchangeable roots $\beta$ and $\beta'$ and a sign 
$\varepsilon$, we have 
\[ 
(\tau^{(k)}_\varepsilon)^{-1} 
(\tau^{(k)}_\varepsilon (\beta)+\tau^{(k)}_\varepsilon (\beta')) 
= 
\begin{cases} 
\beta + \beta' & 
\text{if $0 \leq k \leq \min(k_\varepsilon(\beta),k_\varepsilon(\beta'))$;} \\[.05in] 
\beta \subplus \beta' & 
\text{if $\min(k_\varepsilon(\beta),k_\varepsilon(\beta')) 
< k$} \\ 
&\hspace{1in}\leq \max(k_\varepsilon(\beta),k_\varepsilon(\beta')). 
\end{cases} 
\] 
\end{lemma} 
 
\begin{proof} 
This is a consequence of \cite[Lemma~3.2]{cfz}. 
\end{proof}

We next establish the following symmetry property. 
 
\begin{lemma} 
\label{lem:B-are-tau-invariant} 
For $\tau\in\{\tau_+\,,\tau_-\}$, we have 
$b_{\tau \alpha,\tau\beta}(\tau C) = - b_{\alpha\beta}(C)$. 
\end{lemma} 
 
\pagebreak[2] 
 
\begin{proof} 
Let us introduce some notation. 
For every pair of exchangeable roots 
$(\beta, \beta')$, we denote by $S_+ (\beta, \beta')$ and 
$S_- (\beta, \beta')$ the two elements of $Q$ given by 
\begin{equation} 
\label{eq:Spm-def} 
\begin{array}{r} 
S_{\varepsilon(\beta,\beta')}(\beta, \beta') = \beta + \beta', \\[.05in] 
S_{-\varepsilon(\beta,\beta')}(\beta, \beta') = \beta \subplus 
\beta'.\, 
\end{array} 
\end{equation} 
In this notation, 
(\ref{eq:b-matrix-via-roots-2}) takes the form 
\begin{equation} 
\label{eq:b-matrix-via-roots-3} 
b_{\alpha\beta}(C) = 
[S_+ (\beta, \beta'): \alpha]_{\rm clus} - 
[S_- (\beta, \beta'): \alpha]_{\rm clus} \ . 
\end{equation} 
The functions $S_\pm$ satisfy 
the following symmetry property: 
\begin{equation} 
\label{eq:tau-Spm} 
\tau\,S_\varepsilon(\beta,\beta')=S_{-\varepsilon}(\tau\beta,\tau\beta') 
\text{\ \ for  $\tau\in\{\tau_+,\tau_-\}$}; 
\end{equation} 
this follows by comparing (\ref{eq:tau-sign}) with 
Lemma~\ref{lem:cluster-expansion-properties-6}. 
The lemma follows by combining (\ref{eq:b-matrix-via-roots-3}) 
and (\ref{eq:tau-Spm}) with Proposition~\ref{prop:inv-cluster-expansion-2}. 
\end{proof} 
 
We are finally ready for the task at hand: verifying that each 
matrix $B(C)$ is sign-skew-symmetric. 
In view of (\ref{eq:b-matrix-via-roots-2}) and 
Lemma~\ref{lem:cluster-expansion-properties}, Part~2, 
the signs of its entries are given as follows: 
\begin{equation} 
\label{eq:sign-of-b} 
{\rm sgn}(b_{\alpha \beta}(C))= 
\begin{cases} 
\varepsilon(\beta, \beta') 
& \text{if $\alpha$ is a cluster component of $\beta + \beta'$;}\\ 
-\varepsilon(\beta, \beta') 
& \text{if $\alpha$ is a cluster component of $\beta \subplus \beta'$;}\\ 
0 & \text{otherwise.} 
\end{cases} 
\end{equation} 
Since $\beta$ is incompatible with $\beta'$, 
Lemma~\ref{lem:cluster-expansion-properties}, Part~3, ensures that 
$\beta$ cannot be a cluster component of 
$\beta + \beta'$ or $\beta \subplus \beta'$, 
so all diagonal entries of $B(C)$ are equal to~$0$. 
 
\begin{lemma} 
\label{lem:b-vanishing} 
Let $\alpha$ and $\beta$ be two different 
elements of a cluster $C$, and let the corresponding adjacent clusters be 
$C - \{\alpha\} \cup \{\alpha'\}$ and $C - \{\beta\} \cup \{\beta'\}$. 
Then we have $b_{\alpha\beta}(C) = 0$ if and only if 
$\alpha'$ is compatible with $\beta'$. 
\end{lemma} 
 
\proof 
By (\ref{eq:sign-of-b}), condition $b_{\alpha\beta}(C) = 0$ is 
equivalent to 
\begin{equation} 
\label{eq:both-sums-vanish} 
[\beta + \beta': \alpha]_{\rm clus}= 
[\beta \subplus \beta': \alpha]_{\rm clus}=0. 
\end{equation} 
Since all cluster components of $\beta + \beta'$ and 
$\beta \subplus \beta'$ belong to~$C$ 
(see Corollary~\ref{cor:cluster-components-lie-in-cluster}), and since $\alpha'$ is 
compatible with all elements of $C - \{\alpha\}$, 
condition (\ref{eq:both-sums-vanish}) holds if and only if 
 $\alpha'$ is compatible with 
all cluster components of $\beta + \beta'$ and 
$\beta \subplus \beta'$. 
By Lemma~\ref{lem:cluster-expansion-properties}, Part~4, 
this is in turn equivalent to 
$\alpha'$ being compatible with~$\beta'$. 
\endproof

Since the compatibility relation is symmetric, Lemma~\ref{lem:b-vanishing} 
implies that \hbox{$b_{\alpha\beta}(C)=0$} is equivalent to $b_{\beta \alpha}(C)=0$, 
as needed. 
 
It remains to consider the case where both $b_{\alpha\beta}(C)$ and 
$b_{\beta \alpha}(C)$ are nonzero. 
By (\ref{eq:sign-of-b}), this means that $\alpha$ is a cluster 
component of $\beta + \beta'$ or $\beta \subplus \beta'$, while 
$\beta$ is a cluster component of $\alpha + \alpha'$ or $\alpha \subplus \alpha'$. 
By Lemma~\ref{lem:B-are-tau-invariant}, 
it is enough to consider the special case $\alpha=-\alpha_j\in-\Pi$. 
Let us abbreviate $\varepsilon_\circ = \varepsilon(\beta, \beta')$. 
In view of Lemma~\ref{lem:cluster-expansion-properties}, Part~1, 
(\ref{eq:sign-of-b}) yields 
${\rm sgn}(b_{\alpha \beta}(C))= -\varepsilon_\circ$. 
Interchanging $\alpha$ and $\beta$, and applying the same rule 
(\ref{eq:sign-of-b}), we obtain, taking into account 
(\ref{eq:sign-neg-simple:1}), that 
\begin{equation*} 
{\rm sgn}(b_{\beta \alpha}(C))= 
\begin{cases} 
-\varepsilon(j) 
& \text{if $\beta$ is a cluster component of $\alpha + \alpha'$;}\\ 
\varepsilon(j) 
& \text{if $\beta$ is a cluster component of $\alpha \subplus \alpha'$.} 
\end{cases} 
\end{equation*} 
By Lemma~\ref{lem:cluster-expansion-properties}, Part~1, 
and Proposition~\ref{prop:inv-cluster-expansion-2}, 
the element of the set $\{\alpha + \alpha',\alpha \subplus \alpha'\}$ 
that does not have $\beta$ as a cluster component is equal to 
\[ 
(\tau^{(k_\circ)}_{\varepsilon_\circ})^{-1} 
(\tau^{(k_\circ)}_{\varepsilon_\circ} (\alpha) + 
\tau^{(k_\circ)}_{\varepsilon_\circ}(\alpha')) , 
\] 
where $k_\circ = k_{\varepsilon_\circ} (\beta)$. 
Our goal is to prove that ${\rm sgn}(b_{\beta 
  \alpha}(C))=\varepsilon_\circ$, 
which can now be restated as follows: 
\begin{equation} 
\label{eq:sum-or-subsum} 
(\tau^{(k_\circ)}_{\varepsilon_\circ})^{-1} 
(\tau^{(k_\circ)}_{\varepsilon_\circ} (\alpha) + 
\tau^{(k_\circ)}_{\varepsilon_\circ}(\alpha')) 
=\begin{cases} 
\alpha + \alpha' & \text{if $\varepsilon_\circ = \varepsilon(j)$;} \\[.05in] 
\alpha \subplus \alpha' & \text{if 
$\varepsilon_\circ = -\varepsilon(j)$.} 
\end{cases} 
\end{equation} 
It was shown in the proof of Lemma~\ref{lem:sign-beta-beta'} that 
$k_\circ \!<\! k_{\varepsilon_\circ} (\beta')$. 
As a main step towards proving (\ref{eq:sum-or-subsum}), we show that 
$k_\circ \!<\! k_{\varepsilon_\circ} (\alpha')$. 
Suppose on the contrary that 
$k_{\varepsilon_\circ} (\alpha')\! \leq\! k_\circ$. 
Applying Lemma~\ref{lem:cluster-expansion-properties-6} with 
$k = k_{\varepsilon_\circ} (\alpha') \leq 
\min(k_{\varepsilon_\circ} (\beta),k_{\varepsilon_\circ}(\beta'))$, 
we see that 
$$\tau^{(k)}_{\varepsilon_\circ} (\beta+\beta') = 
\tau^{(k)}_{\varepsilon_\circ} (\beta)+\tau^{(k)}_{\varepsilon_\circ} (\beta')\ .$$ 
To arrive at a contradiction, notice that $\alpha'$ is compatible 
with all cluster components of $\beta + \beta'$ but is 
incompatible with $\beta'$ (see Lemma~\ref{lem:b-vanishing}). 
It follows that the root 
$\tau^{(k)}_{\varepsilon_\circ} (\alpha') \in - \Pi$ is compatible 
with all cluster components of 
$\tau^{(k)}_{\varepsilon_\circ} (\beta)+\tau^{(k)}_{\varepsilon_\circ}(\beta')$ 
but incompatible with $\tau^{(k)}_{\varepsilon_\circ} (\beta')$. 
But this contradicts Lemma~\ref{lem:cluster-expansion-properties}, Parts~4--5. 
 
Having established the inequality $k_\circ < k_{\varepsilon_\circ}(\alpha')$, 
we see that (\ref{eq:sum-or-subsum}) becomes a special case of 
Lemma~\ref{lem:cluster-expansion-properties-6}. 
Indeed, if $\varepsilon_\circ = \varepsilon(j)$ then 
$k_{\varepsilon_\circ}(\alpha) = k_{\varepsilon_\circ}(-\alpha_j) = h+1$, 
so  $k_\circ < \min(k_{\varepsilon_\circ}(\alpha), k_{\varepsilon_\circ}(\alpha'))$; 
and if $\varepsilon_\circ = -\varepsilon(j)$ then 
$k_{\varepsilon_\circ}(\alpha) = k_{\varepsilon_\circ}(-\alpha_j) = 0$, 
so  $\min(k_{\varepsilon_\circ}(\alpha), k_{\varepsilon_\circ}(\alpha')) 
\!\leq\! k_\circ \!<\! \max(k_{\varepsilon_\circ}(\alpha), 
k_{\varepsilon_\circ}(\alpha'))$. 
This concludes our proof of~(\ref{eq:some-matrix-is-ss}). 
\endproof 
 
\credit{Proof of {\rm(\ref{eq:B-mutation-gamma})}} 
Let $C$ and $\overline C = C-\{\gamma\}\cup\{\overline\gamma\}$ be 
adjacent clusters. 
Our task is to show that every matrix entry 
$b_{\alpha \beta}(\overline C)$ is obtained from the entries of 
$B(C)$ by the procedure described in (\ref{eq:B-mutation-gamma}). 
We already know that all diagonal entries in $B(C)$ and 
$B(\overline C)$ are equal to $0$ which is consistent with the 
matrix mutation rule (\ref{eq:mutation}) (since we have already proved 
that these matrices are sign-skew-symmetric). 
So let us assume that $\alpha \neq \beta$. 
If $\beta = \overline \gamma$ then the desired equality 
$b_{\alpha \overline \gamma}(\overline C) = - b_{\alpha \gamma}(C)$ 
is immediate from (\ref{eq:b-matrix-via-roots-2}) 
and (\ref{eq:sign-symmetry}). 
Thus, it remains to treat the case where $\alpha \neq \beta$ and 
$\beta \neq \overline \gamma$. 
By Lemma~\ref{lem:B-are-tau-invariant}, 
it is enough to consider the special case $\beta=-\alpha_j\in-\Pi$. 
 
 
In view of (\ref{eq:b-matrix-via-roots}), (\ref{eq:subplus-special}) 
and (\ref{eq:sign-neg-simple:1}), we have 
\begin{equation} 
\label{eq:b-special} 
\begin{array}{l} 
b_{\alpha \beta}(C) = - \varepsilon(j) [\delta : \alpha]_C\,, \\[.05in] 
b_{\alpha \beta}(\overline C) = 
- \varepsilon(j) [\delta : \alpha]_{\overline C}\,, 
\end{array} 
\end{equation} 
where we use the notation 
\[ 
\delta = \delta_j = - \textstyle\sum_{i \neq j} a_{ij} \alpha_i \,. 
\] 
Since the basis $\overline C$ is obtained from $C$ by replacing 
$\gamma$ by $\overline \gamma$, we can use the expansion 
\[ 
\gamma + \overline \gamma = 
\sum_{\alpha \in C \cap \overline C} 
[\gamma + \overline \gamma : \alpha]_{\rm clus} \cdot  \alpha 
\] 
(cf.\  Corollary~\ref{cor:cluster-components-lie-in-cluster}) 
to obtain: 
\begin{equation} 
\label{eq:C-Cbar-transition} 
\begin{array}{rcl} 
[\delta : \overline \gamma]_{\overline C} &=& 
-[\delta : \gamma]_{C} \\[.05in] 
[\delta : \alpha]_{\overline C} &=& [\delta : \alpha]_{C} 
+ [\gamma + \overline \gamma : \alpha]_{\rm clus} 
\cdot [\delta : \gamma]_{C} \,\, \text{\ \ for 
$\alpha \in C \cap \overline C$} \,. 
\end{array} 
\end{equation} 
Combining (\ref{eq:b-special}) and (\ref{eq:C-Cbar-transition}), we get 
\begin{align} 
\label{eq:C-Cbar-b-transition-1} 
b_{\overline \gamma \beta}(\overline C) &= 
- b_{\gamma \beta}(C) \\ 
\label{eq:C-Cbar-b-transition-2} 
b_{\alpha \beta}(\overline C) &= b_{\alpha \beta}(C) 
+ [\gamma + \overline \gamma : \alpha]_{\rm clus} 
\cdot b_{\gamma \beta}(C) \,\, \text{\ \ for 
$\alpha \in C \cap \overline C$} \,. 
\end{align} 
Formula (\ref{eq:C-Cbar-b-transition-1}) takes care of 
the first case in the mutation rule (\ref{eq:mutation}). 
If $b_{\gamma \beta}(C) = 0$, then 
(\ref{eq:C-Cbar-b-transition-2}) again agrees with (\ref{eq:mutation}). 
So let us assume that $b_{\gamma \beta}(C) \neq 0$. 
Applying 
(\ref{eq:b-matrix-via-roots-2}) to $b_{\alpha \gamma} (C)$ 
and $b_{\beta \gamma} (C)$, we obtain: 
\begin{align} 
\label{eq:coeff-of-alpha-in-gamma+gamma} 
[\gamma + \overline \gamma : \alpha]_{\rm clus} &= 
\max(\varepsilon(\gamma, \overline \gamma)b_{\alpha \gamma} (C), 0), 
\\ 
\label{eq:coeff-of-beta-in-gamma+gamma} 
[\gamma + \overline \gamma : \beta]_{\rm clus} &= 
\max(\varepsilon(\gamma, \overline \gamma)b_{\beta \gamma} (C), 0). 
\end{align} 
By Lemma~\ref{lem:cluster-expansion-properties}, Part~1, the left-hand 
side of (\ref{eq:coeff-of-beta-in-gamma+gamma}) equals~0. 
Hence 
\begin{equation} 
\label{eq:varepsilon(gamma,gamma)} 
\varepsilon(\gamma, \overline \gamma) = - {\rm sgn}(b_{\beta \gamma}(C)) 
= {\rm sgn}(b_{\gamma \beta}(C)) 
\end{equation} 
(using that $B(C)$ is sign-skew-symmetric). 
Now (\ref{eq:C-Cbar-b-transition-2}), 
(\ref{eq:coeff-of-alpha-in-gamma+gamma}), and 
(\ref{eq:varepsilon(gamma,gamma)}) give 
$$b_{\alpha \beta}(\overline C) =  b_{\alpha \beta}(C) + 
\max({\rm sgn}(b_{\gamma \beta}(C)) \cdot b_{\alpha \gamma} (C), 0) 
\cdot b_{\gamma \beta}(C),$$ 
which is easily seen to be equivalent to the second case in 
(\ref{eq:mutation}). 
This completes the verification that the matrices $B(C)$ satisfy 
(\ref{eq:B-mutation-gamma}). 
\endproof


\credit{Proof of {\rm(\ref{eq:geodesic-loops})}} 
We proceed by induction on the rank~$n$ of a root system~$\Phi$. 
For induction purposes, we need to allow $\Phi$ to be 
reducible; this is possible in view of 
Remark~\ref{rem:allow-decomposable}. 
For $n=2$, the generalized associahedra of types $A_1 \times A_1$, 
$A_2$, $B_2$ and $G_2$ are convex polygons with $4$, 
$5$, $6$, and~$8$ sides, respectively, matching the claim. 
 
For the induction step, consider a geodesic loop $L$ in the dual graph 
$\Gamma = \Gamma_{\Delta (\Phi)}$ for a root system 
$\Phi$ of rank $n \geq 3$. 
According to the definition of a geodesic, 
the clusters lying on $L$ are obtained from some initial cluster 
$C$ by fixing \hbox{$n-2$} of the $n$ roots and alternately exchanging the 
remaining two roots. 
By (\ref{eq:compatibility-2}), 
every transformation 
from $\langle \tau_+, \tau_- \rangle$ sends  geodesics to geodesics. 
Furthermore, the Coxeter number associated to a geodesic does not 
change, by Lemma~\ref{lem:B-are-tau-invariant}. 
Using Theorem~\ref{th:dihedral}, Part~2, we may therefore assume, 
without loss of generality, 
that one of the $n-2$ fixed roots in the initial cluster $C$ is $-\alpha_j \in-\Pi$. 
 
\begin{lemma} 
\label{lem:restriction} 
Let $\Phi'$ be the rank $n-1$ root subsystem of~$\Phi$ spanned by the 
simple roots $\alpha_i$ for $i\neq j$. 
 
\begin{enumerate} 
 
\item 
The correspondence $C' \mapsto \{-\alpha_j\} \cup C'$ is a 
bijection between the clusters in $\Delta (\Phi')$ and the clusters in 
$\Delta (\Phi)$ that contain $-\alpha_j$. 
Thus, it identifies the cluster complex $\Delta (\Phi')$ 
with the link $\Delta_{\{-\alpha_j\}}$ of $\{-\alpha_j\}$ 
in $\Delta(\Phi)$ (see Section~\ref{sec:pseudomanifolds}). 
 
\item 
Assume that the sign function {\rm (\ref{eq:sign-function})} 
for $\Phi'$ is a restriction of the sign function for~$\Phi$. 
The matrix $B(C')$ 
associated to a cluster $C' \subset \Phi'_{\geq -1}$ can be 
obtained from the matrix $B(\{-\alpha_j\} \cup C')$ by crossing 
out the row and column corresponding to the root~$-\alpha_j\,$. 
\end{enumerate} 
\end{lemma} 
 
\proof 
Part 1 follows from~\cite[Proposition~3.5 (3)]{yga}. 
The assertion in Part 2 is immediately checked in the special case $C' 
=\{-\alpha_i: i \neq j\}$ (see (\ref{eq:B(C_0)})). 
It is then extended to an arbitrary cluster $C'$ because 
the graph $\Gamma' = \Gamma_{\Delta (\Phi')}$ is 
connected, and the propagation rules for the matrices 
 $B(C')$ and $B(\{-\alpha_j\} \cup C')$ are easily seen to be exactly the same. 
(Here we use the fact that conditions 
(\ref{eq:pseudomanifold-1}), (\ref{eq:some-matrix-is-ss}), and (\ref{eq:B-mutation-gamma}) 
have already been checked for $\Gamma$ and $\Gamma'$ alike.) 
\endproof 
 
By Lemma~\ref{lem:restriction}, Part~1, 
$L$ can be viewed as a geodesic loop in 
$\Gamma' = \Gamma_{\Delta (\Phi')}$. 
To complete the induction step, it remains to notice that, 
in view of Lemma~\ref{lem:restriction}, Part~2, 
the Coxeter number associated 
with $L$ in $\Gamma'$ coincides with the original value 
in~$\Gamma$. 
This completes the proof of (\ref{eq:geodesic-loops}). 
\endproof 
 
Theorem~\ref{th:finite-type-class-new} is proved. 
 
\pagebreak[3] 
 
\section{Proofs of Theorems~\ref{th:cluster-variable-denominators-1} 
and \ref{th:exch-rel-dep-on-vars-only}--\ref{th:finite-type-complex}} 
\label{sec:proofs-laurent-etc} 
 
\credit{Proof of 
  Theorems~\ref{th:exch-rel-dep-on-vars-only}--\ref{th:finite-type-complex} 
modulo Theorem~\ref{th:cluster-variable-denominators-1}} 
To take 
Theorems~\ref{th:exch-rel-dep-on-vars-only}--\ref{th:finite-type-complex} 
out of the way, we begin by deducing 
them from Theorem~\ref{th:cluster-variable-denominators-1}. 
We adopt all the conventions and notation of 
Section~\ref{sec:construction-of-finite-type}. 
In particular, we assume, without loss of generality, 
that the Cartan matrix $A$ is indecomposable, so the corresponding 
(finite) root system $\Phi$ is irreducible. 
We have proved that the complex $\Delta (\Phi)$ on the ground set 
$\Phi_{\geq -1}$ of almost positive roots, together with the 
family of matrices $B(C)$ introduced in 
Definition~\ref{def:b-matrices-via-roots}, satisfy conditions 
(\ref{eq:pseudomanifold-1})--(\ref{eq:geodesic-loops}).

Let $\AA=\AA(B_\circ,\pp_\circ)$ be the 
cluster algebra of finite type appearing in 
Theorem~\ref{th:finite-type-class-new}. 
Here we choose an initial seed $\Sigma_\circ=(\xx_\circ,\pp_\circ,B_\circ)$ 
for~$\AA$ by identifying the matrix $B_\circ$ 
with the matrix $B(C_\circ)$ at the cluster $C_\circ=-\Pi$ in 
$\Delta(\Phi)$ (see (\ref{eq:B(C_0)})). 
This gives us a seed attachment at $C_\circ$. 
Applying Lemma~\ref{lem:attachments}, we obtain a surjection $\alpha 
\mapsto x[\alpha]$ from $\Phi_{\geq -1}$ onto the set of all cluster 
variables in~$\AA$. 
Note that at this point, we have not yet proved that the variables 
$x[\alpha]$ are all distinct.

Assume for a moment that 
Theorem~\ref{th:cluster-variable-denominators-1} has been 
established. 
Then the map $\alpha \mapsto x[\alpha]$ is a bijection, 
and Theorems~\ref{th:pseudomanifold-clusters} and~\ref{th:finite-type-complex} 
follow by Remark~\ref{rem:isom-cluster-complex}. 
 
As for Theorem~\ref{th:exch-rel-dep-on-vars-only}, it 
becomes a consequence of 
Lemma~\ref{lem:b-thru-cluster-expansion}. 
To be more precise, let us associate to every lattice vector $\gamma \in Q$ 
a monomial in the cluster variables by setting 
\[ 
x[\gamma] = \prod_\alpha x[\alpha]^{m_\alpha}, \quad 
m_\alpha = [\gamma : \alpha]_{\rm clus} \, . 
\] 
In view of (\ref{eq:b-matrix-via-roots-2}), 
every exchange relation (\ref{eq:exchange-rel-xx}) 
corresponding to adjacent clusters 
$C$ and $C-\{\beta\}\cup\{\beta'\}$ 
can be written in the form 
\begin{equation} 
\label{eq:exch-rel-fin-type} 
x[\beta]\, x[\beta'] = 
p^{\varepsilon(\beta,\beta')}_{\beta} (C) \, x[\beta+\beta'] + 
p^{-\varepsilon(\beta,\beta')}_{\beta}(C) \, x[\beta \subplus \beta'], 
\end{equation} 
for some coefficients $p^{\pm}_{\beta} (C) \in \PP$. 
Thus, the set of cluster variables and the respective nonzero 
exponents that appear in the right-hand 
side of (\ref{eq:exch-rel-fin-type}) 
are uniquely determined by $\beta$ and~$\beta'$. 
The same holds for the coefficients $p^{\pm}_{\beta} (C)$, 
since the cluster variables appearing in the right-hand 
side are algebraically independent. 
\endproof 
 
We denote $p^\pm_{\beta,\beta'}=p^{\pm}_{\beta} (C)$. 
This notation is justified in view of 
Theorem~\ref{th:exch-rel-dep-on-vars-only}. 
 
\begin{remark} 
\label{rem:exchange-relations-special} 
In view of Corollary~\ref{cor:exch-roots-converse}, the exchange relation 
(\ref{eq:exch-rel-fin-type}) holds for \emph{every pair} $(\beta,\beta')$ 
of exchangeable roots. 
Also note that, in view of 
(\ref{eq:subplus-special}) and (\ref{eq:sign-neg-simple:1}), 
the exchange relation (\ref{eq:exch-rel-fin-type}) takes the 
following more explicit form if $\beta'$ is negative simple: 
\begin{align} 
\label{eq:special-cluster-expansion} 
x[\beta]\, x[- \alpha_j] 
&= p^{\varepsilon(j)}_{\beta,-\alpha_j} \, x[\beta-\alpha_j] 
+ p^{-\varepsilon(j)}_{\beta,-\alpha_j} \, x[\beta\subplus (-\alpha_j)] 
\\ 
\nonumber 
&= 
p^{\varepsilon(j)}_{\beta,-\alpha_j} \, x[\beta-\alpha_j] 
 + p^{-\varepsilon(j)}_{\beta,-\alpha_j} \, 
 x[\beta-\alpha_j + \textstyle\sum_{i\neq j} a_{ij}\alpha_i] \ . 
\end{align} 
For the classical types, the list of all exchangeable pairs 
$(\beta, - \alpha_j)$, together with the explicitly given cluster 
expansions for $\beta-\alpha_j$ and 
$\beta\subplus (-\alpha_j) = 
\beta-\alpha_j + \textstyle\sum_{i\neq j} a_{ij}\alpha_i$, 
was given in \cite[Section~4]{cfz}. 
\end{remark}

\credit{Proof of Theorem~\ref{th:cluster-variable-denominators-1}} 
We prove (\ref{eq:laurent-finite-type}) by induction on 
$$k(\alpha) = \min (k_+ (\alpha), k_- (\alpha)) \geq 0$$ 
(see the proof of Lemma~\ref{lem:sign-beta-beta'}). 
If $k(\alpha) = 0$, then $\alpha$ is a negative simple root, and 
there is nothing to prove. 
So we assume that $k(\alpha) = k \geq 1$, and that (\ref{eq:laurent-finite-type}) 
holds for all roots $\alpha'$ with $k(\alpha') < k$. 
 
By the definition of $k(\alpha)$, we have 
$$\alpha = \tau^{(k)}_{\varepsilon (j)} (- \alpha_j) = 
\tau^{(k-1)}_{-\varepsilon (j)} (\alpha_j)$$ 
for some $j \in I$. 
Since $\alpha_j$ and $-\alpha_j$ are exchangeable, 
so are $\alpha$ and $\tau (-\alpha_j)$, 
where we abbreviate $\tau = \tau^{(k-1)}_{-\varepsilon (j)}$. 
Let us write the corresponding exchange relation. 
Using the $\langle \tau_\pm \rangle$-invariance of 
the exponents appearing in exchange relations 
(Lemma~\ref{lem:B-are-tau-invariant}), together with 
(\ref{eq:b-matrix-via-roots-2}) and (\ref{eq:subplus-special}), 
we obtain: 
\begin{equation} 
\label{eq:exchange-equator} 
x[\alpha]\, x[\tau (-\alpha_j)]= 
q \prod_{i \neq j} x[\tau (-\alpha_i)]^{-a_{ij}} + r, 
\end{equation} 
where $q,r \in \PP$. 
For $k =1$, we have $\alpha = \alpha_j$, and (\ref{eq:exchange-equator}) yields 
$$x[\alpha_j] = 
\frac{q \prod_{i \neq j} x_i^{-a_{ij}} + r}{x_j}\,,$$ 
establishing (\ref{eq:laurent-finite-type}). 
Thus, we may assume that $k \geq 2$. 
In this case, all the roots $\alpha' \neq \alpha$ that appear in 
(\ref{eq:exchange-equator}) are positive with 
$k(\alpha') < k$. 
Abbreviating 
$$\gamma = \sum_{i \neq j} (-a_{ij}) \cdot \tau(-\alpha_i)$$ 
and applying the induction assumption, we can rewrite 
(\ref{eq:exchange-equator}) as 
\begin{equation} 
\label{eq:intermediate} 
x[\alpha] = x^{\tau (-\alpha_j) - \gamma} \cdot 
\frac{q \prod_{i \neq j} P_{\tau(-\alpha_i)}^{-a_{ij}} + r x^\gamma} 
{P_{\tau(-\alpha_j)}} \, , 
\end{equation} 
where all $P_{\alpha'}$ are polynomials over $\ZZ\PP$ 
in the variables from the initial cluster 
$\xx_\circ$ with nonzero constant terms. 
The  next step of the proof relies on the following trivial lemma. 
 
\begin{lemma} 
\label{lem:laurent-implies-poly} 
Let $P$ and $Q$ be two polynomials (in any number of variables) 
with coefficients in a domain~$S$, and with nonzero constant terms $a$ and~$b$, 
respectively. If the ratio $P/Q$ is a Laurent polynomial over~$S$, then it 
is in fact a polynomial over $S$ with the constant term~$a/b$. 
\end{lemma} 
 
By \cite[Theorem~3.1]{fz-clust1}, $x[\alpha]$ is a Laurent polynomial. 
Hence, by Lemma~\ref{lem:laurent-implies-poly}, 
the second factor in (\ref{eq:intermediate}) 
is a polynomial over $\ZZ\PP$ with nonzero constant term. 
To complete the proof of Theorem~\ref{th:cluster-variable-denominators-1}, 
it remains to compare (\ref{eq:intermediate}) with 
(\ref{eq:laurent-finite-type}), and to observe that 
\[ 
\begin{array}{rclll} 
\gamma 
&=& \tau \Bigl(\sum_{i \neq j} a_{ij} \alpha_i\Bigr) 
&&\text{(by Proposition~\ref{prop:inv-cluster-expansion-2})} 
\\[.1in] 
&=& \tau (\alpha_j \subplus (-\alpha_j)) 
&&\text{(by (\ref{eq:subplus-special}))} 
\\[.05in] 
&=& \tau (\alpha_j) + \tau(-\alpha_j) 
&& \text{(by Lemma~\ref{lem:cluster-expansion-properties-6})} 
\\[.05in] 
&=&  \alpha + \tau (-\alpha_j). 
&&\hspace{1in}\qed 
\end{array} 
\] 
 
\begin{remark} 
Unfortunately, the argument above does not establish 
Theorem~\ref{th:cluster-variable-denominators-2} because 
there is no guarantee that the second factor in (\ref{eq:intermediate}) 
is a polynomial with coefficients in $\ZZ_{\geq 0}[\mathcal{P}]$ 
even if we assume that all the polynomials $P_{\alpha'}$ appearing 
there have this property. 
(Recall from Definition~\ref{def:cluster-algebra} 
that $\mathcal{P}$ denotes the set of all coefficients 
$p^\pm_{\beta,\beta'}$ appearing in various exchange relations 
(\ref{eq:exch-rel-fin-type}); we denote by 
$\mathbb{Z}_{\geq 0}[\mathcal{P}]$ the set of 
polynomials with nonnegative integer coefficients 
in the elements of~$\mathcal{P}$.) 
The proof of Theorem~\ref{th:cluster-variable-denominators-2} 
given in Section~\ref{sec:initial-cluster-positivity} below 
does not rely on Theorem~\ref{th:cluster-variable-denominators-1}, 
thus providing an alternative proof of the latter. 
\end{remark} 
 
\pagebreak[3] 
 
\section{Proof of Theorem~\ref{th:cluster-variable-denominators-2}} 
\label{sec:initial-cluster-positivity} 
 
We use the nomenclature of root systems given in Bourbaki~\cite{bourbaki}, 
including the labeling of the simple roots in $\Phi$ by the indices $1,\dots, n$. 
On the other hand, our convention on associating a Cartan matrix $A$ 
to a root system~$\Phi$, as described in 
Section~\ref{sec:laurent-finite-type}, 
is \emph{transposed} to that in~\cite{bourbaki}---and the same as 
that in Kac~\cite{kac}. 
 
We abbreviate $x_i = x[-\alpha_i]$ for $i = 1, \dots, n$. 
Our goal is to prove that, for every almost positive root $\alpha$, 
we can write $x[\alpha]$ as a Laurent polynomial in $x_1,\dots,x_n$ 
with coefficients in $\mathbb{Z}_{\geq  0}[\mathcal{P}]$. 
This time we will proceed by induction on the height of~$\alpha$ 
(recall that ${\rm ht}(\alpha) = \sum_i [\alpha:\alpha_i]$). 
The base case $\alpha\in-\Pi$ is trivial. 
The induction step will follow from the lemma below. 
 
\begin{lemma} 
\label{lem:exceptional-exchange-with-neg-simple} 
For every positive root $\alpha$, there exists an index $j\in I$ 
such that 
\begin{equation} 
x_j x[\alpha] = F(x[\beta_1],\dots,x[\beta_m]) \,, 
\end{equation} 
where $F$ is a polynomial with coefficients in $\mathbb{Z}_{\geq 
    0}[\mathcal{P}]$ in some cluster variables 
$x[\beta_1],\dots,x[\beta_m]$ such that 
  $\height(\beta_i)<\height(\alpha)$ for all~$i$. 
\end{lemma} 
 
The rest of this section is devoted to the proof of 
Lemma~\ref{lem:exceptional-exchange-with-neg-simple}.

We call a positive root $\alpha$ \emph{non-exceptional} if there exists 
a negative simple root $-\alpha_j$ exchangeable with~$\alpha$; 
otherwise, $\alpha$ will be called \emph{exceptional}. 
If the root $\alpha$ in Lemma~\ref{lem:exceptional-exchange-with-neg-simple} 
is non-exceptional, and $-\alpha_j$ is a negative simple root 
exchangeable with~$\alpha$, then one easily sees that all cluster 
components of the vectors $\alpha-\alpha_j$ and $\alpha\subplus (-\alpha_j)$ 
appearing in the right-hand side of (\ref{eq:special-cluster-expansion}) 
are of smaller height than~$\alpha$, and we are done. 
Thus, it remains to prove 
Lemma~\ref{lem:exceptional-exchange-with-neg-simple} for the 
exceptional roots. First, we identify them explicitly. 
 
\begin{lemma} 
\label{lem:every-root-except-some} 
The complete list of all exceptional positive roots is as follows: 
\begin{align} 
\label{eq:E8-exceptional} 
&\text{$\Phi$ is of type~$E_8\,$, and $\alpha=\alpha_{\rm max}$ is the 
highest root in~$\Phi$;} 
\\ 
\label{eq:F4-exceptional-1} 
&\text{$\Phi$ is of type~$F_4\,$, and 
  $\alpha=\alpha_{\rm max}=2\alpha_1+3\alpha_2+4\alpha_3+2\alpha_4\,$; 
} 
\\ 
\label{eq:F4-exceptional-2} 
&\text{$\Phi$ is of type~$F_4\,$, 
and 
  $\alpha=\alpha_1+2\alpha_2+3\alpha_3+2\alpha_2\,$;} 
\\ 
\label{eq:G2-exceptional-1} 
&\text{$\Phi$ is of type~$G_2\,$, and 
  $\alpha=\alpha_{\rm max}=3\alpha_1+2\alpha_2\,$;} 
\\ 
\label{eq:G2-exceptional-2} 
&\text{$\Phi$ is of type~$G_2\,$, 
and 
$\alpha=2\alpha_1+\alpha_2\,$.} 
\end{align} 
\end{lemma}

\begin{proof} 
As noted in \cite[Remark~1.16]{cfz}, $\alpha$ and $-\alpha_j$ 
are exchangeable if and only if 
\begin{equation} 
\label{eq:coefficient-1} 
[\alpha : \alpha_j] = [\alpha^\vee : \alpha_j^\vee] = 1 \, , 
\end{equation} 
where $\alpha^\vee$ is the coroot corresponding to $\alpha$ under 
the natural bijection between $\Phi$ and the dual system~$\Phi^\vee$. 
Let $(\alpha,\beta)$ denote a $W$-invariant scalar product on the root 
lattice~$Q$. 
Then $[\alpha^\vee : 
\alpha_j^\vee]=\frac{(\alpha_j,\alpha_j)}{(\alpha,\alpha)} 
[\alpha:\alpha_j]$, 
so (\ref{eq:coefficient-1}) is equivalent to 
\begin{equation} 
\label{eq:coefficient-1a} 
[\alpha : \alpha_j] = 1\,,\quad 
(\alpha,\alpha) = (\alpha_j,\alpha_j) \,. 
\end{equation} 
Thus, we need to verify that for every positive root~$\alpha$, 
there exists a simple root $\alpha_j$ satisfying 
(\ref{eq:coefficient-1a}), unless $\alpha$ appears 
on the list (\ref{eq:E8-exceptional})--(\ref{eq:G2-exceptional-2}), 
in which case there is no such simple root. 
This is checked by direct inspection using, e.g., 
the tables in~\cite{bourbaki}. 
(In all classical types, the list of all pairs $(\alpha, -\alpha_j)$ 
satisfying (\ref{eq:coefficient-1}) was given in~\cite{cfz}.) 
\end{proof} 
 
\pagebreak[2] 
 
\begin{figure}[ht] 
\setlength{\unitlength}{1.5pt} 
\begin{picture}(140,28)(0,-20) 
\put(0,0){\line(1,0){120}} 
\put(40,0){\line(0,-1){20}} 
\put(40,-20){\circle*{2}} 
\multiput(0,0)(20,0){7}{\circle*{2}} 
\put(0,4){\makebox(0,0){$1$}} 
\put(44,-20){\makebox(0,0){$2$}} 
\put(20,4){\makebox(0,0){$3$}} 
\put(40,4){\makebox(0,0){$4$}} 
\put(60,4){\makebox(0,0){$5$}} 
\put(80,4){\makebox(0,0){$6$}} 
\put(100,4){\makebox(0,0){$7$}} 
\put(120,4){\makebox(0,0){$8$}} 
\end{picture} 
\begin{picture}(40,24)(10,-20) 
\put(20,0){\line(1,0){40}} 
\put(40,0){\line(0,-1){20}} 
\put(40,-20){\circle*{2}} 
\multiput(20,0)(20,0){3}{\circle*{2}} 
\put(44,-20){\makebox(0,0){$3$}} 
\put(20,4){\makebox(0,0){$1$}} 
\put(40,4){\makebox(0,0){$2$}} 
\put(60,4){\makebox(0,0){$4$}} 
\end{picture} 
\caption{Dynkin diagrams of types $E_8$ and $D_4$} 
\label{fig:dynkin-diagram-E8} 
\end{figure}
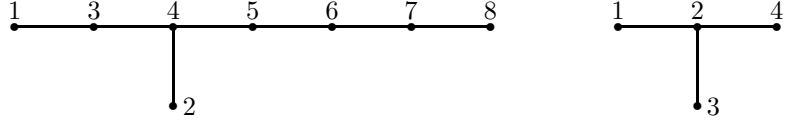

\credit{Proof of Lemma~\ref{lem:exceptional-exchange-with-neg-simple} 
for the type~$E_8$ and $\alpha=\alpha_{\rm max}$} 
This case is by far the hardest among 
(\ref{eq:E8-exceptional})--(\ref{eq:G2-exceptional-2}), so we will 
treat it in detail. 
We will prove that in this special case, 
Lemma~\ref{lem:exceptional-exchange-with-neg-simple} 
holds with $j=8$, in the standard numeration of simple roots 
(see Figure~\ref{fig:dynkin-diagram-E8}). 
 
We will need the following construction. 
In view of 
Lemma~\ref{lem:B-are-tau-invariant}, 
any transformation $\sigma\in\langle\tau_+,\tau_-\rangle$ gives rise 
to a ``twisted" cluster algebra $\sigma (\AA)$ whose seeds are the 
transfers by $\sigma$ of the seeds of~$\AA$; 
if $\sigma$ is written in terms of $\tau_+$ and $\tau_-$ as a product of 
an odd number of factors, this transfer involves the change of signs for 
the matrices $B$ and the corresponding interchange of $p^+$ and $p^-$ for 
the coefficients, as in Remark~\ref{rem:sign-symmetry}. 
This twist preserves the Cartan-Killing type.

Direct computation shows that for $\sigma=(\tau_-\tau_+)^{8}=(\tau_+\tau_-)^{8}$ 
(cf.\ Theorem~\ref{th:dihedral}), we have 
$\sigma(\alpha_{\rm max})=-\alpha_4$ and 
$\sigma(-\alpha_8)=\alpha_2+\alpha_3+2\alpha_4+\alpha_5\,$. 
To prove Lemma~~\ref{lem:exceptional-exchange-with-neg-simple} 
for the type~$E_8$ and $\alpha=\alpha_{\rm max}$, it is 
therefore sufficient to show that, in the twisted cluster algebra 
$\sigma (\AA)$, we have 
\begin{equation} 
\label{eq:x[alpha4]...} 
x_4 x[\alpha_2+\alpha_3+2\alpha_4+\alpha_5] 
=\tilde F(x[\beta_1],\dots,x[\beta_m]) \,, 
\end{equation} 
where $\tilde F$ is a polynomial with coefficients in $\ZZ_{\geq 
  0}[\mathcal{P}]$, 
and each $\beta_i$ is different from $\alpha_2+\alpha_3+2\alpha_4+\alpha_5\,$.

Let $J=\{2,3,4,5\}\subset I$, 
and let $\Phi(J)$ denote the type $D_4$ root subsystem of $\Phi$ 
spanned by the simple roots $\alpha_j$ with~$j\in J$. 
Applying Lemma~\ref{lem:restriction} four times, 
we conclude that the correspondence 
\[ 
C' \mapsto (-\Pi(I-J)) \cup C' 
\] 
identifies the cluster complex $\Delta(\Phi(J))$ 
with the link of $-\Pi(I-J)$ in 
the cluster complex $\Delta(\Phi)$; 
here we use the notation 
\[ 
-\Pi(I-J) = \{-\alpha_i: i \in I - J\}. 
\] 
The exchange graph $\Gamma(J) = \Gamma_{\Delta(\Phi(J))}$ 
is therefore identified with the induced subgraph in the exchange 
graph of $\AA$ whose vertices are all the clusters 
containing~$-\Pi(I-J)$. 
Let $\AA'$ denote the subring in $\AA$ generated by the cluster variables 
$x[\alpha]$ for $\alpha \in \Phi(J)_{\geq -1}$, 
together with the ``coefficients'' in all exchange relations corresponding 
to the edges in~$\Gamma(J)$, where by a ``coefficient" 
we mean the part of a monomial that does \emph{not} involve the 
variables $x[\alpha]$ for $\alpha \in \Phi(J)_{\geq -1}\,$. 
(Thus, each ``coefficient'' is a product of an element of $\mathcal{P}$ 
and a monomial in the variables $x_i$ for $i \in I - J$.) 
By Lemma~\ref{lem:restriction}, Part~2, the ring $\AA'$ 
is a normalized cluster algebra of type~$D_4$ 
(cf.\ \cite[Proposition~2.6]{fz-clust1}). 
The claim (\ref{eq:x[alpha4]...}) now becomes a consequence of the following lemma. 
 
\begin{lemma} 
In 
the case of type $D_4\,$, with the notation as in 
Figure~\ref{fig:dynkin-diagram-E8}, we have 
\[ 
x_2 x[\alpha_{\max}] 
=G(x[\gamma_1],\dots,x[\gamma_k]) \,. 
\] 
where $\alpha_{\max}=\alpha_1+2\alpha_2+\alpha_3+\alpha_4\,$, the 
almost positive roots $\gamma_1,\dots,\gamma_k$ are different from 
$\alpha_{\max}\,$, 
and $G$ is a polynomial with coefficients in 
$\mathbb{Z}_{\geq 0}[\mathcal{P}]$. 
\end{lemma} 
 
We note that the roots $\beta_1,\dots,\beta_m$ appearing in 
(\ref{eq:x[alpha4]...}) are of two kinds: first, the images of 
$\gamma_1,\dots,\gamma_k$ under the embedding $D_4\to E_8\,$, 
and second, (some of) the ``frozen'' roots 
$-\alpha_1,-\alpha_6,-\alpha_7,-\alpha_8\,$. 
 
\begin{proof} 
Figure~\ref{fig:pentagons-D4} shows a fragment of the 
exchange graph in type~$D_4\,$, 
with each vertex~$C$ representing a cluster containing the $4$ roots written 
into the regions adjacent to~$C$. 
The mutual compatibility of the roots in each of these quadruples is 
easily checked from the definitions. 
 
\begin{figure}[ht] 
\begin{center} 
\setlength{\unitlength}{3.4pt} 
\begin{picture}(91,42)(-55,-4) 
  \put(0,6){\line(0,1){20}} 
  \put(19,0){\line(1,-4){1.5}} 
  \put(19,32){\line(1,4){1.5}} 
  \put(31,16){\line(1,0){6}} 
  \qbezier(0,6)(9.5,3)(19,0) 
  \qbezier(0,26)(9.5,29)(19,32) 
  \qbezier(19,0)(25,8)(31,16) 
  \qbezier(19,32)(25,24)(31,16) 
 
  \put(0,6){\circle*{1}} 
  \put(0,26){\circle*{1}} 
  \put(19,0){\circle*{1}} 
  \put(19,32){\circle*{1}} 
  \put(31,16){\circle*{1}} 
 
  \put(0,6){\line(-1,0){20}} 
  \put(0,26){\line(-1,0){20}} 
 
  \put(-20,6){\circle*{1}} 
  \put(-20,26){\circle*{1}} 
  \put(-39,0){\circle*{1}} 
  \put(-39,32){\circle*{1}} 
  \put(-51,16){\circle*{1}} 
 
  \put(-20,6){\line(0,1){20}} 
  \put(-39,0){\line(-1,-4){1.5}} 
  \put(-39,32){\line(-1,4){1.5}} 
  \put(-51,16){\line(-1,0){6}} 
  \qbezier(-20,6)(-29.5,3)(-39,0) 
  \qbezier(-20,26)(-29.5,29)(-39,32) 
  \qbezier(-39,0)(-45,8)(-51,16) 
  \qbezier(-39,32)(-45,24)(-51,16) 
 
\put(30,0){\makebox(0,0){$\alpha_{\max}$}} 
\put(30,32){\makebox(0,0){$\alpha_1+\alpha_2+\alpha_4$}} 
\put(-50,0){\makebox(0,0){$-\alpha_4$}} 
\put(-50,32){\makebox(0,0){$-\alpha_2$}} 
\put(-10,-1){\makebox(0,0){$\alpha_2+\alpha_3$}} 
\put(-10,33){\makebox(0,0){$\alpha_4$}} 
 
\put(14,16){\makebox(0,0) 
{$\begin{array}{c}\alpha_2+\alpha_4\\[.1in] 
\alpha_2+\alpha_3+\alpha_4\end{array}$}} 
 
\put(-34,16){\makebox(0,0) 
{$\begin{array}{c}-\alpha_1\\[.1in] 
\alpha_3\end{array}$}} 
 
\put(-10,16){\makebox(0,0) 
{$\begin{array}{c}-\alpha_1\\[.1in] 
\alpha_2+\alpha_3+\alpha_4\end{array}$}} 
 
 
 
\end{picture} 
\end{center} 
\caption{\hbox{Fragment of the exchange graph in the 
type~$D_4$}} 
\label{fig:pentagons-D4} 
\end{figure}
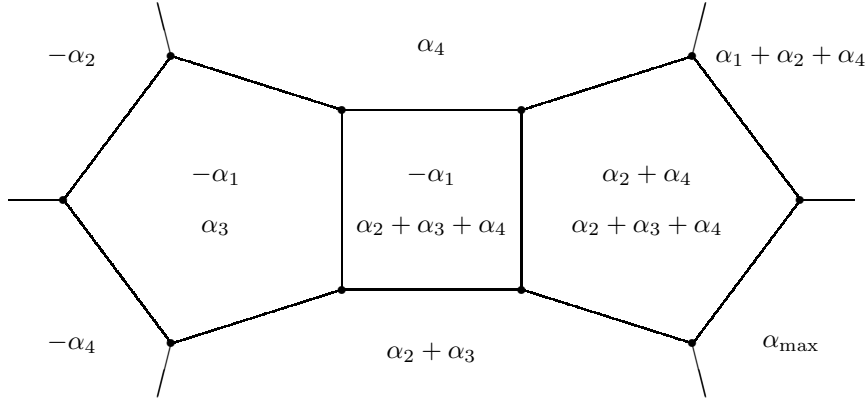 
 
\pagebreak[2] 
 
We next write the exchange relations 
for some pairs of adjacent clusters shown in Figure~\ref{fig:pentagons-D4}. 
In doing so, we use: 
\begin{itemize} 
\item 
(implicitly) the  combinatorial 
interpretation of almost positive roots of type $D_n$ 
given in \cite[Section~3.5]{yga} and reproduced in Section~\ref{sec:ca-type-D} below; 
see specifically \cite[Figure~7]{yga} for the type~$D_4$; 
 
\item 
the resulting explicit expressions for the exchange relations 
which are consequences of \cite[Lemma~4.6]{cfz} 
(see Proposition~\ref{pr:exchanges-D} below); 
 
\item 
the monomial relations among the coefficients of exchange relations 
along a geodesic of type~$A_2$, as given in 
\cite[(6.11)]{fz-clust1}; 
our notation is patterned after \cite[Figure~3]{fz-clust1}. 
\end{itemize} 
The exchange relations for the left pentagonal geodesic 
in Figure~\ref{fig:pentagons-D4} can be written in the following form, 
with $p_1,\dots,p_5\in\mathcal{P}$: 
\begin{align} 
\nonumber 
x_4 x[\alpha_2\!+\!\alpha_3\!+\!\alpha_4] 
&= p_1 x[\alpha_2\!+\!\alpha_3] + p_3 p_4 x[\alpha_3], 
\\[.05in] 
\label{eq:x2x23} 
x_2 x[\alpha_2\!+\!\alpha_3] 
&=p_2 x_1x_4 + p_4 p_5 x[\alpha_3], 
\\[.05in] 
\nonumber 
x_2 x[\alpha_2] 
&=p_3 x_2 + p_5 p_1\,, 
\\[.05in] 
\label{eq:x2x234} 
x_2 x[\alpha_2\!+\!\alpha_3\!+\!\alpha_4] 
&= p_4 x[\alpha_3]x[\alpha_4]+ p_1 p_2 x_1\,, 
\\[.05in] 
\label{eq:4.23-2} 
x[\alpha_4] x[\alpha_2\!+\!\alpha_3] 
&= p_5 x[\alpha_2\!+\!\alpha_3\!+\!\alpha_4] + p_2 p_3 x_1\,.  \hspace{.7in} 
\end{align} 
(Among these relations, only 
(\ref{eq:x2x23}), (\ref{eq:x2x234}), and (\ref{eq:4.23-2}) 
are needed in the proof; 
we wrote all five relations for the sake of clarity.) 
Similarly, the 
exchange relations for the right pentagonal geodesic can 
be written as follows, with $q_1,\dots,q_5\in\mathcal{P}$: 
\begin{align} 
\label{eq:x1x-max} 
x_1 x[\alpha_{\max}] &=q_1  x[\alpha_2\!+\!\alpha_3]x[\alpha_2\!+\!\alpha_4] 
+q_3 q_4 x[\alpha_2\!+\!\alpha_3\!+\!\alpha_4], 
\\[.05in] 
\nonumber 
x[\alpha_2\!+\!\alpha_3]x[\alpha_1\!+\!\alpha_2\!+\!\alpha_4] 
&=q_2 x[\alpha_{\max}] + q_4 q_5\,, 
\\[.05in] 
\nonumber 
x[\alpha_4] x[\alpha_{\max}] 
&= q_3 
x[\alpha_1\!+\!\alpha_2\!+\!\alpha_4]x[\alpha_2\!+\!\alpha_3\!+\!\alpha_4] 
+ q_5 q_1 x[\alpha_2+\alpha_4], 
\\[.05in] 
\label{eq:1-1234} 
x_1 x[\alpha_1\!+\!\alpha_2\!+\!\alpha_4] 
&=q_4 x[\alpha_4] + q_1 q_2 x[\alpha_2\!+\!\alpha_4], 
\\[.05in] 
\label{eq:4.23-1} 
x[\alpha_4] x[\alpha_2\!+\!\alpha_3] 
&= q_5 x_1 + q_2 q_3 x[\alpha_2\!+\!\alpha_3\!+\!\alpha_4]. 
\end{align} 
Comparing (\ref{eq:4.23-1}) to (\ref{eq:4.23-2}), we conclude that 
\begin{equation} 
\label{eq:p5=q2q3} 
p_5=q_2 q_3\,. 
\end{equation} 
Successively applying (\ref{eq:x1x-max}), 
(\ref{eq:x2x23})--(\ref{eq:x2x234}), 
(\ref{eq:p5=q2q3}), and (\ref{eq:1-1234}), we obtain: 
\begin{align*} 
&x_1 x_2 x[\alpha_{\max}] 
\\[.05in] 
=&q_1 x_2 x[\alpha_2\!+\!\alpha_3]x[\alpha_2\!+\!\alpha_4] 
+q_3 q_4 x_2 x[\alpha_2\!+\!\alpha_3\!+\!\alpha_4] 
\\[.05in] 
=& q_1 p_2 x_1 x_4 x[\alpha_2\!+\!\alpha_4] 
+ q_1 p_4 p_5 x[\alpha_3]x[\alpha_2\!+\!\alpha_4] 
+q_3 q_4 p_4 x[\alpha_3]x[\alpha_4] 
+ q_3 q_4 p_1 p_2 x_1 
\\[.05in] 
=& q_1 p_2 x_1 x_4 x[\alpha_2\!+\!\alpha_4] 
+ q_1 q_2 q_3 p_4 x[\alpha_3]x[\alpha_2\!+\!\alpha_4] 
+q_3 q_4 p_4 x[\alpha_3]x[\alpha_4] 
+ q_3 q_4 p_1 p_2 x_1 
\\[.05in] 
=& q_1 p_2 x_1 x_4 x[\alpha_2\!+\!\alpha_4] 
+ q_3 p_4 x[\alpha_3] x_1 x[\alpha_1\!+\!\alpha_2\!+\!\alpha_4] 
+ q_3 q_4 p_1 p_2 x_1 \,, 
\end{align*} 
which implies 
\begin{equation} 
x_2 x[\alpha_{\max}] 
= q_1 p_2 x_4 x[\alpha_2\!+\!\alpha_4] 
+ q_3 p_4 x[\alpha_3] x[\alpha_1\!+\!\alpha_2\!+\!\alpha_4] 
+ q_3 q_4 p_1 p_2 \,, 
\end{equation} 
and we are done. 
\end{proof} 
 

\credit{Proof of Lemma~\ref{lem:exceptional-exchange-with-neg-simple} 
in the types~$F_4$ and~$G_2$} 
One way of handling the non-simply-laced cases is to deduce them from 
the simply-laced ones by means of the ``folding'' technique (see, 
e.g., \cite[Section~2.4]{yga}). 
Alternatively, one can 
perform direct computations, which show that, in the type~$F_4$, we have 
\begin{align} 
\label{eq:F4-1232} 
x_1 x[\alpha_1\!+\!2\alpha_2\!+\!3\alpha_3\!+\!2\alpha_4] 
&=P_1(x[\alpha_3\!+\!\alpha_4], x[\alpha_2\!+\!\alpha_3], \\ 
\nonumber 
&\qquad\quad x[\alpha_2\!+\!2\alpha_3\!+\!2\alpha_4], 
x[\alpha_2\!+\!2\alpha_3\!+\!\alpha_4]), 
\\[.05in] 
\label{eq:F4-xmax} 
 x_4 x[2\alpha_1\!+\!3\alpha_2\!+\!4\alpha_3\!+\!2\alpha_4] 
&= P_2(x[\alpha_1\!+\!2\alpha_2\!+\!3\alpha_3\!+\!\alpha_4], 
x[\alpha_1\!+\!2\alpha_2\!+\!2\alpha_3], \\ 
\nonumber 
&\qquad\quad x[\alpha_1+\alpha_2+\alpha_3+\alpha_4], x[\alpha_1+\alpha_2+\alpha_3]) \,, 
\end{align} 
and,  in the type~$G_2$, we have 
\begin{align} 
x_2 x[2\alpha_1+\alpha_2] 
&=P_3(x[\alpha_1],x[\alpha_1+\alpha_2]), 
\\[.05in] 
x_1 x[3\alpha_1+2\alpha_2] 
&=P_4(x[\alpha_2],x[3\alpha_1+\alpha_2]), 
\end{align} 
where $P_1,P_2,P_3,P_4$ are polynomials with coefficients in 
$\mathbb{Z}_{\geq 0}[\mathcal{P}]$. 
Details are left to the reader. 
 
This completes our proofs of 
Lemma~\ref{lem:exceptional-exchange-with-neg-simple} 
and Theorem~\ref{th:cluster-variable-denominators-2}. 
 
\pagebreak[3] 
 
 
\section{$2$-finite matrices} 
\label{sec:2-finite} 
 
In accordance with the plan outlined in Section~\ref{sec:organization-ca2}, 
our next task is to prove the implication 
${\rm(iii)}\!\Longrightarrow\!{\rm(iv)}$ 
in Theorem~\ref{th:fin-type-characterizations}, 
which will in turn imply 
Theorem~\ref{th:finite-type-class-exhaust}. 
As a first step, we restate the claim at hand as a purely combinatorial result 
(see Theorem~\ref{th:2-fin-CK} below) 
on matrix mutations~(\ref{eq:mutation}). 
 
We shall write $B' = \mu_k (B)$ to denote that a matrix $B'$ is 
obtained from $B$ by a matrix mutation in direction~$k$. 
Note that $\mu_k$ preserves integrality of entries, 
and is an involution: $\mu_k(\mu_k(B))=B$. 
If two matrices 
can be obtained from each other by a 
sequence of matrix mutations followed by a simultaneous permutation of 
rows and columns, we will say that they are \emph{mutation equivalent}. 

A real square matrix $B = (b_{ij})$ is 
\emph{sign-skew-symmetric} (cf.\ (\ref{eq:sss})) 
if, for any $i$ and~$j$, 
either $b_{ij} = b_{ji} = 0$, or else $b_{ij}b_{ji} < 0$; 
in particular, $b_{ii} = 0$ for all~$i$. 
Furthermore, we say that $B$ is \emph{$2$-finite} if it has integer entries, and 
any matrix $B'=(b'_{ij})$ mutation equivalent to 
$B$ is sign-skew-symmetric and satisfies $|b'_{ij} b'_{ji}| \leq 3$ for 
all $i$ and~$j$. 
 
 
 
 
In the language just introduced, 
the implication 
${\rm(iii)}\!\Longrightarrow\!{\rm(iv)}$ 
in Theorem~\ref{th:fin-type-characterizations} 
can be formulated as follows. 
 
\begin{theorem} 
\label{th:2-fin-CK} 
Every $2$-finite matrix $B$ is mutation equivalent to a matrix $B_\circ$ 
from Theorem~~\ref{th:finite-type-class-new}. 
\end{theorem} 
 
 
The converse of Theorem~\ref{th:2-fin-CK} also holds: 
by Theorem~\ref{th:finite-type-class-new} 
(which has already been proved), $B_\circ$ 
is $2$-finite. 
 
Our proof of Theorem~\ref{th:2-fin-CK} occupies the rest of Sections 
\ref{sec:2-finite}--\ref{sec:diagram-CK-proof} below. 
The main result of Section~\ref{sec:2-finite} is the following proposition. 
 
\begin{proposition} 
\label{pr:2-finite-is-ss} 
Every $2$-finite matrix is skew-symmetrizable. 
\end{proposition} 
 
(Recall that a square matrix $B$ is skew-symmetrizable if there exists 
a diagonal matrix $D$ with positive diagonal entries such that $DB$ is 
skew-symmetric.) 
 
The rest of this section is devoted to the proof of 
Proposition~\ref{pr:2-finite-is-ss}. 
 
The crucial role in the sequel will be played by a combinatorial 
construction that associates with a sign-skew-symmetric matrix 
its \emph{diagram}, whose role is parallel to that of the Dynkin 
diagram 
for a generalized Cartan matrix. 
 
\begin{definition} 
\label{def:diagram-of-B} 
The \emph{diagram} of a sign-skew-symmetric 
matrix~$B=(b_{ij})_{i,j\in I}$ is the weighted directed graph $\Gamma (B)$ 
with the vertex set $I$ such that there is a directed edge from $i$ to $j$ 
if and only if $b_{ij} > 0$, and this edge is assigned the weight 
$|b_{ij}b_{ji}|\,$. 
\end{definition} 
 
More generally, we will use the term \emph{diagram} to denote a finite 
directed graph without loops and multiple edges, whose edges are 
assigned positive real weights. 
By some abuse of notation, we denote by the same symbol $\Gamma$ 
the underlying directed graph of a diagram. 
If two vertices of $\Gamma$ are not joined by an edge, 
we may also say that they are joined by an edge of weight~$0$.

The following lemma is an analogue of the well-known 
symmetrizability criterion \cite[Exercise~2.1]{kac}. 
 
\begin{lemma} 
A matrix $B = (b_{ij})$ is skew-symmetrizable if and only if, first, 
it is sign-skew-symmetric and, second, for all $k \geq 3$ and all 
$i_1, \dots, i_k\,$, it satisfies 
\begin{equation} 
\label{eq:cycle=cycle} 
b_{i_1 i_2} b_{i_2 i_3} \cdots b_{i_k i_1} = 
(-1)^k b_{i_2 i_1} b_{i_3 i_2} \cdots b_{i_1 i_k}\,. 
\end{equation} 
\end{lemma}

\proof 
The ``only if" part 
is trivial. 
Thus, let us assume that $B$ is sign-skew-symmetric and satisfies 
(\ref{eq:cycle=cycle}). 
Without loss of generality, we also assume that $B$ is 
indecomposable, i.e., cannot be represented as a direct 
(block-diagonal) sum of two proper submatrices. 
It follows that the graph $\Gamma (B)$ is connected. 
Let $T$ be one of its spanning trees. 
There exists a diagonal matrix $D=(d_{ij})$ with positive diagonal entries such that 
$d_{ii} b_{ij}=-d_{jj} b_{ji}$ for every edge $(i,j)$ in~$T$. 
(Such a matrix can be constructed inductively by setting $d_{ii}$ 
equal to an arbitrary positive number for some vertex $i$, and moving 
within the tree~$T$ away from this vertex.) 
Then $DB$ is skew-symmetric, for the following reason: 
by definition of a spanning tree, any edge $(i,j)$ of $\Gamma (B)$ 
which is not in $T$  belongs to a cycle in which the rest of the edges belong to~$T$; then use 
(\ref{eq:cycle=cycle}). 
\endproof

\begin{lemma} 
\label{lem:2-finite-3-cycle} 
Let $B$ be a $2$-finite matrix. 
Then the edges of every triangle in $\Gamma(B)$ are 
oriented in a cyclic way. 
\end{lemma} 
 
\proof 
Suppose on the contrary that 
$b_{ij}, b_{ik}, b_{kj} > 0$ for some distinct $i,j,k$. 
Then in the matrix $B' = \mu_k (B)$, we have 
$b'_{ij} = b_{ij} + b_{ik} b_{kj} \geq 2$ and 
$b'_{ji} = b_{ji} - b_{jk} b_{ki} \leq -2$, 
violating $2$-finiteness. 
\endproof 
 
\begin{lemma} 
\label{lem:2-finite-3-ss} 
Let $B$ be a $2$-finite matrix. Then 
\begin{equation} 
\label{eq:3by3-ss-ijk} 
b_{ij} b_{jk} b_{ki} = -b_{ji} b_{kj} b_{ik} 
\end{equation} 
for any distinct $i,j,k$. 
Also, in every triangle in $\Gamma(B)$, the edge weights are either 
$\{1,1,1\}$ or $\{2,2,1\}$. 
\end{lemma} 
 
\proof 
In view of Lemma~\ref{lem:2-finite-3-cycle}, 
we may assume without loss of generality that 
$B$ is a $3 \times 3$ matrix 
\begin{equation} 
\label{eq:3by3matrix} 
\left[\begin{array}{ccc} 
0   & a_1 & - c_2 \\ 
-a_2 & 0   & b_1 \\ 
c_1 & -b_2 & 0 \\ 
\end{array}\right] \,, 
\end{equation} 
where $a_1,b_1,c_1,a_2,b_2,c_2$ are positive integers. 
(If one of these entries is~$0$, then (\ref{eq:3by3-ss-ijk}) is 
automatic.) 
Again without loss of generality, we may assume that 
the entry of maximal absolute value in $B$ is~$-c_2$. 
We claim that, under this assumption, 
\begin{equation} 
\label{eq:strong-3-cyclic} 
c_1 = a_2 b_2, \,\, c_2 = a_1 b_1 \, , 
\end{equation} 
implying $a_1 b_1 c_1 = a_2 b_2 c_2$, 
and hence proving (\ref{eq:3by3-ss-ijk}). 
 
Indeed, we have 
\begin{equation} 
\label{eq:3by3matrix-mutated} 
\mu_2(B)= 
\left[\begin{array}{ccc} 
0   & -a_1 & a_1 b_1 -c_2 \\ 
a_2 & 0   & -b_1 \\ 
-a_2 b_2+ c_1 & b_2 & 0 \\ 
\end{array}\right] \,. 
\end{equation} 
Applying Lemma~\ref{lem:2-finite-3-cycle} to 
$\mu_2(B)$, we conclude that 
$$a_1 b_1 -c_2 \geq 0, \,\, a_2 b_2 -c_1 \geq 0 ,$$ 
where either both inequalities are strict, or both are equalities. 
We need to show that the former case is impossible. 
Indeed, otherwise we would have had 
$a_2 b_2 > c_1 \geq 1$, implying 
$\max (a_2, b_2) \geq 2$; 
also, $a_1 b_1 > c_2 \geq \max (a_1, b_1)$, implying 
$a_1 \geq 2$ and $b_1 \geq 2$. 
But then 
$\max(a_1 a_2, b_1 b_2) 
\geq 4$, 
contradicting the $2$-finiteness of~$B$. 
 
It remains to show that the set of edge weights 
$\{a_1a_2,b_1b_2,c_1c_2\}$ 
is either $\{1,1,1\}$, or $\{2,2,1\}$. 
The only other option consistent with both (\ref{eq:strong-3-cyclic}) and 
the inequalities $a_1a_2\leq 3,b_1b_2\leq 3,c_1c_2\leq 3$ is 
$c_2=3$, $\{a_1,b_1\}=\{3,1\}$, $c_1=a_2=b_2=1$. 
Say $a_1=3$ and $b_1=1$ (the other case is analogous). 
Then $B'=\mu_1(B)$ has $|b'_{23}b'_{32}|=4$, 
violating $2$-finiteness. 
\endproof

Now everything is ready for the proof of 
Proposition~\ref{pr:2-finite-is-ss}. 
It suffices to check that every $2$-finite matrix satisfies the criterion (\ref{eq:cycle=cycle}). 
Suppose this is not the case. 
Among all instances where (\ref{eq:cycle=cycle}) is violated for some 
$2$-finite matrix $B$, pick one with the smallest value of~$k$. 
Then $b_{i_j, i_m}=0$ for any pair of subscripts $(i_j,i_m)$ not 
appearing in (\ref{eq:cycle=cycle}). 
(Otherwise we could obtain (\ref{eq:cycle=cycle}) as a 
corollary of its counterparts for two smaller cycles.) 
In other words, the diagram $\Gamma(B)$ restricted to the vertices 
$i_1,\dots,i_k$ must be a cycle. 
Pick any two consecutive edges on this cycle that form an oriented 
$2$-path (that is, $b_{i_{j-1} i_j} b_{i_j i_{j+1}}>0$). 
(If there is no such pair, we will need to first apply a mutation at an 
arbitrary vertex~$i_j\,$.) 
By Lemma~\ref{lem:2-finite-3-ss}, we have $k \geq 4$, hence 
$b_{i_{j-1} i_{j+1}} = 0$. 
Now apply the mutation~$\mu_{i_j}\,$. 
In the resulting matrix, condition (\ref{eq:cycle=cycle}) for the 
sequence of indices $i_1,\dots,i_{j-1},i_{j+1},\dots,i_k$ will be 
equivalent to (\ref{eq:cycle=cycle}) in the original matrix; 
hence it must fail, contradicting our choice of~$k$. 
\endproof

\section{Diagram mutations} 
\label{sec:diagram-mutations} 
 
Let $B = (b_{ij})_{i,j \in I}$ be a skew-symmetrizable matrix. 
Notice that the diagram $\Gamma (B)$ does \emph{not} determine $B$: for 
instance, the matrix $(-B^T)$ has the same diagram as $B$. 
However, the following important property holds. 
 
\begin{figure}[ht] 
\setlength{\unitlength}{1.5pt} 
\begin{picture}(30,17)(-5,0) 
\put(0,0){\line(1,0){20}} 
\put(0,0){\line(2,3){10}} 
\put(0,0){\vector(2,3){6}} 
\put(10,15){\line(2,-3){10}} 
\put(10,15){\vector(2,-3){6}} 
\put(0,0){\circle*{2}} 
\put(20,0){\circle*{2}} 
\put(10,15){\circle*{2}} 
\put(2,10){\makebox(0,0){$a$}} 
\put(18,10){\makebox(0,0){$b$}} 
\put(10,-4){\makebox(0,0){$c$}} 
\put(10,19){\makebox(0,0){$k$}} 
\end{picture} 
$ 
\begin{array}{c} 
\stackrel{\textstyle\mu_k}{\longleftrightarrow} 
\\[.3in] 
\end{array} 
$ 
\setlength{\unitlength}{1.5pt} 
\begin{picture}(30,17)(-5,0) 
\put(0,0){\line(1,0){20}} 
\put(0,0){\line(2,3){10}} 
\put(10,15){\vector(-2,-3){6}} 
\put(10,15){\line(2,-3){10}} 
\put(20,0){\vector(-2,3){6}} 
\put(0,0){\circle*{2}} 
\put(20,0){\circle*{2}} 
\put(10,15){\circle*{2}} 
\put(2,10){\makebox(0,0){$a$}} 
\put(18,10){\makebox(0,0){$b$}} 
\put(10,-4){\makebox(0,0){$c'$}} 
\put(10,19){\makebox(0,0){$k$}} 
\end{picture} 
 
\vspace{-.2in} 
\caption{Diagram mutation} 
\label{fig:diagram-mutation-general} 
\end{figure}
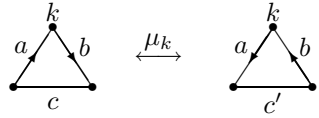

\begin{proposition} 
\label{pr:diagram-mutation} 
For a skew-symmetrizable matrix $B$, 
the diagram $\Gamma'\!\!=\!\Gamma(\mu_k(B))$ is uniquely determined by 
the diagram $\Gamma=\Gamma(B)$ and an index $k\in I$. 
Specifically, $\Gamma'$ is obtained from $\Gamma$ as follows: 
\begin{itemize} 
\item The orientations of all edges incident to~$k$ are reversed, 
their weights intact. 
\item 
For any vertices $i$ and $j$ which are connected in 
$\Gamma$ via a two-edge oriented path going through~$k$ (refer to 
Figure~\ref{fig:diagram-mutation-general} for the rest of notation), 
the direction of the edge $(i,j)$ in $\Gamma'$ and its weight $c'$ 
are uniquely determined by the rule 
\begin{equation} 
\label{eq:weight-relation-general} 
\pm\sqrt {c} \pm\sqrt {c'} = \sqrt {ab} \,, 
\end{equation} 
where the sign before $\sqrt {c}$ 
(resp., before $\sqrt {c'}$) 
is ``$+$'' if $i,j,k$ form an oriented cycle 
in~$\Gamma$ (resp., in~$\Gamma'$), and is ``$-$'' otherwise. 
Here either $c$ or $c'$ can be equal to~$0$. 
 
\item 
The rest of the edges and their weights in $\Gamma$ 
remain unchanged. 
\end{itemize} 
\end{proposition} 
 
\begin{remark} 
If $B$ has integer entries, then all edge weights in $\Gamma$ 
are positive integers. 
The rule (\ref{eq:weight-relation-general}) ensures that the same 
is true for $\Gamma'$: indeed, the fact that 
$c'=(\sqrt {ab}\mp\sqrt {c})^2=ab+c\mp 2\sqrt{abc}$ 
is an integer (that is, $abc$ is a perfect square) 
is an easy consequence of the skew-symmetrizability of~$B$ 
(more specifically, of the identity (\ref{eq:cycle=cycle}) with $k=3$). 
\end{remark} 
 
Our proof of Proposition~\ref{pr:diagram-mutation} is based on the 
following construction. 
 
\pagebreak[3] 
 
\begin{lemma} 
\label{lem:ss-conjugation} 
Let $B$ be a skew-symmetrizable matrix. 
Then there exists a diagonal matrix $H$ with 
positive diagonal entries such that 
$HBH^{-1}$ is skew-symmetric. 
Furthermore, the matrix $S(B)=(s_{ij})=HBH^{-1}$ is 
uniquely determined by~$B$. 
Specifically, the matrix entries of $S(B)$ are given by 
\begin{equation} 
\label{eq:s-matrix} 
s_{ij} = 
{\rm sgn}(b_{ij}) \textstyle\sqrt{|b_{ij} b_{ji}|}\,. 
\end{equation} 
\end{lemma} 
 
\begin{proof} 
Let $D$ be a skew-symmetrizing matrix for $B$, i.e., 
a diagonal matrix with positive diagonal entries such that $DB$ is 
skew-symmetric. 
Setting $H =  D^{1/2}$, we see that 
$HBH^{-1} = H^{- 1} (DB) H^{-1}$ 
is skew-symmetric. 
To prove~(\ref{eq:s-matrix}), note that 
\begin{align*} 
{\rm sgn} (s_{ij}) &= {\rm sgn} (h_i b_{ij} h_j^{-1}) = {\rm sgn} 
(b_{ij}),\\ 
s_{ij}^2  = |s_{ij} s_{ji}| &= 
|(h_i b_{ij} h_j^{-1}) \cdot (h_j b_{ji} h_i^{-1})| = 
 |b_{ij} b_{ji}|, 
\end{align*} 
where the $h_i$ are the diagonal entries of~$H$. 
\end{proof}

\begin{lemma} 
\label{lem:S(t)} 
Let $B$ be a skew-symmetrizable matrix. 
Then, for any $k\in I$, we have $S(\mu_k(B))=\mu_k(S(B))$. 
\end{lemma}

\proof 
Follows from Lemma~\ref{lem:ss-conjugation}, together with the 
directly checked fact that the mutation rules are invariant under 
conjugation by a diagonal matrix with positive entries. 
\endproof 
 
\credit{Proof of Proposition~\ref{pr:diagram-mutation}} 
Formula (\ref{eq:s-matrix}) shows that the diagram $\Gamma(B)$ and the 
skew-symmetric matrix $S(B)$ encode the same information about a 
skew-symmetr\-izable matrix~$B$: 
having an edge in $\Gamma (B)$ directed from $i$ to $j$ and supplied 
with weight $c$ is the same as saying that $s_{ij} = \sqrt {c}$ and 
$s_{ji} = - \sqrt {c}$. 
Lemma~\ref{lem:S(t)} asserts that, as $B$ undergoes a mutation 
$\mu_k\,$, so does the matrix~$S(B)$. 
Translating this statement into the language of diagrams, 
we obtain Proposition~\ref{pr:diagram-mutation}. 
\endproof

In the situation of Proposition~\ref{pr:diagram-mutation}, 
we write $\Gamma' = \mu_k (\Gamma)$, and call the transformation 
$\mu_k$ a \emph{diagram mutation} in the direction~$k$. 
Two diagrams $\Gamma$ and $\Gamma'$ related by a sequence of diagram mutations are called 
\emph{mutation equivalent,} and we write $\Gamma\sim\Gamma'$. 
A diagram $\Gamma$ is called 
\emph{$2$-finite} if any diagram $\Gamma'\sim\Gamma$ 
has all edge weights equal to $1,2$ or~$3$. 
Thus a matrix $B$ is $2$-finite if and only if its diagram $\Gamma(B)$ 
is $2$-finite. 
(Here we rely on Proposition~\ref{pr:2-finite-is-ss}.) 
Note that a diagram is $2$-finite if and only if so are all its 
connected components.


In the case of $2$-finite diagrams, 
Lemmas~\ref{lem:2-finite-3-cycle} and \ref{lem:2-finite-3-ss} 
ensure that every triangle 
is oriented in a cyclic way, and has edge weights $(1,1,1)$ or $(2,2,1)$. 
As a result, the rules of diagram mutations (as given in 
Proposition~\ref{pr:diagram-mutation}) simplify as follows. 
 
\vbox{ 
\begin{lemma} 
\label{lem:diagram-mutation-2-finite} 
Let $\Gamma$ be a $2$-finite diagram, and $k$ a vertex of~$\Gamma$. 
Then the diagram $\mu_k (\Gamma)$ is obtained from $\Gamma$ as follows: 
\begin{itemize} 
\item The orientations of all edges incident to~$k$ are reversed, 
their weights intact. 
\item 
For any vertices $i$ and $j$ which are connected in 
$\Gamma$ via a two-edge oriented path going through~$k$, 
the diagram mutation $\mu_k$ affects the edge connecting $i$ and $j$ in the 
way shown in 
Figure~\ref{fig:diagram-mutation}, where the weights $c$ and $c'$ are 
related by 
\begin{equation} 
\label{eq:weight-relation} 
\sqrt {c} + \sqrt {c'} = \sqrt {ab} \ ; 
\end{equation} 
here either $c$ or $c'$ can be equal to $0$. 
\item 
The rest of the edges and their weights in $\Gamma$ 
remain unchanged. 
\end{itemize} 
\end{lemma} 
} 
 
\pagebreak[3] 
 
\begin{center} 
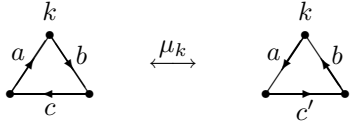
\begin{figure}[ht] 
\setlength{\unitlength}{1.5pt} 
\begin{picture}(30,23)(-5,0) 
\put(0,0){\line(1,0){20}} 
\put(20,0){\vector(-1,0){12}} 
\put(0,0){\line(2,3){10}} 
\put(0,0){\vector(2,3){6}} 
\put(10,15){\line(2,-3){10}} 
\put(10,15){\vector(2,-3){6}} 
\put(0,0){\circle*{2}} 
\put(20,0){\circle*{2}} 
\put(10,15){\circle*{2}} 
\put(2,10){\makebox(0,0){$a$}} 
\put(18,10){\makebox(0,0){$b$}} 
\put(10,-4){\makebox(0,0){$c$}} 
\put(10,21){\makebox(0,0){$k$}} 
\end{picture} 
\ \ 
$\begin{array}{c} 
\stackrel{\textstyle\mu_k}{\longleftrightarrow} 
\\[.3in] 
\end{array} 
$ 
\ \ 
\setlength{\unitlength}{1.5pt} 
\begin{picture}(30,23)(-5,0) 
\put(0,0){\line(1,0){20}} 
\put(0,0){\vector(1,0){12}} 
\put(0,0){\line(2,3){10}} 
\put(10,15){\vector(-2,-3){6}} 
\put(10,15){\line(2,-3){10}} 
\put(20,0){\vector(-2,3){6}} 
\put(0,0){\circle*{2}} 
\put(20,0){\circle*{2}} 
\put(10,15){\circle*{2}} 
\put(2,10){\makebox(0,0){$a$}} 
\put(18,10){\makebox(0,0){$b$}} 
\put(10,-4){\makebox(0,0){$c'$}} 
\put(10,21){\makebox(0,0){$k$}} 
\end{picture} 
\vspace{-.1in} 
\caption{Mutation of $2$-finite diagrams} 
\label{fig:diagram-mutation} 
\end{figure} 
\end{center} 
\vspace{-.2in} 
 
Taking into account Propositions~\ref{pr:2-finite-is-ss} and~\ref{pr:diagram-mutation}, 
we see that Theorem~\ref{th:2-fin-CK} becomes a consequence of the 
following classification of $2$-finite diagrams.

\begin{figure}[ht] 
\vspace{-.2in} 
\[ 
\begin{array}{ccl} 
A_n && 
\setlength{\unitlength}{1.5pt} 
\begin{picture}(140,17)(0,-2) 
\put(0,0){\line(1,0){140}} 
\multiput(0,0)(20,0){8}{\circle*{2}} 
\end{picture}\\
B_n
&& 
\setlength{\unitlength}{1.5pt} 
\begin{picture}(140,17)(0,-2) 
\put(0,0){\line(1,0){140}} 
\multiput(0,0)(20,0){8}{\circle*{2}} 
\put(10,4){\makebox(0,0){$2$}} 
\end{picture} 
\\[.2in] 
D_n 
&& 
\setlength{\unitlength}{1.5pt} 
\begin{picture}(140,17)(0,-2) 
\put(20,0){\line(1,0){120}} 
\put(0,10){\line(2,-1){20}} 
\put(0,-10){\line(2,1){20}} 
\multiput(20,0)(20,0){7}{\circle*{2}} 
\put(0,10){\circle*{2}} 
\put(0,-10){\circle*{2}} 
\end{picture} 
\\[.2in] 
E_6 
&& 
\setlength{\unitlength}{1.5pt} 
\begin{picture}(140,17)(0,-2) 
\put(0,0){\line(1,0){80}} 
\put(40,0){\line(0,-1){20}} 
\put(40,-20){\circle*{2}} 
\multiput(0,0)(20,0){5}{\circle*{2}} 
\end{picture} 
\\[.25in] 
E_7 
&& 
\setlength{\unitlength}{1.5pt} 
\begin{picture}(140,17)(0,-2) 
\put(0,0){\line(1,0){100}} 
\put(40,0){\line(0,-1){20}} 
\put(40,-20){\circle*{2}} 
\multiput(0,0)(20,0){6}{\circle*{2}} 
\end{picture} 
\\[.25in] 
E_8 
&& 
\setlength{\unitlength}{1.5pt} 
\begin{picture}(140,17)(0,-2) 
\put(0,0){\line(1,0){120}} 
\put(40,0){\line(0,-1){20}} 
\put(40,-20){\circle*{2}} 
\multiput(0,0)(20,0){7}{\circle*{2}} 
\end{picture} 
\\[.3in] 
F_4 
&& 
\setlength{\unitlength}{1.5pt} 
\begin{picture}(140,17)(0,-2) 
\put(0,0){\line(1,0){60}} 
\multiput(0,0)(20,0){4}{\circle*{2}} 
\put(30,4){\makebox(0,0){$2$}} 
\end{picture} 
\\[.1in] 
G_2 
&& 
\setlength{\unitlength}{1.5pt} 
\begin{picture}(140,17)(0,-2) 
\put(0,0){\line(1,0){20}} 
\multiput(0,0)(20,0){2}{\circle*{2}} 
\put(10,4){\makebox(0,0){$3$}} 
\end{picture} 
\end{array} 
\] 
\vspace{-.1in} 
\caption{Dynkin diagrams} 
\label{fig:dynkin-diagrams} 
\end{figure}
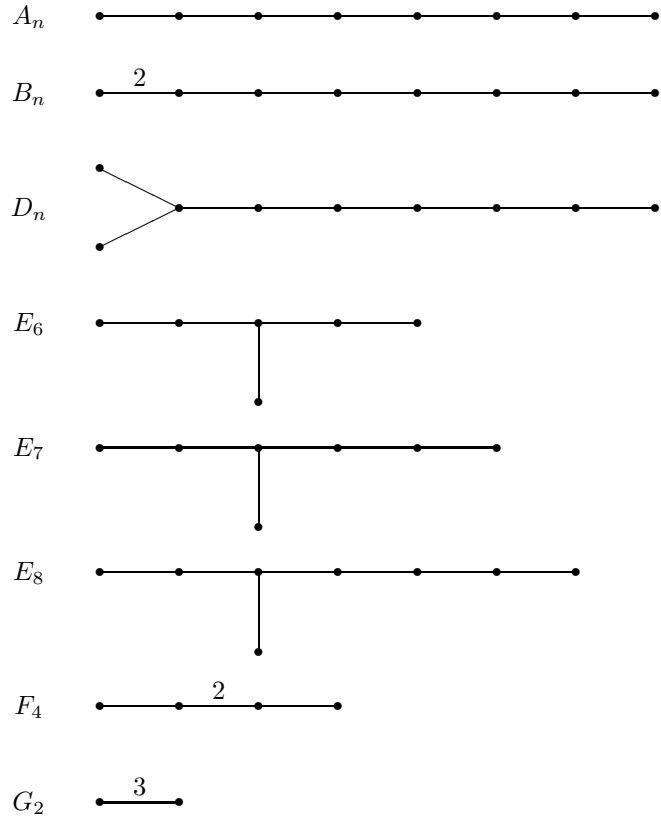

\begin{theorem} 
\label{th:diagrams-fin-CK} 
Any connected $2$-finite diagram is mutation equivalent to 
an orientation of a Dynkin diagram. 
(Cf.\ Figure~\ref{fig:dynkin-diagrams}, where all unspecified 
weights are equal to~$1$.) 
Furthermore, all orientations of the same Dynkin diagram are 
mutation equivalent to each other. 
\end{theorem} 
 
As already noted following Theorem~\ref{th:2-fin-CK}, 
the converse is true as well: any diagram 
mutation equivalent to an orientation of a Dynkin diagram 
is $2$-finite. 
 
\pagebreak[3]

\section{Proof of Theorem~\ref{th:diagrams-fin-CK}} 
\label{sec:diagram-CK-proof} 
 
Throughout this section, all diagrams are presumed connected, and 
all edge weights are positive integers. 
With some abuse of notation, we use the same symbol $\Gamma$ to denote 
a diagram and the set of its vertices. 
A diagram that is not $2$-finite will be called \emph{$2$-infinite.} 
 
\begin{definition} 
A \emph{subdiagram} of a diagram $\Gamma$ is a diagram $\Gamma'$ 
obtained from $\Gamma$ by taking an induced directed subgraph on a subset of vertices and 
keeping all its edge weights the same as in $\Gamma$. 
We will denote this by $\Gamma\supset \Gamma'$. 
\end{definition} 
 
We will repeatedly use the following obvious fact: 
any subdiagram of a $2$-finite diagram is $2$-finite. 
Equivalently, any diagram that has a $2$-infinite subdiagram is 
$2$-infinite.

The proof of Theorem~\ref{th:diagrams-fin-CK} will proceed in 
several steps.

\subsection{Shape-preserving diagram mutations} 
Let $k$ be a sink (resp., source) of a diagram $\Gamma$, 
that is, a vertex such that all edges incident to $k$ 
are directed towards~$k$ (resp., away from~$k$). 
Then a diagram mutation at~$k$ reverses the orientations of 
all edges incident to~$k$, leaving the rest of the graph and 
all the edge weights unchanged. 
We shall refer to such mutations as \emph{shape-preserving}. 
 
\begin{proposition} 
\label{pr:tree-orientations} 
Let $T$ be a subdiagram of a diagram $\Gamma$ such that: 
\begin{itemize} 
\item[{\rm (i)}] $T$ is a tree. 
 
\item[{\rm (ii)}] $T$ is attached to the rest of $\Gamma$ by a single 
vertex $v \in T$, i.e., no vertex in $T - \{v\}$ is joined by an 
edge with a vertex in $\Gamma - T$. 
\end{itemize} 
Then any diagram obtained from $\Gamma$ by arbitrarily re-orienting 
the edges of $T$ (while keeping the rest of $\Gamma$ intact) 
is mutation equivalent to~$\Gamma$. 
 
In particular, any two orientations of a tree diagram are 
mutation equivalent. 
\end{proposition} 
 
(A \emph{tree diagram} is a diagram whose underlying graph is an 
orientation of a tree.) 
 
\proof 
Using induction on the size of $T$, we will show that one can 
arbitrarily re-orient the edges of $T$ by applying a sequence of 
shape-preserving mutations at the vertices of $T - \{v\}$. 
If $T$ consists of a single vertex $v$, there is nothing to prove. 
Otherwise, pick a leaf $l \in T$ different from $v$, and apply the 
inductive assumption to the diagram $\Gamma' = \Gamma - \{l\}$ and 
its subdiagram $T' = T - \{l\}$. 
So we are able to arbitrarily re-orient the edges of $T'$ 
by a sequence of 
shape-preserving mutations of $\Gamma'$ at the vertices of $T' - \{v\}$. 
To do the same for $T$, we lift this sequence from $\Gamma'$ to $\Gamma$ 
as follows: each time right before we need to perform a mutation 
at the unique vertex $k \in T'$ adjacent to $l$, 
we first mutate at $l$ if necessary to 
make $k$ a source or sink in~$T$, rather than just in $T'$. 
This way, we can achieve an arbitrary re-orientation of 
the edges of $T'$ by a sequence of 
shape-preserving mutations of $\Gamma$ at 
the vertices of $T - \{v\}$. 
The remaining edge $(k,l)$ can then be given an arbitrary 
orientation by a (shape-preserving) mutation at~$l$. 
\endproof 
 
As a practical consequence of Proposition~\ref{pr:tree-orientations}, 
in drawing a diagram $\Gamma$, we do not have to specify 
orientations of edges in any subdiagram $T\subset \Gamma$ 
satisfying the conditions  of the proposition. 
Figure~\ref{fig:dynkin-diagrams} provides an example of this 
(see also Figure~\ref{fig:Spqrs} below). 
 
Proposition~\ref{pr:tree-orientations} also justifies notation of the 
form $\Gamma\sim A_m\,$, $\Gamma\supset A_m\,$, etc.

\subsection{Taking care of the trees} 
 
\begin{proposition} 
\label{pr:trees-Dynkin-vs-infty} 
Any $2$-finite tree diagram is an orientation of a Dynkin diagram. 
\end{proposition} 
 
\begin{figure}[ht] 
\[ 
\begin{array}{ccl} 
B_n^{(1)} 
&& 
\setlength{\unitlength}{1.5pt} 
\begin{picture}(140,15)(0,-2) 
\put(20,0){\line(1,0){120}} 
\put(0,10){\line(2,-1){20}} 
\put(0,-10){\line(2,1){20}} 
\multiput(20,0)(20,0){7}{\circle*{2}} 
\put(0,10){\circle*{2}} 
\put(0,-10){\circle*{2}} 
\put(130,4){\makebox(0,0){$2$}} 
\end{picture} 
\\[.2in] 
C_n^{(1)}
&& 
\setlength{\unitlength}{1.5pt} 
\begin{picture}(140,17)(0,-2) 
\put(0,0){\line(1,0){140}} 
\multiput(0,0)(20,0){8}{\circle*{2}} 
\put(10,4){\makebox(0,0){$2$}} 
\put(130,4){\makebox(0,0){$2$}} 
\end{picture} 
\\[.1in] 
D_n^{(1)} 
&& 
\setlength{\unitlength}{1.5pt} 
\begin{picture}(140,17)(0,-2) 
\put(20,0){\line(1,0){100}} 
\put(0,10){\line(2,-1){20}} 
\put(0,-10){\line(2,1){20}} 
\put(120,0){\line(2,-1){20}} 
\put(120,0){\line(2,1){20}} 
\multiput(20,0)(20,0){6}{\circle*{2}} 
\put(0,10){\circle*{2}} 
\put(0,-10){\circle*{2}} 
\put(140,10){\circle*{2}} 
\put(140,-10){\circle*{2}} 
\end{picture} 
\\[.2in] 
E_6^{(1)} 
&& 
\setlength{\unitlength}{1.5pt} 
\begin{picture}(140,17)(0,-2) 
\put(0,0){\line(1,0){80}} 
\put(40,0){\line(0,-1){40}} 
\put(40,-20){\circle*{2}} 
\put(40,-40){\circle*{2}} 
\multiput(0,0)(20,0){5}{\circle*{2}} 
\end{picture} 
\\[.7in] 
E_7^{(1)} 
&& 
\setlength{\unitlength}{1.5pt} 
\begin{picture}(140,17)(0,-2) 
\put(0,0){\line(1,0){120}} 
\put(60,0){\line(0,-1){20}} 
\put(60,-20){\circle*{2}} 
\multiput(0,0)(20,0){7}{\circle*{2}} 
\end{picture} 
\\[.25in] 
E_8^{(1)} 
&& 
\setlength{\unitlength}{1.5pt} 
\begin{picture}(140,17)(0,-2) 
\put(0,0){\line(1,0){140}} 
\put(40,0){\line(0,-1){20}} 
\put(40,-20){\circle*{2}} 
\multiput(0,0)(20,0){8}{\circle*{2}} 
\end{picture} 
\\[.3in] 
F_4^{(1)} 
&& 
\setlength{\unitlength}{1.5pt} 
\begin{picture}(140,17)(0,-2) 
\put(0,0){\line(1,0){80}} 
\multiput(0,0)(20,0){5}{\circle*{2}} 
\put(30,4){\makebox(0,0){$2$}} 
\end{picture} 
\\[.1in] 
G_2^{(1)} 
&& 
\setlength{\unitlength}{1.5pt} 
\begin{picture}(140,17)(0,-2) 
\put(0,0){\line(1,0){40}} 
\multiput(0,0)(20,0){3}{\circle*{2}} 
\put(10,4){\makebox(0,0){$3$}} 
\put(30,4){\makebox(0,0){$a$}} 
\end{picture} 
\end{array} 
\] 
\caption{Extended Dynkin tree diagrams} 
\label{fig:extended-dynkin-diagrams} 
\end{figure}
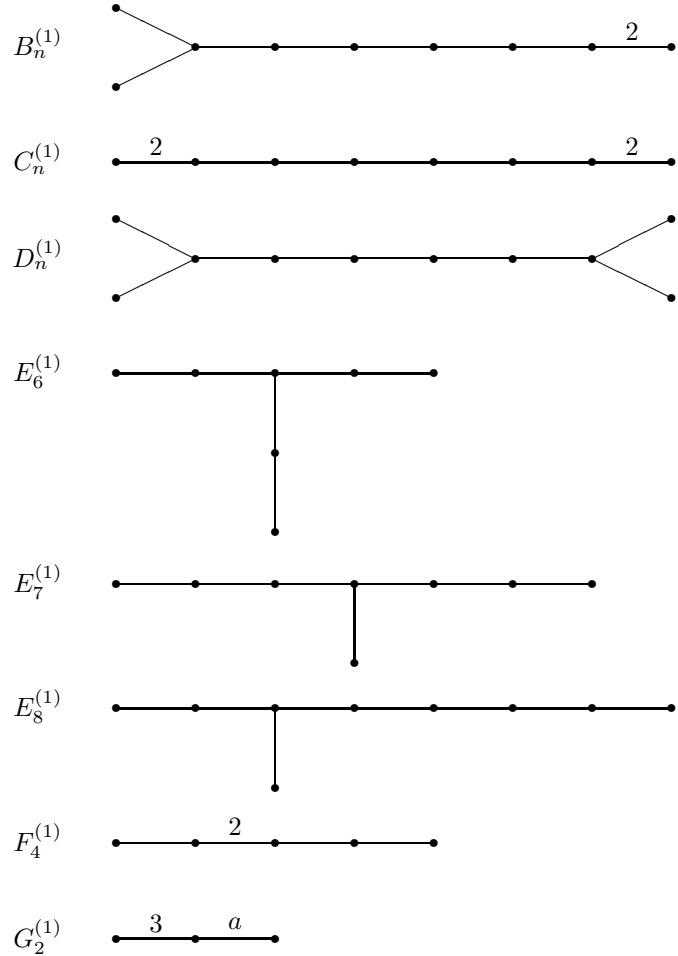

\proof 
A diagram $\Gamma$ is called an \emph{extended Dynkin tree diagram} if 
\begin{itemize} 
\item 
$\Gamma$ is a tree diagram with edge weights $\leq 3$; 
\item 
$\Gamma$ is not on the Dynkin diagram list; 
\item 
every proper subdiagram of $\Gamma$ is a disjoint union of Dynkin diagrams. 
\end{itemize} 
(In this definition, we ignore the orientations of the edges.) 
To prove the proposition, it is enough to show that any extended 
Dynkin tree diagram is $2$-infinite. 
Direct inspection shows that Figure~\ref{fig:extended-dynkin-diagrams} 
provides a complete list of such diagrams. 
Here each tree $X_n^{(1)}$ has $n+1$ vertices. 
As before, all unspecified edge weights are equal to $1$; in the 
diagram~$G_2^{(1)}$, we have $a\in\{1,2,3\}$. 
We note that all these diagrams are associated with untwisted affine Lie algebras 
and can be found in the tables in \cite{bourbaki} or in 
\cite[Chapter 4, Table Aff~1]{kac}. 
The only diagram from those tables that is missing in 
Figure~\ref{fig:extended-dynkin-diagrams} is $A_n^{(1)}$, which is 
an $(n+1)$-cycle; it will be treated in Section~\ref{sec:take-care-of-cycles}. 
 
In showing that an extended Dynkin tree diagram is 
$2$-infinite, we can choose its orientation arbitrarily, 
by Proposition~\ref{pr:tree-orientations}. 
Let us start with the three infinite series $B_n^{(1)}, C_n^{(1)}$, 
and $D_n^{(1)}$, and in each case let us orient all the edges 
left to right. 
Let us denote the diagram in question by $X_n^{(1)}$; thus, if 
$X = D$ (resp., $B,C$) then the minimal value of $n$ is equal to 
$4$ (resp., $3,2$). 
If $n$ is greater than this minimal value, then 
performing a mutation at the second vertex from the left, and 
subsequently removing this vertex (together with all incident edges) 
leaves us with a subdiagram of type~$X_{n-1}^{(1)}$. 
Using induction on~$n$, it remains to check the basic cases 
$D_4^{(1)}, B_3^{(1)}$ and $C_2^{(1)}$. 
For $C_2^{(1)}$, the mutation at the middle vertex produces a 
triangle with edge weights $(2,2,2)$, 
which is $2$-infinite by Lemma~\ref{lem:2-finite-3-ss}. 
For $B_3^{(1)}$, mutating at the branching vertex and then 
removing it leaves us with the subdiagram $C_2^{(1)}$ which was 
just shown to be $2$-infinite. 
Finally, for $D_4^{(1)}$, let the branching point be labeled by~$1$, 
and let it be joined 
with vertices $2$ and $4$ by incoming edges, and with $3$ and $5$ by 
outgoing edges. Then the composition of mutations 
$\mu_3 \circ \mu_2 \circ \mu_1$ makes the subdiagram on the vertices 
$2, 4$ and $5$ a $2$-infinite triangle. 
 
To see that $G_2^{(1)}$ is $2$-infinite, orient the two edges 
left to right and mutate at the middle vertex to obtain a 
$2$-infinite triangle. To see that $F_4^{(1)}$ is 
$2$-infinite, again orient all the edges left to right, label 
the vertices also left to right, and apply $\mu_1 
\circ \mu_2 \circ \mu_3 \circ \mu_4$ to obtain a subdiagram 
$C_2^{(1)}$ on the vertices $1,3$ and~$5$.

The remaining three cases $E_6^{(1)}, E_7^{(1)}$ and $E_8^{(1)}$ 
can be treated in a similar manner but we prefer another approach. 
To describe it, we will need to introduce some notation. 
 
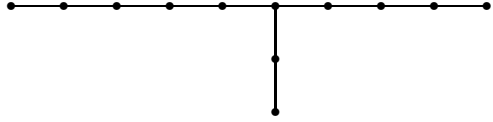
\begin{figure}[ht] 
\setlength{\unitlength}{1pt} 
\begin{picture}(180,43)(0,0) 
 
\put(0,40){\line(1,0){180}} 
\put(100,40){\line(0,-1){40}}

\multiput(0,40)(20,0){10}{\circle*{3}} 
\multiput(100,20)(0,-20){2}{\circle*{3}} 
 
\end{picture} 
\caption{The tree diagram $T_{5,4,2}\,$.} 
\label{fig:Tabc} 
\end{figure}

\begin{definition} 
\label{def:Tabc} 
For $p,q,r\in\ZZ_{\geq 0}$, we denote by 
$T_{p,q,r}$ the tree diagram (with 
all edge weights equal to $1$) on $p+q+r+1$ vertices obtained by 
connecting an endpoint of each of the three chains $A_p$, $A_q$ and~$A_r$ 
to a single extra vertex (see Figure~\ref{fig:Tabc}). 
\end{definition} 
 
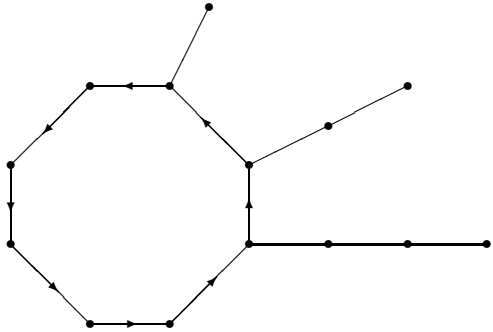
\begin{figure}[ht] 
\setlength{\unitlength}{1.5pt} 
\begin{picture}(60,83)(0,-2) 
\put(20,0){\line(1,0){20}} 
\put(20,60){\line(1,0){20}} 
\put(0,20){\line(0,1){20}} 
\put(60,20){\line(0,1){20}} 
\put(0,20){\line(1,-1){20}} 
\put(40,60){\line(1,-1){20}} 
\put(0,40){\line(1,1){20}} 
\put(40,0){\line(1,1){20}} 
 
\put(20,0){\vector(1,0){12}} 
\put(40,60){\vector(-1,0){12}} 
\put(0,40){\vector(0,-1){12}} 
\put(60,20){\vector(0,1){12}} 
\put(0,20){\vector(1,-1){12}} 
\put(60,40){\vector(-1,1){12}} 
\put(20,60){\vector(-1,-1){12}} 
\put(40,0){\vector(1,1){12}}

\multiput(20,0)(20,0){2}{\circle*{2}} 
\multiput(20,60)(20,0){2}{\circle*{2}} 
\multiput(0,20)(60,0){2}{\circle*{2}} 
\multiput(0,40)(60,0){2}{\circle*{2}} 
\put(50,80){\circle*{2}} 
 
\multiput(80,20)(20,0){3}{\circle*{2}} 
\multiput(80,50)(20,10){2}{\circle*{2}} 
\put(40,60){\line(1,2){10}} 
\put(60,20){\line(1,0){60}} 
\put(60,40){\line(2,1){40}}

 
\end{picture} 
\caption{The diagram $S_{4,3,2}^5\,$. 
} 
\label{fig:Spqrs} 
\end{figure}

\begin{definition} 
\label{def:Spqrs} 
For $p,q,r \in \ZZ_{> 0}$ and $s \in \ZZ_{\geq 0}$, 
let $S_{p,q,r}^s$ denote the diagram (with all edge weights equal 
to $1$) on $p+q+r+s$ vertices obtained by attaching three branches 
$A_{p-1}$, $A_{q-1}$, and $A_{r-1}$ to three consecutive vertices 
on a \emph{cyclically oriented} $(s+3)$-cycle (see Figure~\ref{fig:Spqrs}). 
\end{definition} 
 
\vspace{-.1in} 
 
\begin{lemma} 
\label{lem:crown} 
The diagram $S_{p,q,r}^s$ is mutation equivalent to $T_{p+r-1,q,s}$. 
\end{lemma} 
 
\proof 
Let us consider the subdiagram of $S_{p,q,r}^s$ obtained by removing the 
middle branch~$A_q\,$. 
This subdiagram is a copy of $A_{p+s+r}\,$. 
We label its vertices consecutively by $1,\dots,p+s+r$, starting with 
the endpoint of~$A_r$; and we orient the edges of 
$A_p$ and $A_r$ so that all the edges of $A_{p+s+r}$ point at the 
same direction. 
Now a direct check shows that 
$\mu_1 \circ \mu_2 \circ \cdots \circ \mu_{s+r}$ transforms 
$S_{p,q,r}^s$ into $T_{p+r-1,q,s}$. 
\endproof 
 
The proof of Proposition~\ref{pr:trees-Dynkin-vs-infty} can now be 
completed as follows: 
\[ 
\begin{array}{l} 
E_6^{(1)} = T_{2,2,2} \sim S_{2,2,1}^2 \supset D_5^{(1)} \,;\\ 
E_7^{(1)} = T_{3,1,3} \sim S_{3,1,1}^3 \supset E_6^{(1)} \,;\\ 
E_8^{(1)} = T_{2,1,5} \sim S_{2,1,1}^5 \supset E_7^{(1)} \,. \qed 
\end{array} 
\] 
 
\subsection{Taking care of the cycles} 
\label{sec:take-care-of-cycles} 
 
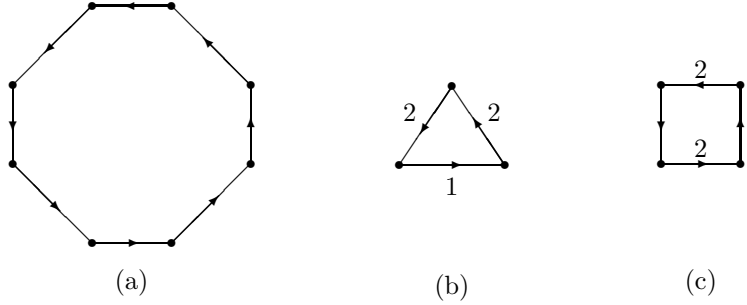
\begin{figure} 
\setlength{\unitlength}{1.5pt} 
\begin{picture}(60,71)(0,-7) 
\put(20,0){\line(1,0){20}} 
\put(20,60){\line(1,0){20}} 
\put(0,20){\line(0,1){20}} 
\put(60,20){\line(0,1){20}} 
\put(0,20){\line(1,-1){20}} 
\put(40,60){\line(1,-1){20}} 
\put(0,40){\line(1,1){20}} 
\put(40,0){\line(1,1){20}} 
 
\put(20,0){\vector(1,0){12}} 
\put(40,60){\vector(-1,0){12}} 
\put(0,40){\vector(0,-1){12}} 
\put(60,20){\vector(0,1){12}} 
\put(0,20){\vector(1,-1){12}} 
\put(60,40){\vector(-1,1){12}} 
\put(20,60){\vector(-1,-1){12}} 
\put(40,0){\vector(1,1){12}}

\multiput(20,0)(20,0){2}{\circle*{2}} 
\multiput(20,60)(20,0){2}{\circle*{2}} 
\multiput(0,20)(60,0){2}{\circle*{2}} 
\multiput(0,40)(60,0){2}{\circle*{2}} 
 
\put(30,-10){\makebox(0,0){(a)}} 
 
\end{picture} 
\hspace{.5in} 
\setlength{\unitlength}{2pt} 
\begin{picture}(30,17)(-5,-17) 
\put(0,3){\line(1,0){20}} 
\put(0,3){\vector(1,0){12}} 
\put(0,3){\line(2,3){10}} 
\put(10,18){\vector(-2,-3){6}} 
\put(10,18){\line(2,-3){10}} 
\put(20,3){\vector(-2,3){6}} 
\put(0,3){\circle*{1.5}} 
\put(20,3){\circle*{1.5}} 
\put(10,18){\circle*{1.5}} 
\put(2,13){\makebox(0,0){$2$}} 
\put(18,13){\makebox(0,0){$2$}} 
\put(10,-1){\makebox(0,0){$1$}} 
\put(10,-20){\makebox(0,0){(b)}} 
\end{picture} 
\hspace{.5in} 
\setlength{\unitlength}{1.5pt} 
\begin{picture}(25,27)(-2,-27) 
\multiput(0,0)(0,20){2}{\line(1,0){20}} 
\put(20,20){\vector(-1,0){12}} 
\put(0,0){\vector(1,0){12}} 
\multiput(0,0)(20,0){2}{\line(0,1){20}} 
\put(0,20){\vector(0,-1){12}} 
\put(20,0){\vector(0,1){12}} 
\multiput(0,0)(20,0){2}{\circle*{2}} 
\multiput(0,20)(20,0){2}{\circle*{2}} 
\multiput(10,4)(0,20){2}{\makebox(0,0){$2$}} 
\put(10,-30){\makebox(0,0){(c)}} 
\end{picture} 
\vspace{.1in} 
\caption{$2$-finite cycles 
} 
\label{fig:3-4-cycles} 
\end{figure}

\begin{proposition} 
\label{pr:cycles-Dynkin-vs-infty} 
Let $\Gamma$ be a $2$-finite diagram whose underlying graph is an 
$n$-cycle for some $n \geq 3$ 
(with some orientation of edges). 
Then 
$\Gamma$ must be one of the diagrams shown in 
Figure~\ref{fig:3-4-cycles}. 
More precisely, one of the following holds: 
\begin{itemize} 
 
\item[{\rm(a)}] $\Gamma$ is an oriented cycle with all weights equal 
 to~$1$. \\ 
In this case, $\Gamma\sim D_n$ (with the understanding that $D_3 =A_3$). 
 
\item[{\rm(b)}] $\Gamma$ is an oriented triangle with 
edge weights $2,2,1$ shown in 
 Figure~\ref{fig:3-4-cycles}(b). \\ 
In this case, $\Gamma\sim B_3$. 
 
\item[{\rm(c)}] $\Gamma$ is an oriented $4$-cycle with 
edge weights $2,1,2,1$ shown in 
 Figure~\ref{fig:3-4-cycles}(c). \\ 
In this case, $\Gamma\sim F_4$. 
 
\end{itemize} 
In particular, the edges in $\Gamma$ must be cyclically oriented. 
\end{proposition} 
 
\proof 
The case $n=3$ of Proposition~\ref{pr:cycles-Dynkin-vs-infty} follows 
from Lemmas~\ref{lem:2-finite-3-cycle} and~\ref{lem:2-finite-3-ss}, 
so for the rest of the proof we assume that $n\geq 4$. 
 
We begin by proving the last claim of 
Proposition~\ref{pr:cycles-Dynkin-vs-infty} by induction on~$n$. 
Invoking if necessary a shape-preserving mutation, we may assume that there 
is a vertex $v \in \Gamma$ that has one incoming and one outgoing edge. 
Let $\Gamma' = \mu_v (\Gamma)$. 
Then the subdiagram $\Gamma''=\Gamma' - \{v\}$ is an $(n-1)$-cycle, which must be 
cyclically oriented by the induction assumption. Backtracking to~$\Gamma$, we 
obtain the desired claim. 
 
Furthermore, observe that the product of edge weights of $\Gamma''$ is 
the same as in~$\Gamma$. Again using induction together with 
Lemma~\ref{lem:2-finite-3-ss}, we conclude that this product is either 
1 or~4. 
In the former case, $\Gamma$ is an oriented 
$n$-cycle, and we apply Lemma~\ref{lem:crown} to obtain 
$\Gamma = S_{1,1,1}^{n-3}\sim T_{1,1,n-3} = D_n\,$, as needed. 
In the latter case, $\Gamma$ has two edge weights equal to~2, and the 
rest of them are equal to~1. 
Then either $\Gamma$ is one of the two diagrams (b) and (c) in 
Figure~\ref{fig:3-4-cycles}, or else it contains a 
$2$-infinite subdiagram $C_m^{(1)}$ for some $m \geq 2$. 
It remains to show that the diagrams in 
Figure~\ref{fig:3-4-cycles}(b)-(c) are mutation equivalent to $B_3$ 
and $F_4\,$, respectively. This is straightforward. 
\endproof

\subsection{Completing the proof of Theorem~\ref{th:diagrams-fin-CK}} 
The second claim in Theorem~\ref{th:diagrams-fin-CK} follows from 
Proposition~\ref{pr:tree-orientations}, so we only need to show 
that a connected $2$-finite diagram~$\Gamma$ is 
mutation equivalent to some Dynkin diagram. 
We proceed by induction on $n$, the number of vertices in 
$\Gamma$. 
If $n \leq 3$, then $\Gamma$ is either a tree or a cycle, and the 
theorem follows by 
Propositions~\ref{pr:trees-Dynkin-vs-infty} and 
\ref{pr:cycles-Dynkin-vs-infty}. 
So let us assume that the statement is already known for some 
$n \geq 3$; we need to show that it holds for a diagram $\Gamma$ on 
$n+1$ vertices. 
Pick a vertex $v\in\Gamma$ such that the subdiagram 
$\Gamma' = \Gamma - \{v\}$ is connected. 
Since $\Gamma'$ is $2$-finite, it is 
mutation equivalent to some Dynkin diagram~$X_n\,$. 
Furthermore, we may assume that $\Gamma'$ is (isomorphic to) 
our favorite representative of the mutation equivalence class 
of~$X_n\,$. For each $X_n\,$, we will choose a representative that is 
most convenient for the purposes of this proof. 
 
\credit{Case 1} 
\emph{$\Gamma'$ is a Dynkin diagram 
with no branching point, i.e., is 
of one of the types $A_n$, $B_n$, $F_4$, or $G_2$.} 
Let us orient the edges of $\Gamma'$ so that they all point in the 
same direction. 
If $v$ is adjacent to exactly one vertex of $\Gamma'$, then 
$\Gamma$ is a tree, and we are done by 
Proposition~\ref{pr:trees-Dynkin-vs-infty}. 
If $v$ is adjacent to more than $2$ vertices of $\Gamma'$, then 
$\Gamma$ has a cycle subdiagram whose edges are not cyclically 
oriented, contradicting Proposition~\ref{pr:cycles-Dynkin-vs-infty}. 
Therefore we may assume that $v$ is adjacent to precisely two vertices 
$v_1$ and $v_2$ of $\Gamma'$. 
Thus, $\Gamma$ has precisely one cycle $\mathcal{C}$, which 
furthermore must be of one of the types (a)--(c) described in 
Proposition~\ref{pr:cycles-Dynkin-vs-infty} and Figure~\ref{fig:3-4-cycles}. 
 
\credit{Subcase 1.1} 
\emph{$\mathcal{C}$ is an oriented cycle with unit edge weights.} 
If $\Gamma$ has an edge of weight $\geq 2$, then it contains a 
subdiagram of type $B_m^{(1)}$ or $G_2^{(1)}$, unless $\mathcal{C}$ is 
a $3$-cycle, in which case $\mu_v(\Gamma)\sim B_{n+1}\,$. 
On the other hand, if all edges in $\Gamma$ are of weight~$1$, then it is one of 
the diagrams $S_{p,q,r}^s$ in Lemma~\ref{lem:crown} (with $q=0$). 
Hence $\Gamma$ is mutation equivalent to a tree, and we are done by 
Proposition~\ref{pr:trees-Dynkin-vs-infty}. 
 
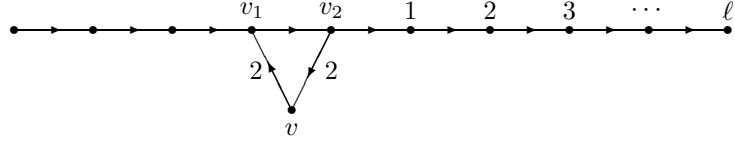
\begin{figure} 
\setlength{\unitlength}{1.5pt} 
\begin{picture}(180,30)(0,-20) 
\put(0,0){\line(1,0){180}} 
\multiput(0,0)(20,0){10}{\circle*{2}} 
\put(70,-20){\circle*{2}} 
\put(70,-20){\line(-1,2){10}} 
\put(70,-20){\line(1,2){10}} 
 
\put(80,0){\vector(-1,-2){6}} 
\put(70,-20){\vector(-1,2){6}} 
 
\multiput(0,0)(20,0){9}{\vector(1,0){12}} 
 
\put(70,-25){\makebox(0,0){$v$}} 
\put(60,5){\makebox(0,0){$v_1$}} 
\put(61,-10){\makebox(0,0){$2$}} 
\put(80,-10){\makebox(0,0){$2$}} 
\put(80,5){\makebox(0,0){$v_2$}} 
\put(100,5){\makebox(0,0){$1$}} 
\put(120,5){\makebox(0,0){$2$}} 
\put(140,5){\makebox(0,0){$3$}} 
\put(160,5){\makebox(0,0){$\cdots$}} 
\put(180,5){\makebox(0,0){$\ell$}} 
\end{picture} 
\caption{Subcase 1(b) 
} 
\label{fig:B-triangle-with-ears} 
\end{figure}

\credit{Subcase 1.2} 
\emph{$\mathcal{C}$ is as in 
Figure~\ref{fig:3-4-cycles}(b).} 
If one of the edges $(v,v_1)$ and $(v,v_2)$ has weight~$1$, 
then $\mu_v$ removes the edge $(v_1, v_2)$, resulting in a tree, and 
we are done again. 
So assume that both $(v,v_1)$ and $(v,v_2)$ have weight~$2$. 
If at least one edge outside $\mathcal{C}$ 
has weight $\geq 2$, then $\Gamma\supset C_m^{(1)}$ or 
$\Gamma\supset G_2^{(1)}\,$. 
It remains to consider the case shown in Figure~\ref{fig:B-triangle-with-ears} 
(as before, unspecified edge weights are equal to~$1$). 
Direct check shows that 
$\mu_l \circ \cdots \circ \mu_2 \circ \mu_1 \circ \mu_{v_2} 
\circ \mu_v(\Gamma)=B_{n+1}\,$, and we are done.

\credit{Subcase 1.3} \emph{$\mathcal{C}$ is as in 
Figure~\ref{fig:3-4-cycles}(c).} 
It suffices to show that any diagram ob\-tained from $\mathcal{C}$ 
by adjoining a single vertex adjacent to one of its vertices 
is $2$-infinite. 
If this extra edge has weight $1$ (resp., $2$, $3$), then 
the resulting $5$-vertex diagram has a $2$-infinite subdiagram of type 
$B_3^{(1)}$ (resp., $C_2^{(1)}$, $G_2^{(1)}$), proving the claim. 
 
\credit{Case 2} 
\emph{$\Gamma'\sim D_n$ ($n\geq 4$).} 
By Proposition~\ref{pr:cycles-Dynkin-vs-infty}(a), we may assume that 
$\Gamma'$ is an oriented $n$-cycle with unit edge weights. 
 
\credit{Subcase 2.1} \emph{$v$ is adjacent to a single 
vertex $v_1\in\Gamma'$.} 
If the edge $(v,v_1)$ has weight $\geq 2$, then 
$\Gamma$ has a subdiagram $B_3^{(1)}$ or 
$G_2^{(1)}$. 
If the edge $(v,v_1)$ has weight $1$, then 
by Lemma~\ref{lem:crown}, $\Gamma$ is mutation equivalent to a 
tree, and we are done again. 
 
\credit{Subcase 2.2} \emph{$v$ is adjacent to exactly two 
vertices $v_1$ and $v_2$ of $\Gamma'$, which are adjacent to each other.} 
Then the triangle $(v,v_1,v_2)$ is either an oriented $3$-cycle with 
unit edge weights or the diagram in Figure~\ref{fig:3-4-cycles}(b). 
In the former case, $\mu_v(\Gamma)$ is an oriented $(n+1)$-cycle, so 
$\Gamma\sim D_{n+1}$. 
In the latter case, $\mu_v$ reverses the orientation of the edge 
$(v_1,v_2)$, transforming $\Gamma'$ into an improperly oriented (hence 
$2$-infinite) cycle (cf.\ Proposition~\ref{pr:cycles-Dynkin-vs-infty}). 
 
\credit{Subcase 2.3} \emph{$v$ is adjacent to two 
non-adjacent vertices of $\Gamma'$ (and maybe to some other vertices).} 
In this case, $\Gamma$ contains a subdiagram which is a 
non-cyclically-oriented cycle, contradicting 
Proposition~\ref{pr:cycles-Dynkin-vs-infty}.

\credit{Case 3} 
\emph{$\Gamma'\sim E_n=T_{1,2,n-4}\,$, for $n\!\in\!\{6,7,8\}$.} 
By Lemma~\ref{lem:crown}, we may assume~that $\Gamma'=S_{1,2,1}^{n-4}\,$. 
In other words, $\Gamma'$ is a cyclically oriented $(n-1)$-cycle 
$\mathcal{C}$ with unit edge weights, and an extra edge of weight $1$ 
connecting a vertex in 
$\mathcal{C}$ to a vertex 
$v_1\notin\mathcal{C}$. 
 
\credit{Subcase 3.1} \emph{$v$ is adjacent to $v_1$, 
and to no other vertices in $\Gamma'$.} 
If the edge $(v,v_1)$ has weight $\geq 2$, then 
$\Gamma$ has a $2$-infinite subdiagram $B_3^{(1)}$ or 
$G_2^{(1)}$. 
If $(v,v_1)$ has weight~$1$, then 
by Lemma~\ref{lem:crown}, $\Gamma$ is mutation equivalent to a 
tree. 
 
\credit{Subcase 3.2} \emph{$v$ is adjacent to a 
vertex $v_2\in\mathcal{C}$, 
and to no other vertices in $\Gamma'$. 
} 
Then $\Gamma$ has a subdiagram 
of type $D_m^{(1)}$ or $B_3^{(1)}$ or $G_2^{(1)}$. 
 
\credit{Subcase 3.3} \emph{$v$ is adjacent to at least two 
vertices in~$\mathcal{C}$. 
} 
By the analysis in Subcases 2.2--2.3 (with $\Gamma'$ replaced 
by~$\mathcal{C}$), we must have $\Gamma-\{v_1\}\sim D_n\,$, 
and the problem reduces to Case~2 already treated above 
(the role of $v$ now played by~$v_1$). 
 
\credit{Subcase 3.4:} \emph{$v$ is adjacent to $v_1$ and a single 
vertex $v_2\in\mathcal{C}$.} 
Let $v_0$ be the only vertex on $\mathcal{C}$ adjacent to $v_1$ 
in~$\Gamma'$. 
If $v_2$ is neither $v_0$ nor a vertex adjacent to $v_0$, then the 
three cycles in $\Gamma$ cannot be simultaneously oriented. 
If $v_2=v_0$, then $\mu_{v_1}$ removes the edge $(v,v_2)$, 
transforming $\Gamma$ into a diagram mutation equivalent to a tree by 
Lemma~\ref{lem:crown}. 
If $v_2$ is adjacent to $v_0$, then $\Gamma-\{v_0\}$ has no branching 
point, and the problem reduces to Case~1 
already treated above, with the role of $v$ now played by~$v_0\,$.

This concludes the proof of 
Theorem~\ref{th:diagrams-fin-CK}. 
As a consequence, we obtain Theorems~\ref{th:2-fin-CK} 
and~\ref{th:finite-type-class-exhaust}. 
\endproof 
 
\pagebreak[2] 
 
\section{Proof of 
  Theorem~\ref{th:series-CK} 
} 
\label{sec:proof-CK-exactly} 
 
 
Let $B$ and $B'$ be sign-skew-symmetric matrices such that both $A\!=\! A(B)$ 
and~$A'\!=\! A(B')$ are Cartan matrices of finite type. 
We already proved that both $B$ and $B'$ are $2$-finite. 
We need to show that $B$ and  $B'$ are mutation equivalent if and 
only~if $A$ and $A'$ are of the same 
type. 
Without loss of generality, we may assume that~$A$ and $A'$ are 
indecomposable, i.e., the corresponding root systems $\Phi$ and 
$\Phi'$ are irreducible. 
 
We first prove the ``only if'' part. 
If $B$ and $B'$ are mutation 
equivalent, then the simplicial complexes $\Delta (\Phi)$ and 
$\Delta (\Phi')$ are isomorphic to each other, 
by Theorem~\ref{th:finite-type-complex}. 
In particular, $\Phi$ and $\Phi'$ have the same rank and the same cardinality. 
A direct check using the tables in \cite{bourbaki} shows 
that the only different Cartan-Killing types with this property 
are $B_n$ and $C_n$ for all $n \geq 3$, and also $E_6$, which has 
the same data as $B_6$ and~$C_6$. 
To distinguish between these types, note that 
mutation-equivalent skew-symmetrizable matrices 
share the same skew-symmetrizing matrix~$D$. 
Furthermore, $D$ is skew-symmetrizing for $B$ if and only if 
it is symmetrizing for~$A$; thus, the diagonal entries of $D$ are given by 
$d_i = (\alpha_i,\alpha_i)$, where $(\alpha, \beta)$ is a 
$W$-invariant scalar product on the root lattice. 
Since the root system of type $B_n$ (resp.,~$C_n$) 
has one short simple root and $n-1$ long ones 
(resp.,~one long and $n-1$ short), the 
corresponding matrices $B$ and $B'$ cannot be mutation equivalent. 
The same is true for $E_6$ and $B_6$ (or $C_6$) since all simple 
roots for $E_6$ are of the same length.

To prove the ``if'' part, 
suppose that 
$A$ and $A'$ are of the same Cartan-Killing type. 
By Proposition~\ref{pr:tree-orientations}, we may assume without 
loss of generality that $B$ and~$B'$ have the same diagram. 
By Lemma~\ref{lem:ss-conjugation}, we have $S(B) = S(B')$. 
Since $B$ and $B'$ share a skew-symmetrizing matrix~$D$, 
the proof of Lemma~\ref{lem:ss-conjugation} shows that 
$S(B) = HBH^{-1}$ and $S(B') = HB'H^{-1}$ for $H = D^{1/2}$. 
Hence $B = B'$, and we are done. 
\endproof 
 
\pagebreak[3] 
 
\section{On cluster algebras of geometric type} 
\label{sec:geometric-type-ca} 
In this section we present two general results 
on 
cluster algebras of geometric type in the sense of 
\cite[Definition~5.7]{fz-clust1}. 
These algebras are \emph{not} assumed to be of finite type, so all 
the 
necessary background 
is contained in 
Section~\ref{sec:def-of-normalized-ca}. 
 
Recall that a cluster algebra is of \emph{geometric type} if it 
satisfies the following two conditions: 
\begin{align} 
\label{eq:geometric-1} 
&\text{The coefficient semifield $\PP$ is of the form $\Trop (p_j: j 
  \in  J)$. 
That is,} \\ 
\nonumber 
&\text{the multiplicative group of $\PP$ is a free abelian group 
with a finite set} \\ 
\nonumber 
&\text{of generators $p_j \, (j \in J)$, and the 
auxiliary addition $\oplus$ is given by 
(\ref{eq:tropical addition}).} 
\\ 
&\label{eq:geometric-2} 
\text{Every element $p \in \mathcal{P}$, i.e., every coefficient in one 
of the exchange}\\ 
\nonumber 
&\text{relations (\ref{eq:exchange-rel-xx}), is a 
monomial in the $p_j$ with all exponents nonnegative.} 
\end{align} 
 
We note a little discrepancy between our choice of the ground ring 
$\ZZ[\mathcal{P}]$ in Definition~\ref{def:cluster-algebra} 
and the choice described in \cite[Section~5]{fz-clust1}, where the 
ground ring was taken to be the polynomial ring $\ZZ[p_j: j \in J]$. 
The following additional assumption 
guarantees that these two choices coincide: 
\begin{equation} 
\label{eq:geometric-3} 
\text{Every generator $p_j$ of $\PP$ belongs to $\mathcal{P}$.} \hspace{2.05in} 
\end{equation} 
 
\subsection{Geometric realization criterion} 
Our first result gives sufficient 
conditions under which a cluster algebra of geometric type 
can be realized as  a $\ZZ$-form of the coordinate ring $\CC[X]$ 
of some algebraic variety $X$. 
 
We make the following assumptions on $X$: 
\begin{equation} 
\label{eq:X-rational-etc} 
\text{$X$ is a rational quasi-affine irreducible algebraic variety 
over $\CC$.} \hspace{.5in} 
\end{equation} 
Irreducibility implies that the ring of regular functions $\CC[X]$ is a 
domain, so its fraction field is well defined. 
Quasi-affine means Zariski open in some affine variety; 
this condition is imposed to ensure that the fraction field of 
$\CC[X]$ coincides with the usual field $\CC(X)$ of rational functions 
on~$X$. 
Rationality means that $X$ is birationally isomorphic 
to an affine space, i.e., 
$\CC(X)$ 
is isomorphic to the field of rational functions 
over $\CC$ in $\dim (X)$ independent variables. 
 
Let $\AA$ be a cluster algebra  of rank $n$ whose coefficient 
system satisfies conditions (\ref{eq:geometric-1})--(\ref{eq:geometric-3}), and 
let $\mathcal{X}$ be the set of cluster variables in~$\AA$. 
Suppose the variety $X$ satisfies 
\begin{equation} 
\label{eq:X-dimension} 
\dim (X) = n + |J|; 
\end{equation} 
also suppose we are given a family of functions 
$$\{\varphi_y : y \in \mathcal{X}\} \cup 
\{\varphi_j : j \in J\}$$ 
in $\CC[X]$ satisfying the following conditions: 
\begin{eqnarray} 
\label{eq:phi-generate-X} 
&&\text{the functions $\varphi_y$ and $\varphi_j$ 
generate~$\CC[X]$;}\\ 
\label{eq:same-relations} 
&&\text{every exchange relation (\ref{eq:exchange-rel-xx}) becomes 
an identity in $\CC[X]$ if we replace}\\ 
\nonumber 
&&\text{each cluster 
variable $y$ by $\varphi_y$, and each coefficient $p_z^\pm \!=\! \prod_{j 
  \in J} p_j^{a_j}$ by $\prod_{j \in J} \varphi_j^{a_j}$.} \hspace{-.2in} 
\end{eqnarray}

\begin{proposition} 
\label{pr:geom-realization-criterion} 
Under conditions {\rm (\ref{eq:geometric-1})--(\ref{eq:same-relations})}, 
the correspondence 
\begin{equation} 
\label{eq:varphi} 
y \mapsto \varphi_y \,\, (y \in \mathcal{X}), 
\quad p_j \mapsto \varphi_j \,\, (j \in J) 
\end{equation} 
extends uniquely to an algebra isomorphism between the cluster algebra 
$\AA$ and the $\ZZ$-form of $\CC[X]$ generated by all $\varphi_y$ 
and~$\varphi_j\,$. 
\end{proposition} 
 
\proof 
Pick an arbitrary cluster $\xx$ of $\AA$, and let 
$\tilde \xx = \xx \cup \{p_j: j \in J\}$. 
Since $\xx$ is a transcendence basis of the ambient field $\FFcal$ 
over $\ZZ \PP$, the set $\tilde \xx$ is a transcendence basis of $\FFcal$ 
over $\QQ$. 
Furthermore, every  cluster variable is uniquely expressed as a rational function in 
$\tilde \xx$ by iterating the exchange relations away from a seed 
containing $\xx$ in the exchange graph of $\AA$. 
In view of (\ref{eq:same-relations}), we can apply the same 
procedure to express all functions $\varphi_y$ and $\varphi_j$ 
inside the field $\CC(X)$ as rational functions in the set 
$$\varphi (\tilde \xx) = \{\varphi_x: x \in \xx\} \cup \{\varphi_j : j \in J\}\, .$$ 
Furthermore, we have $|\varphi (\tilde \xx)| = \dim (X)$ by 
(\ref{eq:X-dimension}). 
Since $X$ is rational, we conclude from (\ref{eq:phi-generate-X}) 
that $\varphi (\tilde \xx)$ is a transcendence basis of the field 
of rational functions $\CC(X)$, and that the correspondence (\ref{eq:varphi}) extends 
to an embedding of fields $\FFcal \to \CC(X)$, and hence to an 
embedding of algebras $\AA \to \CC[X]$. 
This proves Proposition~\ref{pr:geom-realization-criterion}. 
\endproof 
 
\subsection{Sharpening the Laurent phenomenon} 
As mentioned in Section~\ref{sec:laurent-finite-type}, 
the \emph{Laurent phenomenon}, established in~\cite{fz-clust1} for 
arbitrary cluster algebras, says that 
every cluster variable can be written as a Laurent polynomial in the 
variables of an arbitrary fixed cluster, with coefficients in~$\ZZ \PP$. 
For the cluster algebras of geometric type, this result can be 
sharpened as follows. 
 
\begin{proposition} 
\label{rem:Laurent-sharper} 
In any cluster algebra with the coefficient system satisfying 
conditions {\rm (\ref{eq:geometric-1})--(\ref{eq:geometric-3})}, 
every cluster variable is expressed in terms of an arbitrary 
cluster $\xx$ as a Laurent polynomial with coefficients in 
$\ZZ[\mathcal{P}]$. 
\end{proposition} 
 
\proof 
Fix some generator $p \!=\! p_{j_\circ}$ of the coefficient semifield 
$\PP \!=\! \Trop (p_j\!:\! j \!\in\!  J)$. 
We will think of any cluster variable $z$ as a Laurent polynomial 
$z(p)$ whose coefficients are integral Laurent polynomials in the 
set $\xx \cup \{p_j: j \in J, j \neq j_\circ\}$. 
Our goal is to show that $z(p)$ is in fact a polynomial in~$p$; 
Proposition~\ref{rem:Laurent-sharper} will then follow by varying a 
distinguished index~$j_\circ$ over the index set~$J$. 
 
Define the \emph{distance} $d(z,\xx)$ 
between $z$ and $\xx$ as the shortest distance 
in the exchange graph between a seed containing $z$ and a seed whose cluster 
is~$\xx$. 
We will use induction on $d(z,\xx)$ to show the following 
strengthening of the desired statement: 
\begin{itemize} 
 
\item 
$z(p)$ is a polynomial in $p$ 
whose constant term $z(0)$ is a subtraction-free rational 
expression in $\xx \cup \{p_j: j \in J, j \neq j_\circ\}$ 
(in particular, $z(0)\neq 0$). 
\end{itemize} 
 
If $d(z,\xx) = 0$, then $z \in \xx$, and there is nothing to prove. 
If $d(z,\xx) > 0$, then, by the definition of the distance, $z$ 
participates in an exchange relation (\ref{eq:exchange-rel-xx}) 
such that all the other participating cluster variables are at a 
smaller distance from $\xx$ than~$z$. 
Applying the inductive assumption to all these cluster variables 
and using Lemma~\ref{lem:laurent-implies-poly} together with 
the fact that, by the normalization condition, $p$ appears in at most 
one of the monomials on the right hand side of 
(\ref{eq:exchange-rel-xx}), we obtain our claim for~$z$. 
\endproof

\section{Examples of geometric realizations of cluster algebras} 
\label{sec:geometric-realization} 
In this section, we present some examples of concrete 
geometric realizations of cluster algebras $\AA = \AA(B,\pp)$ of finite type. 
In all these examples, the Cartan counterpart of $B$ 
is a Cartan matrix of one of the classical types 
$A_n, B_n, C_n, D_n$, and the coefficient system of $\AA$ 
satisfies conditions (\ref{eq:geometric-1})--(\ref{eq:geometric-3}).

\subsection{Type $A_1$} 
\label{sec:ca-type-A1} 
We start by presenting four natural geometric realizations of 
cluster algebras of type $A_1$. 
Such an algebra $\AA$ has only two one-element clusters~$\{x\}$ 
and~$\{\overline x\}$, and a single exchange relation 
\begin{equation} 
\label{eq:exchange-A1} 
x \overline x = p^+ + p^- \, , 
\end{equation} 
where $p^+$ and $p^-$ belong to the coefficient semifield $\PP$. 
By Definition~\ref{def:cluster-algebra}, $\AA$ is a subalgebra 
of the ambient field $\FFcal$ generated by $x, \overline x, p^+$, 
and $p^-$. 
 
\begin{example} 
\label{ex:SL2-A1} 
Let $\AA$ have the coefficient semifield $\PP = \Trop (p)$ 
(the free abelian group with one generator), 
and let the coefficients in (\ref{eq:exchange-A1}) be given 
by $p^+ = p$ and $p^- = 1$. 
Let $G = SL_2 (\CC)$ be the group of complex matrices 
$$\bmat{a}{b}{c}{d}$$ 
with $ad - bc =1$. 
The correspondence 
$$x \mapsto a, \quad \overline x \mapsto d, \quad p \mapsto bc$$ 
identifies $\AA$ with the subring of the coordinate 
ring $\CC[G]$ generated by $a, d$, and $bc$. 
It is easy to see that this ring is a $\ZZ$-form of the ring of invariants 
$\CC[G]^H$, where $H$ is the maximal torus of diagonal matrices in 
$G$ acting on $G$ by conjugation. 
\end{example} 
 
The next three examples give three different realizations of the 
same cluster algebra $\AA$ for which the coefficients in 
(\ref{eq:exchange-A1}) are the generators of $\PP = \Trop(p^+,p^-)$. 
 
\begin{example} 
\label{ex:N3-A1} 
Let $N$ be the group of complex matrices of the form 
$$\left[\!\!\begin{array}{ccc} 
1 & a & c \\ 
0 & 1 & b \\ 
0 & 0 & 1\\ 
\end{array}\!\!\right] \, .$$ 
The correspondence 
$$x \mapsto a, \quad \overline x \mapsto b, 
\quad p^+ \mapsto c, \quad p^- \mapsto ab-c$$ 
identifies $\AA$ with a $\ZZ$-form 
$\ZZ[N] = \ZZ[a,b,c]$ of the ring $\CC[N]$. 
\end{example} 
 
\begin{example} 
\label{ex:GmodN3-A1} 
Let $G = SL_3 (\CC)$, and let $N \subset G$ be the same as in 
Example~\ref{ex:N3-A1}. 
Let $X = G/N$ be the \emph{base affine space} of $G$ taken in the 
standard embedding into 
$\CC^3 \times \bigwedge^2 \CC^3$. 
Let $(\Delta_1, \Delta_2, \Delta_3)$ and $(\Delta_{12}, \Delta_{13}, \Delta_{23})$ 
(here $\Delta_{ij}=\Delta_i\wedge \Delta_j$) 
be the standard (Pl\"ucker) coordinates in~$\CC^3$ and~$\bigwedge^2 
\CC^3$, respectively. 
In these coordinates, 
the coordinate ring of $X$ is given by 
\[ 
\CC[X] = \CC[\Delta_1, \Delta_2, \Delta_3, \Delta_{12}, \Delta_{13}, \Delta_{23}]/ 
\langle \Delta_1 \Delta_{23} - \Delta_2 \Delta_{13} + \Delta_3 \Delta_{12}\rangle . 
\] 
The correspondence 
$$x \mapsto \Delta_2, \quad \overline x \mapsto \Delta_{13}, 
\quad p^+ \mapsto \Delta_1 \Delta_{23}, \quad p^- \mapsto \Delta_3 \Delta_{12}$$ 
identifies $\AA$ with the subring of $\CC[X]$ generated by 
$\Delta_2, \Delta_{13}, \Delta_1 \Delta_{23}$, and $\Delta_3 \Delta_{12}$. 
It is easy to see that this ring is a $\ZZ$-form of the ring of invariants 
$\CC[X]^T$, where 
$T \subset G$ is the torus of all 
diagonal matrices of the form 
$$\left[\!\!\begin{array}{ccc} 
t & 0 & 0 \\ 
0 & 1 & 0 \\ 
0 & 0 & t^{-1}\\ 
\end{array}\!\!\right] \, ,$$ 
acting on $X$ by left translations. 
\end{example} 
 
\begin{example} 
\label{ex:Gr24-A1} 
Let $X 
\subset \bigwedge^2 \CC^4$ 
be the affine cone over the Grassmannian ${\rm Gr}_{2,4}$ 
taken in its Pl\"ucker embedding. 
In the standard coordinates $(\Delta_{ij}: 1 \leq i < j \leq 4)$ 
on $\bigwedge^2 \CC^4$, 
the coordinate ring of $X$ is given by 
$$\CC[X] = \CC[(\Delta_{ij})]/ 
\langle\Delta_{12} \Delta_{34} - \Delta_{13} \Delta_{24} + \Delta_{14} 
\Delta_{23}\rangle\,. 
$$ 
The correspondence 
$$x \mapsto \Delta_{13}, \quad \overline x \mapsto \Delta_{24}, 
\quad p^+ \mapsto \Delta_{12} \Delta_{34}, \quad p^- \mapsto \Delta_{14} \Delta_{23}$$ 
identifies $\AA$ with the subring of $\CC[X]$ generated by 
$\Delta_{13}, \Delta_{24}, \Delta_{12} \Delta_{34}$, and $\Delta_{14} \Delta_{23}$. 
This ring is a $\ZZ$-form of the ring of invariants 
$\CC[X]^T$, where 
$T \subset SL_4$ is the torus of all 
diagonal matrices of the form 
$$\left[\!\!\begin{array}{cccc} 
t_1 & 0 & 0 & 0\\ 
0 & t_2 & 0 & 0\\ 
0 & 0 & t_1^{-1} & 0\\ 
0 & 0 & 0 & t_2^{-1}\\ 
\end{array}\!\!\right] \, ,$$ 
naturally acting on $X$. 
\end{example} 
 
\subsection{Type $A_n \,\, (n \geq 2)$} 
\label{sec:ca-type-A}

Here we present a geometric realization of a cluster algebra of type 
$A_n$ for all $n \geq 2$, for a special choice of a coefficient system, 
to be specified below. 
 
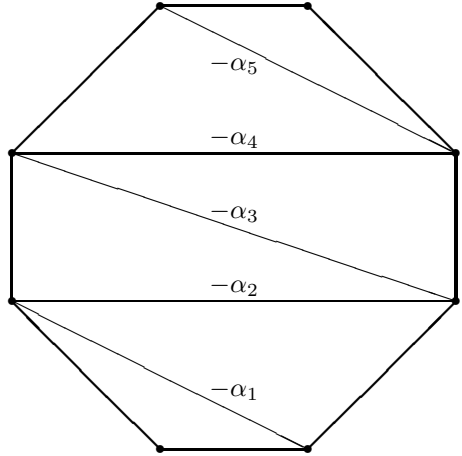
\begin{figure}[ht] 
\begin{center} 
\setlength{\unitlength}{2.8pt} 
\begin{picture}(60,60)(0,0) 
\thicklines 
  \multiput(0,20)(60,0){2}{\line(0,1){20}} 
  \multiput(20,0)(0,60){2}{\line(1,0){20}} 
  \multiput(0,40)(40,-40){2}{\line(1,1){20}} 
  \multiput(20,0)(40,40){2}{\line(-1,1){20}} 
 
  \multiput(20,0)(20,0){2}{\circle*{1}} 
  \multiput(20,60)(20,0){2}{\circle*{1}} 
  \multiput(0,20)(0,20){2}{\circle*{1}} 
  \multiput(60,20)(0,20){2}{\circle*{1}} 
 
\thinlines \put(0,20){\line(1,0){60}} \put(0,40){\line(1,0){60}} 
\put(0,20){\line(2,-1){40}} \put(0,40){\line(3,-1){60}} 
\put(20,60){\line(2,-1){40}} 
 
\put(30,8){\makebox(0,0){$-\alpha_1$}} 
\put(30,22){\makebox(0,0){$-\alpha_2$}} 
\put(30,32){\makebox(0,0){$-\alpha_3$}} 
\put(30,42){\makebox(0,0){$-\alpha_4$}} 
\put(30,52){\makebox(0,0){$-\alpha_5$}} 
 

\end{picture} 
\end{center} 
\caption{The ``snake'' in type $A_5$} 
\label{fig:octagon-snake} 
\end{figure}

\begin{figure}[ht] 
\begin{center} 
\setlength{\unitlength}{6pt} 
\begin{picture}(40,32)(-4,0) 
  \put(6,0){\line(1,0){20}} 
  \put(0,19){\line(1,0){32}} 
  \qbezier(6,0)(11,15.5)(16,31) 
  \qbezier(26,0)(21,15.5)(16,31) 
  \qbezier(6,0)(3,9.5)(0,19) 
  \qbezier(26,0)(29,9.5)(32,19) 
  \qbezier(0,19)(13,9.5)(26,0) 
  \qbezier(32,19)(19,9.5)(6,0) 
  \qbezier(0,19)(8,25)(16,31) 
  \qbezier(32,19)(24,25)(16,31) 
 
  \put(6,0){\circle*{.5}} 
  \put(26,0){\circle*{.5}} 
  \put(0,19){\circle*{.5}} 
  \put(32,19){\circle*{.5}} 
  \put(16,31){\circle*{.5}} 
 
\put(16,20.5){\makebox(0,0){$\alpha_1+\alpha_2$}} 
\put(23,16){\makebox(0,0){$-\alpha_2$}} 
\put(9,16){\makebox(0,0){$-\alpha_1$}} 
\put(12,8.7){\makebox(0,0){$\alpha_1$}} 
\put(20,8.7){\makebox(0,0){$\alpha_2$}} 
 
\end{picture} 
\end{center} 
\caption{Labeling of the diagonals in type $A_2$} \label{fig:a2} 
\end{figure}
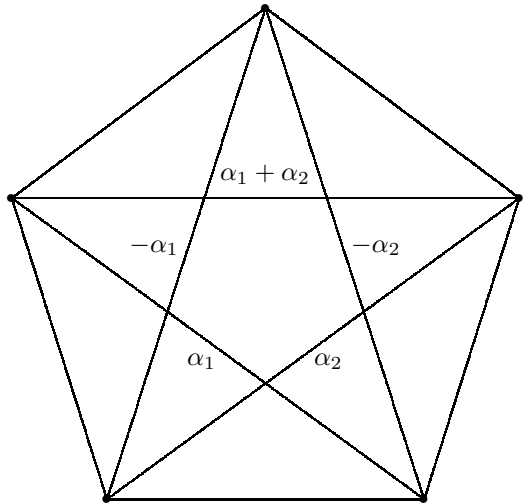

First, we reproduce the concrete description of the cluster complex of 
type~$A_n$ given in~\cite[Section~3.5]{yga}. 
We identify $\Phi_{\geq -1}$ with the set of all diagonals of a 
regular $(n+3)$-gon $\Poly_{n+3}$. Under this identification, the 
roots in $- \Pi$ correspond to the diagonals on the ``snake'' 
shown in Figure~\ref{fig:octagon-snake}. 
 Non-crossing diagonals represent compatible roots, while crossing 
diagonals correspond to roots 
whose compatibility degree is~1. 
(Here and in the sequel, two 
diagonals are called \emph{crossing} if they are distinct and have 
a common interior point.) 
Thus, each positive root $\alpha[i,j]= 
\alpha_i + \alpha_{i+1} + \cdots + \alpha_j$ corresponds to the 
unique diagonal that crosses precisely the diagonals $- \alpha_i, 
- \alpha_{i+1}, \ldots,- \alpha_j$ from the snake (see 
Figure~\ref{fig:a2}). 
 
 The clusters are in bijection with the triangulations of 
$\Poly_{n+3}$ by non-crossing diagonals. 
The cluster complex is the dual complex of the ordinary associahedron. 
Two triangulations are joined by an edge in the exchange graph if and 
only if they are obtained from each other by a ``flip'' that replaces 
a diagonal in a quadrilateral formed by two triangles of the 
triangulation by another diagonal of the same quadrilateral. 
See \cite[Section~3.5]{yga} and \cite[Section~4.1]{cfz} for further 
details. 
 
 
We next describe the cluster variables and the exchange relations in 
concrete combinatorial terms. 
For a diagonal $[a,b]$, we denote by $x_{ab}$ the cluster variable 
$x[\alpha]$ associated to the corresponding root. 
We adopt the convention that $x_{ab} = 1$ if $a$ and $b$ are two 
consecutive vertices of $\Poly_{n+3}$. 
Comparing (\ref{eq:b-matrix-via-roots-2}) 
with \cite[Lemma~4.2]{cfz}, we conclude that 
the matrices $B(C)$ can be described as follows. 
 
\begin{proposition} 
\label{pr:exchanges-A-1} 
Let $C$ be the cluster corresponding to  a triangulation $T$ of $\Poly_{n+3}$, 
and let $B(C) = B(T)$ be the corresponding matrix 
with rows and columns indexed by the diagonals in~$T$. 
Then each matrix entry $b_{\alpha \beta}$ is equal to $0$ 
unless $\alpha$ and $\beta$ are two sides of some triangle 
$(a,b,c)$ in $T$; in the latter case, $b_{\alpha \beta} = 1$ 
(resp.,~$-1$) if $\alpha = [a,b], \beta = [a,c]$, and the order of 
points $a,b,c$ is counter-clockwise (resp., clockwise). 
\end{proposition} 
 
\begin{figure}[ht] 
\begin{center} 
\setlength{\unitlength}{2pt} 
\begin{picture}(60,66)(0,-5) 
\thicklines 
  \put(0,20){\line(1,2){20}} 
  \put(0,20){\line(1,-1){20}} 
  \put(0,20){\line(3,1){60}} 
  \put(20,0){\line(0,1){60}} 
  \put(20,0){\line(1,1){40}} 
  \put(20,60){\line(2,-1){40}} 
 
  \put(20,0){\circle*{1}} 
  \put(20,60){\circle*{1}} 
  \put(0,20){\circle*{1}} 
  \put(60,40){\circle*{1}} 
 
\put(24,42){\makebox(0,0){$\beta$}} 
\put(35,35){\makebox(0,0){$\beta'$}}

\put(63,40){\makebox(0,0){$b$}} 
\put(20,63){\makebox(0,0){$c$}} 
\put(-3,20){\makebox(0,0){$d$}} 
\put(20,-3){\makebox(0,0){$a$}} 

\end{picture} 

\end{center} 
\caption{Exchanges in type $A_n$} 
\label{fig:type-a-exch} 
\end{figure}
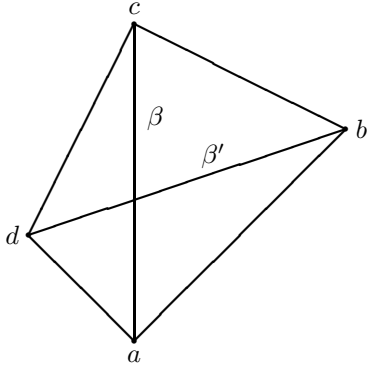

In view of Proposition~\ref{pr:exchanges-A-1}, 
the exchange relations (\ref{eq:exch-rel-fin-type}) in a cluster 
algebra of type $A_n$ have the form 
\begin{equation} 
\label{eq:exchange-relation-A} 
x_{ac} x_{bd} = p^+_{ac,bd}\, x_{ab}\, x_{cd} 
+ p^-_{ac,bd}\, x_{ad}\, x_{bc} \ , 
\end{equation} 
where $a,b,c,d$ are any four vertices of $\Poly_{n+3}$ taken in 
counter-clockwise order, and 
$p^\pm_{ac,bd}$ 
are elements of the coefficient semifield~$\PP$. 
See Figure~\ref{fig:type-a-exch}. 
 

Let $\AA_\circ$ denote the cluster algebra of type $A_n$ 
associated with the following coefficient system. 
We take 
\begin{equation} 
\label{coeff-special-A} 
\PP = \Trop (p_{ab}: \text{$[a,b]$ is a side of $\Poly_{n+3}$}) 
\,, 
\end{equation} 
and define the coefficients in (\ref{eq:exchange-relation-A}) by 
\begin{equation} 
\label{coeff-special-A-2} 
p^+_{ac,bd} = q_{ab}\, q_{cd}, \quad 
p^-_{ac,bd} = q_{ad}\, q_{bc} \, , 
\end{equation} 
where 
\begin{equation} 
\label{coeff-special-A-3} 
q_{ab} = 
\begin{cases} 
1 & \text{if $[a,b]$ is a diagonal};\\ 
p_{ab} & \text{if $[a,b]$ is a side}. 
\end{cases} 
\end{equation} 
Since $p^+_{ac,bd}$ and $p^-_{ac,bd}$ have no common factors, 
they satisfy the normalization condition 
$p^+_{ac,bd} \oplus p^-_{ac,bd} = 1$; 
a direct check using Proposition~\ref{pr:exchanges-A-1} shows that 
this choice of coefficients also satisfies the mutation rule 
(\ref{eq:p-mutation}), making the cluster algebra $\AA_\circ$ well-defined. 
 
\begin{example} 
[\emph{Geometric realization for $\AA_\circ$ in type~$A_n$}] 
\label{ex:realization-An} 
Let $X = X_{n+3}$ 
be the affine cone over the Grassmannian ${\rm Gr}_{2,n+3}$ 
of 2-dimensional subspaces in $\CC^{n+3}$ 
taken in its Pl\"ucker embedding (cf.~Example~\ref{ex:Gr24-A1}); 
simply put, $X$ is the variety of all nonzero decomposable 
bivectors in $\bigwedge^2 \CC^{n+3}$. 
Let $(\Delta_{ab}: 1 \leq a < b \leq n+3)$ be the standard Pl\"ucker 
coordinates on~$X$. 
We identify the indices $1, \dots, n+3$ with the vertices of 
$\Poly_{n+3}$ by numbering these vertices, say, counterclockwise. 
Thus, we associate the Pl\"ucker coordinates with all the sides and 
diagonals of $\Poly_{n+3}$. 
Note that we have previously used the same set 
$\{(ab): 1 \leq a < b \leq n+3)\}$ to label the 
cluster variables $x_{ab}$ and the coefficients~$p_{ab}\,$. 
\end{example} 
 
\begin{proposition} 
\label{pr:geom-realization-An} 
The correspondence sending each cluster variable $x_{ab}$ 
and each coefficient $p_{ab}$ to the corresponding element $\Delta_{ab}$ 
extends uniquely to an algebra isomorphism between the cluster algebra $\AA_\circ$ 
and the $\ZZ$-form $\ZZ[X]$ of $\CC[X]$ generated by all Pl\"ucker 
coordinates. 
\end{proposition} 
 
\proof 
This is a special case of Proposition~\ref{pr:geom-realization-criterion}. 
To see this, we need to verify the conditions 
(\ref{eq:X-rational-etc})--(\ref{eq:same-relations}). 
The fact that $X$ satisfies (\ref{eq:X-rational-etc}) is well 
known (for the rationality property, it is enough to note that 
$X$ has a Zariski open subset isomorphic to an affine space). 
For the dimension count (\ref{eq:X-dimension}), we have 
$$\dim (X) = \dim ({\rm Gr}_{2,n+3}) + 1 = 2n+3 = n + |J|,$$ 
as required. 
The property (\ref{eq:phi-generate-X}) means that 
$\CC[X]$ is generated by all Pl\"ucker coordinates, which is 
trivial. 
Finally, (\ref{eq:same-relations}) follows from the standard fact that 
the Pl\"ucker coordinates satisfy the Grassmann-Pl\"ucker relations 
\[ 
\Delta_{ac} \Delta_{bd} = \Delta_{ab}\, \Delta_{cd} + 
\Delta_{ad}\, \Delta_{bc} 
\] 
for all $1 \leq a < b < c < d \leq n+3$. 
\endproof 
 
We note that the ring $\ZZ[X]$ is naturally identified 
with the ring of $SL_2$-invariant polynomial functions 
with coefficients in $\ZZ$ on the space of $(n+3)$-tuples of vectors in~$\CC^2$. 
Representing these vectors as columns of a $2 \times (n+3)$ 
matrix $Z = (z_{ij})$, we identify the Pl\"ucker coordinates with 
the $2 \times 2$ minors of $Z$: 
$$\Delta_{ab} = z_{1a} z_{2b} - z_{1b} z_{2a} \quad (1 \leq a < b \leq n+3).$$ 
 
\begin{remark} 
\label{rem:compatible-monomial-basis} 
It is classically known that the monomials in the Pl\"ucker coordinates 
that are not divisible by $\Delta_{ac} \Delta_{bd}$ for 
any $a < b < c <d$, form a $\ZZ$-basis in $\ZZ[X]$ 
(see~\cite{kungrota84} or~\cite{stur} for a proof). 
Let us translate this fact into the setting of cluster algebras. 
We shall call a monomial $\prod_\alpha x_\alpha^{m_\alpha}$ in the 
cluster variables \emph{compatible} if $m_\alpha m_\beta = 0$ 
whenever the roots $\alpha$ and $\beta$ are incompatible, i.e., 
whenever the corresponding diagonals cross each other. 
(Equivalently, all variables contributing to a compatible monomial 
belong to a single cluster.) 
In this terminology, the cluster algebra $\AA_\circ$ is a free 
$\ZZ[\mathcal P]$-module with the basis formed by all compatible monomials. 
We believe that this property remains true for an arbitrary cluster 
algebra of finite type (we have checked it for all classical 
types); we plan to investigate it in a separate publication. 
We note that linear independence of compatible monomials is an 
immediate consequence of Theorem~\ref{th:cluster-variable-denominators-1} 
and the uniqueness of cluster expansions 
(Proposition~\ref{prop:inv-cluster-expansion-1}). 
\end{remark}

\subsection{Types~$B_n$ and $C_n$} 
\label{sec:ca-type-C} 
Let $\Phi$ be a root system of type $B_n$ or~$C_n$. 
We identify the set $I$ in a standard way with $[1,n]$. 
As in \cite[Section~4.2]{cfz}, in order to treat both cases at the same time, 
we set $r =1$ for 
$\Phi$ of type $B_n\,$, and $r = 2$ for 
$\Phi$ of type $C_n\,$. 
Once again, our convention for the Cartan matrices is different from 
the one in \cite{bourbaki} but agrees with that in~\cite{kac}: 
thus, we have $a_{n-1,n} = -r$ and $a_{n,n-1} = -2/r$.

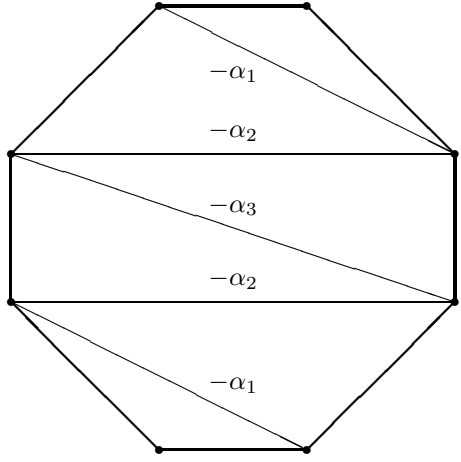
\begin{figure}[ht] 
\begin{center} 
\setlength{\unitlength}{2.8pt} 
\begin{picture}(60,60)(0,0) 
\thicklines 
  \multiput(0,20)(60,0){2}{\line(0,1){20}} 
  \multiput(20,0)(0,60){2}{\line(1,0){20}} 
  \multiput(0,40)(40,-40){2}{\line(1,1){20}} 
  \multiput(20,0)(40,40){2}{\line(-1,1){20}} 
 
  \multiput(20,0)(20,0){2}{\circle*{1}} 
  \multiput(20,60)(20,0){2}{\circle*{1}} 
  \multiput(0,20)(0,20){2}{\circle*{1}} 
  \multiput(60,20)(0,20){2}{\circle*{1}} 
 
\thinlines \put(0,20){\line(1,0){60}} \put(0,40){\line(1,0){60}} 
\put(0,20){\line(2,-1){40}} \put(0,40){\line(3,-1){60}} 
\put(20,60){\line(2,-1){40}} 
 
\put(30,9){\makebox(0,0){$-\alpha_1$}} 
\put(30,23){\makebox(0,0){$-\alpha_2$}} 
\put(30,33){\makebox(0,0){$-\alpha_3$}} 
\put(30,43){\makebox(0,0){$-\alpha_2$}} 
\put(30,51){\makebox(0,0){$-\alpha_1$}} 
 
\end{picture} 
 
\end{center} 
\caption{The ``snake'' for the types $B_3$ and $C_3$} 
\label{fig:b-octagon} 
\end{figure}

We recall the combinatorial description of the cluster complex of 
type~$B_n/C_n$ from~\cite[Section~3.5]{yga}. 
Let $\Theta$ denote the $180^\circ$ rotation of a regular 
$(2n+2)$-gon $\Poly_{2n+2}$. There is a natural action of $\Theta$ 
on the diagonals of $\Poly_{2n+2}$. Each orbit of this action is 
either a diameter (i.e., a diagonal connecting antipodal vertices) 
or an unordered pair of centrally symmetric non-diameter diagonals 
of $\Poly_{2n+2}$. 
Following \cite[Section~3.5]{yga}, we identify almost positive roots 
in $\Phi$ 
with these orbits. 
Under this identification, each of the roots $- 
\alpha_i$ for $i = 1, \dots, n-1$ is represented by a pair of 
diagonals on the ``snake'' shown in Figure~\ref{fig:b-octagon}, 
whereas $- \alpha_n$ is identified with the only diameter on the snake. 
Two $\Theta$-orbits represent compatible roots if and only 
if the diagonals they involve do not cross each other. 
More generally, in type $B_n$ (resp.,~$C_n$), 
for $\alpha, \beta \in \Phi_{\geq -1}$, the 
compatibility degree $(\alpha \| \beta)$ is equal to the number of 
crossings of one of the diagonals representing $\alpha$ (resp.,~$\beta$) 
by the diagonals representing~$\beta$ (resp.,~$\alpha$). 
Thus, each positive root $\beta=\sum_{i} b_i \alpha_i$ is represented by the unique 
$\Theta$-orbit such that every diagonal representing $-\alpha_i$ 
(resp.,~$\beta$) crosses the diagonals representing $\beta$ (resp.,~$-\alpha_i$) 
at $b_i$ points.

The clusters are in bijection with the centrally-symmetric (that is, $\Theta$-invariant) 
triangulations of $\Poly_{2n+2}$ by non-crossing diagonals. 
The cluster complex is the dual complex for the Bott-Taubes 
cyclohedron~\cite{bott-taubes}. 
Two centrally symmetric triangulations are joined by an edge in 
the exchange graph if and only if they are obtained from each other either by a flip 
involving two diameters, or by a pair of centrally symmetric flips. 
See \cite[Section~3.5]{yga} and \cite[Section~4.2]{cfz} for 
details. 
 
For a vertex $a$ of $\Poly_{2n+2}$, let $\overline a$ denote the 
antipodal vertex $\Theta (a)$. 
For a diagonal $[a,b]$, we denote 
by $x_{ab}$ the cluster variable $x[\alpha]$ associated to the 
root corresponding to the $\Theta$-orbit of $[a,b]$. 
Thus, we have $x_{ab} = x_{ba} = x_{\overline a \,\overline b}$. 
Similarly to the type $A_n$, we adopt the convention that $x_{ab} = 1$ if $a$ 
and $b$ are consecutive vertices in~$\Poly_{2n+2}$. 
 
Comparing (\ref{eq:exch-rel-fin-type}) with \cite[Lemma~4.4]{cfz}, 
we obtain the following concrete description of the exchange relations in 
types $B_n$ and~$C_n\,$. 
 
\begin{figure}[ht] 
\begin{center} 
\setlength{\unitlength}{2pt} 
\begin{picture}(60,66)(0,-2) 
\thicklines 
  \multiput(0,20)(60,0){2}{\line(0,1){20}} 
  \multiput(20,0)(0,60){2}{\line(1,0){20}} 
  \multiput(0,40)(40,-40){2}{\line(1,1){20}} 
  \multiput(20,0)(40,40){2}{\line(-1,1){20}} 
 
  \multiput(20,0)(20,0){2}{\circle*{1}} 
  \multiput(20,60)(20,0){2}{\circle*{1}} 
  \multiput(0,20)(0,20){2}{\circle*{1}} 
  \multiput(60,20)(0,20){2}{\circle*{1}} 
 
\thinlines 
\put(20,60){\line(2,-1){40}} 
\put(20,60){\line(1,-1){40}} 
\put(40,60){\line(1,-2){20}} 
 
\put(42,62){\makebox(0,0){$c$}} 
\put(64,20){\makebox(0,0){$a$}} 
\put(64,40){\makebox(0,0){$b$}} 
\put(18,62){\makebox(0,0){$d$}} 
\put(-4,40){\makebox(0,0){$\overline a$}} 

 
\multiput(20,0)(-1,2){20}{\circle*{0.5}} 
\multiput(0,20)(2,-1){20}{\circle*{0.5}} 
\multiput(0,40)(2,-2){20}{\circle*{0.5}} 
 
\put(30,-7){\makebox(0,0){(i)}} 
 
\end{picture} 
\qquad\qquad\qquad 
\begin{picture}(60,66)(0,-2) 
\thicklines 
  \multiput(0,20)(60,0){2}{\line(0,1){20}} 
  \multiput(20,0)(0,60){2}{\line(1,0){20}} 
  \multiput(0,40)(40,-40){2}{\line(1,1){20}} 
  \multiput(20,0)(40,40){2}{\line(-1,1){20}} 
 
  \multiput(20,0)(20,0){2}{\circle*{1}} 
  \multiput(20,60)(20,0){2}{\circle*{1}} 
  \multiput(0,20)(0,20){2}{\circle*{1}} 
  \multiput(60,20)(0,20){2}{\circle*{1}} 
 
\thinlines 
 
\put(40,0){\line(0,1){60}} 
\put(20,60){\line(2,-1){40}} 
\put(40,0){\line(1,2){20}}

\put(20,60){\line(1,-3){20}} 

\put(42,-2){\makebox(0,0){$a$}} 
\put(42,62){\makebox(0,0){$c$}} 
\put(64,40){\makebox(0,0){$b$}} 
\put(18,62){\makebox(0,0){$\overline a$}} 
 
 
\multiput(20,0)(0,2){30}{\circle*{0.5}} 
\multiput(0,20)(2,-1){20}{\circle*{0.5}} 
 
\put(30,-7){\makebox(0,0){(ii)}} 
 
\end{picture} 
\\[.4in] 
\begin{picture}(60,66)(0,-7) 
\thicklines 
  \multiput(0,20)(60,0){2}{\line(0,1){20}} 
  \multiput(20,0)(0,60){2}{\line(1,0){20}} 
  \multiput(0,40)(40,-40){2}{\line(1,1){20}} 
  \multiput(20,0)(40,40){2}{\line(-1,1){20}} 
 
  \multiput(20,0)(20,0){2}{\circle*{1}} 
  \multiput(20,60)(20,0){2}{\circle*{1}} 
  \multiput(0,20)(0,20){2}{\circle*{1}} 
  \multiput(60,20)(0,20){2}{\circle*{1}} 
 
\thinlines 
\put(20,60){\line(2,-1){40}} 
\put(0,20){\line(2,-1){40}} 
\put(20,60){\line(1,-3){20}} 
\put(40,0){\line(1,2){20}} 
\put(0,20){\line(1,2){20}} 
\put(0,20){\line(3,1){60}} 
 
\put(42,-2){\makebox(0,0){$a$}} 
\put(-4,20){\makebox(0,0){$\overline b$}} 
\put(64,40){\makebox(0,0){$b$}} 
\put(18,62){\makebox(0,0){$\overline a$}} 
 
 
\put(30,-7){\makebox(0,0){(iii)}} 
 
\end{picture} 
 
\end{center} 
\caption{
Exchanges in types $B_n$ and $C_n$} 
\label{fig:quadrilateral-D} 
\end{figure}
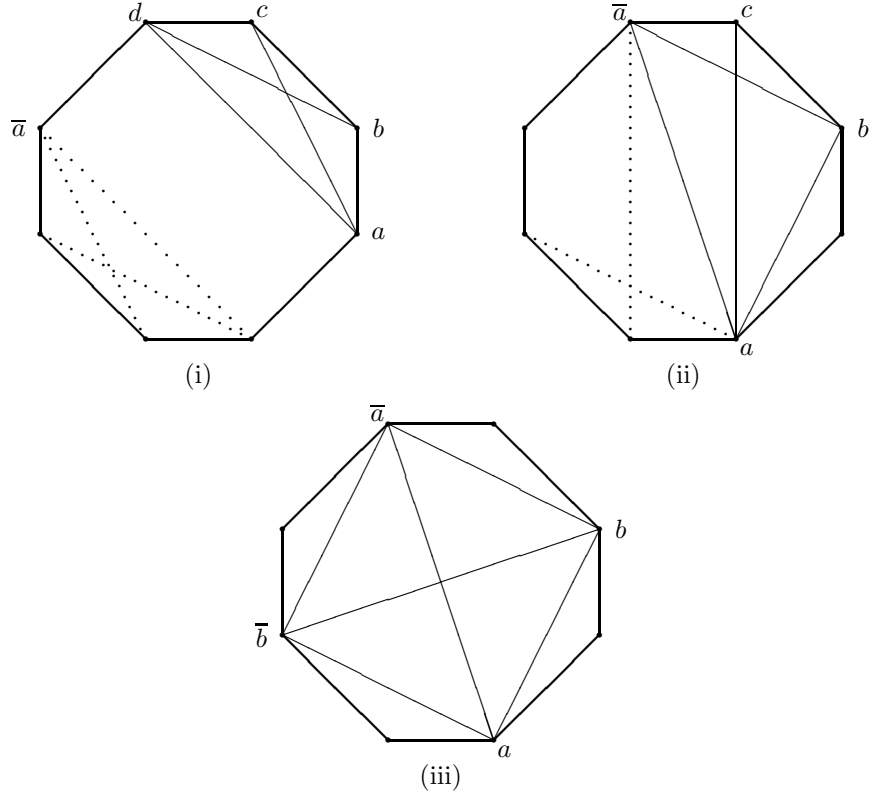

\begin{proposition} 
\label{pr:exchanges-C} 
The exchange relations in a cluster algebra of type $B_n$ or $C_n$ 
have the following form: 
\begin{equation} 
\label{eq:exchange-relation-C-1} 
x_{ac} x_{bd} = 
p^+_{ac,bd}\, 
x_{ab}\, x_{cd} 
+ 
p^-_{ac,bd}\, 
x_{ad}\, x_{bc} \ , 
\end{equation} 
whenever $a,b,c,d,\overline a$ are 
in counter-clockwise order; 
\begin{equation} 
\label{eq:exchange-relation-C-2} 
x_{ac} x_{a \overline b} = 
p^+_{ac,a \overline b}\, 
x_{ab}\, x_{a \overline c} 
+ 
p^-_{ac,a\overline b}\, 
x_{a \overline a}^{2/r} \, x_{bc} \ , 
\end{equation} 
whenever $a,b,c,\overline a$ are 
in counter-clockwise order; 
\begin{equation} 
\label{eq:exchange-relation-C-3} 
x_{a\overline a} x_{b \overline b} = 
p^+_{a\overline a,b \overline b}\, 
x_{ab}^r 
+ 
p^-_{a\overline a,b\overline b}\, 
x_{a \overline b}^r \ , 
\end{equation} 
whenever $a,b,\overline a$ are 
in counter-clockwise order. 
See Figure~\ref{fig:quadrilateral-D}. 
\end{proposition} 
 
We will provide an explicit realization of a special cluster 
algebra $\AA_\circ$ of type $B_n$ or $C_n$ 
similar to its namesake for $A_n$. 
The algebra $\AA_\circ$ corresponds to the following special choice of coefficients. 
We set $\PP = \Trop (\{p_\delta\})$, where 
$\delta$ runs over all centrally-symmetric pairs of sides of the polygon~$\Poly_{2n+2}\,$. 
For such a pair $\delta = \{[a,b], [\overline a, \overline b])\}$, 
we write the corresponding generator of $\PP$ as 
$p_\delta = p_{ab} = p_{\overline a \, \overline b}$. 
The coefficients in 
(\ref{eq:exchange-relation-C-1})--(\ref{eq:exchange-relation-C-3}) 
are specified in a similar way to 
(\ref{coeff-special-A-2})--(\ref{coeff-special-A-3}). 
More precisely, to obtain a coefficient of some monomial in 
(\ref{eq:exchange-relation-C-1})--(\ref{eq:exchange-relation-C-3}), 
take this monomial and replace each of its cluster variables 
$x_{ab}$ by $1$ (resp.,~$p_{ab}$) if $[a,b]$ is a diagonal 
(resp.,~a side) of $\Poly_{2n+2}$. 
The fact that these coefficients satisfy the normalization 
condition is again obvious; we leave to the reader a direct check 
that they also satisfy the mutation rule (\ref{eq:p-mutation}). 
 
So far the material for $B_n$ has been completely parallel to 
that for $C_n$. 
However, our geometric realizations for these two types are quite different from 
each other. 
 
\begin{example} 
[\emph{Geometric realization for $\AA_\circ$ in type~$B_n$}] 
\label{ex:realization-Bn} 
Somewhat surprising\-ly, it turns out that the algebra 
$\AA_\circ$ for the type $B_n$ (for $n\geq 3$) is isomorphic, as a 
ring, to the cluster algebra $\AA_\circ$ for the type $A_{n-1}$. 
Recall that the latter is naturally identified with 
the ring $\ZZ[X_{n+2}]$ generated by the Pl\"ucker coordinates 
$\Delta_{ab}$ (for $1 \leq a < b \leq n+2$) on the Grassmannian ${\rm Gr}_{2,n+2}$ 
(see Proposition~\ref{pr:geom-realization-An}). 
 
Let us label the vertices of 
$\Poly_{2n+2}$ in the counterclockwise order by the indices 
$1, \dots, n+1, \overline 1, \dots, \overline {n+1}$. 
We associate a function from $\ZZ[X_{n+2}]$ to every 
$\Theta$-orbit on the set of all diagonals and sides of 
$\Poly_{2n+2}\,$, as follows: 
\begin{equation} 
\label{eq:Delta-functions-Bn} 
\begin{array}{rcll} 
[a, \overline a]  &\mapsto& \Delta_{a \overline a} = \Delta_{a,n+2} 
  &\quad (1 \leq a \leq n+1),\\[.05in] 
\{[a,b], [\overline a, \overline b]\} &\mapsto& \Delta_{ab} 
  &\quad (1 \leq a < b \leq n+1),\\[.05in] 
\{[a,\overline b], [\overline a, b]\} &\mapsto& 
\Delta_{a \overline b} = 
\Delta_{a,n+2} \Delta_{b,n+2} - \Delta_{ab} 
  &\quad (1 \leq a < b \leq n+1). 
\end{array} 
\end{equation} 
 
\begin{proposition} 
\label{pr:geom-realization-Bn} 
The correspondence sending each cluster variable 
and each coefficient for the cluster algebra $\AA_\circ$ 
of type $B_n$ to the element in {\rm (\ref{eq:Delta-functions-Bn})} 
with the same label extends uniquely to an algebra isomorphism of 
$\AA_\circ$ with $\ZZ[X_{n+2}]$. 
\end{proposition} 
 
\proof 
The proof is similar to that of 
Proposition~\ref{pr:geom-realization-An}. 
It is enough to check that our data satisfy the conditions 
(\ref{eq:X-rational-etc})--(\ref{eq:same-relations}). 
Condition (\ref{eq:X-rational-etc}) was already checked in 
the proof  of Proposition~\ref{pr:geom-realization-An}. 
The dimension count (\ref{eq:X-dimension}) now takes the form: 
$$\dim (X_{n+2}) = 2n+1 = n + |J|,$$ 
as required. 
The property (\ref{eq:phi-generate-X}) is clear since 
all the Pl\"ucker coordinates $\Delta_{ab}$ for $1 \leq a < b \leq n+2$ 
are among the functions (\ref{eq:Delta-functions-Bn}). 
Finally, (\ref{eq:same-relations}) amounts to checking the 
following six identities (with $r=1$) obtained from the exchange relations 
(\ref{eq:exchange-relation-C-1})--(\ref{eq:exchange-relation-C-3}) 
(we have to take into account possible positions of vertices 
in Figure~\ref{fig:quadrilateral-D} among the vertices 
$1, \dots, n+1, \overline 1, \dots, \overline {n+1}$): 
\begin{align} 
\label{eq:exchange-relations-in-deltas-Bn-1} 
\Delta_{ac} \, \Delta_{bd} = 
\Delta_{ab}\, \Delta_{cd} + 
\Delta_{ad}\, \Delta_{bc}  \qquad & (1 \leq a < b < c < d \leq n+1),\\ 
\label{eq:exchange-relations-in-deltas-Bn-2} 
\Delta_{a \overline c} \, \Delta_{bd} = 
\Delta_{a \overline b}\, \Delta_{cd} + 
\Delta_{a \overline d}\, \Delta_{bc} \qquad & (1 \leq a < b < c < d \leq n+1),\\ 
\label{eq:exchange-relations-in-deltas-Bn-3} 
\Delta_{a \overline c} \, \Delta_{b \overline d} = 
\Delta_{ab}\, \Delta_{cd} + 
\Delta_{a \overline d}\, \Delta_{b \overline c} \qquad & (1 \leq a < b 
< c < d \leq n+1),\\ 
\label{eq:exchange-relations-in-deltas-Bn-4} 
\Delta_{ac} \, \Delta_{a \overline b} = 
\Delta_{ab}\, \Delta_{a \overline c} + 
\Delta_{a \overline a}^{2/r} \, \Delta_{bc} \qquad & (1 \leq a < b < c \leq n+1),\\ 
\label{eq:exchange-relations-in-deltas-Bn-5} 
\Delta_{a \overline b} \, \Delta_{b \overline c} = 
\Delta_{ab}\, \Delta_{bc} + 
\Delta_{b \overline b}^{2/r} \, \Delta_{a \overline c} \qquad & 
(1 \leq a < b < c \leq n+1),\\ 
\label{eq:exchange-relations-in-deltas-Bn-6} 
\Delta_{a\overline a}\, \Delta_{b \overline b} = 
\Delta_{ab}^r + \Delta_{a \overline b}^r \qquad\qquad &(1 \leq a < b \leq n+1) \, . 
\end{align} 
Of these identities, (\ref{eq:exchange-relations-in-deltas-Bn-1}) is 
a Grassmann-Pl\"ucker relation, and the rest are reduced 
to this relation by simple algebraic manipulations. 
\endproof 
 
\end{example} 
 
\begin{example} 
[\emph{Geometric realization for $\AA_\circ$ in type~$C_n$}] 
\label{ex:realization-Cn} 
Let $SO_2$ be the group of complex matrices 
\[ 
\bmat{u}{-v}{v}{u} 
\] 
with $u^2+v^2 = 1$. 
Consider the algebra $\mathcal{R} = \CC[{\rm Mat}_{2,n+1}]^{SO_2}$ 
of $SO_2$-invariant polynomial functions 
on the space of $2 \times (n+1)$ complex matrices, 
or, equivalently, on the space of 
$(n+1)$-tuples of vectors in~$\CC^2$. 
Alternatively, $\mathcal{R}$ can be identified with the ring 
of invariants 
$\CC[{\rm Mat}_{2,n+1}]^T$, where $T\subset SL_2$ is the torus of all 
diagonal matrices of the form 
\[ 
\bmat{t}{0}{0}{t^{-1}}\,. 
\] 
Indeed, we have $g(SO_2) g^{-1} = T$, where 
\[ 
g = \bmat{1}{-i}{1}{i}\,, 
\] 
so the map $f \mapsto f^g$  defined by 
$f^g(z) = f(gz)$ is an isomorphism 
\begin{equation} 
\label{eq:SO2=T} 
\CC[{\rm Mat}_{2,n+1}]^T \to \CC[{\rm Mat}_{2,n+1}]^{SO_2}. 
\end{equation} 
 
The ring $\mathcal{R}=\CC[{\rm Mat}_{2,n+1}]^T$ can 
also be viewed as the coordinate ring $\CC[X]$ of the variety 
$X$ of complex $(n+1) \times (n+1)$ matrices of rank $\leq 1$ 
(even more geometrically, $X - \{0\}$ is the affine cone over the 
product of two copies of the projective space $\mathbb{CP}^n$ taken in the 
Segre embedding). 
Specifically, the map 
\[ 
y=\left[\!\!\begin{array}{ccc} 
y_{11} & \cdots & y_{1,n+1}\\ 
y_{21} & \cdots & y_{2,n+1}\\ 
\end{array}\!\!\right] \,\, 
\mapsto \,\, (y_{1a}y_{2b})_{a,b= 1, \dots, n+1} \in X 
\] 
induces an algebra isomorphism 
$\CC[X] \to \CC[{\rm Mat}_{2,n+1}]^T$. 
Combining this with (\ref{eq:SO2=T}), we obtain an isomorphism 
$\CC[X] \to \CC[{\rm Mat}_{2,n+1}]^{SO_2}$ induced by the map 
\[ 
\left[\!\!\begin{array}{ccc} 
z_{11} & \cdots & z_{1,n+1}\\ 
z_{21} & \cdots & z_{2,n+1}\\ 
\end{array}\!\!\right] \,\, 
\mapsto \,\, ((z_{1a} - i z_{2a})(z_{1b} + i z_{2b}))_{a,b= 1, 
\dots, n+1} \in X \,. 
\] 

By analogy with (\ref{eq:Delta-functions-Bn}), 
we associate an element from 
$\mathcal{R} = \CC[X]$ to every 
$\Theta$-orbit on the set of all diagonals and sides of 
$\Poly_{2n+2}\,$, as follows: 
\begin{equation} 
\label{eq:Delta-functions-Cn} 
\begin{array}{ll} 
\{[a,b], [\overline a, \overline b]\} \mapsto \Delta_{ab} 
= z_{1a} z_{2b} - z_{1b} z_{2a} 
= \displaystyle\frac{y_{1a}y_{2b} - y_{1b}y_{2a}}{2i} 
& (1 \leq a < b \leq n+1)\,,\\[.1in] 
\{[a,\overline b], [\overline a, b]\} \mapsto 
\Delta_{a \overline b} 
= z_{1a} z_{1b} + z_{2a} z_{2b} 
= \displaystyle\frac{y_{1a}y_{2b} + y_{1b}y_{2a}}{2} 
& (1 \leq a \leq b \leq n+1) \, . 
\end{array} 
\end{equation} 
The following result is a type $C_n$ counterpart of 
Proposition~\ref{pr:geom-realization-Bn}. 
\end{example} 
 
\begin{proposition} 
\label{pr:geom-realization-Cn} 
The correspondence sending each cluster variable 
and each coefficient for the cluster algebra 
$\AA_\circ$ of type $C_n$ to the element in {\rm (\ref{eq:Delta-functions-Cn})} 
with the same label extends uniquely to 
an algebra isomorphism of $\AA_\circ$ with a $\ZZ$-form of $\mathcal{R}$. 
\end{proposition} 
 
\proof 
The proof is analogous to that of 
Proposition~\ref{pr:geom-realization-Bn}. 
The only work involved is to check that the functions in 
(\ref{eq:Delta-functions-Cn}) satisfy the identities 
(\ref{eq:exchange-relations-in-deltas-Bn-1})--(\ref{eq:exchange-relations-in-deltas-Bn-6}), 
this time with $r= 2$. 
This is completely straightforward; for example, 
(\ref{eq:exchange-relations-in-deltas-Bn-6}) becomes 
\[ 
(z_{1a}^2 + z_{2a}^2)(z_{1b}^2 + z_{2b}^2) = 
(z_{1a} z_{2b} - z_{1b} z_{2a})^2 + 
(z_{1a} z_{1b} + z_{2a}z_{2b})^2 
.\hfil \qed 
\]

\subsection{Type~$D_n$} 
\label{sec:ca-type-D} 
\begin{figure}[ht] 
\begin{center} 
\setlength{\unitlength}{2.8pt} 
\begin{picture}(60,65)(0,-2) 
\thicklines 
  \multiput(0,20)(60,0){2}{\line(0,1){20}} 
  \multiput(20,0)(0,60){2}{\line(1,0){20}} 
  \multiput(0,40)(40,-40){2}{\line(1,1){20}} 
  \multiput(20,0)(40,40){2}{\line(-1,1){20}} 
 
  \multiput(20,0)(20,0){2}{\circle*{1}} 
  \multiput(20,60)(20,0){2}{\circle*{1}} 
  \multiput(0,20)(0,20){2}{\circle*{1}} 
  \multiput(60,20)(0,20){2}{\circle*{1}} 
 
\thinlines 
\put(0,20){\line(1,0){60}} 
\put(0,40){\line(1,0){60}} 
\put(0,20){\line(2,-1){40}} 
\put(1,40){\line(3,-1){59}} 
\put(0,39.5){\line(3,-1){4}} 
\put(6,37.5){\line(3,-1){4}} 
\put(12,35.5){\line(3,-1){4}} 
\put(18,33.5){\line(3,-1){4}} 
\put(24,31.5){\line(3,-1){4}} 
\put(30,29.5){\line(3,-1){4}} 
\put(36,27.5){\line(3,-1){4}} 
\put(42,25.5){\line(3,-1){4}} 
\put(48,23.5){\line(3,-1){4}} 
\put(54,21.5){\line(3,-1){4}} 
\put(20,60){\line(2,-1){40}} 
 
\put(30,8){\makebox(0,0){$-\alpha_1$}} 
\put(30,22){\makebox(0,0){$-\alpha_2$}} 
\put(30,32.5){\makebox(0,0){$-\alpha_3$}} 
\put(30,27.5){\makebox(0,0){$-\alpha_4$}} 
\put(30,42){\makebox(0,0){$-\alpha_2$}} 
\put(30,52){\makebox(0,0){$-\alpha_1$}} 
 
\end{picture} 
\end{center} 
\caption{Representing the roots in $-\Pi$ for the type $D_4$} 
\label{fig:d-octagon} 
\end{figure}
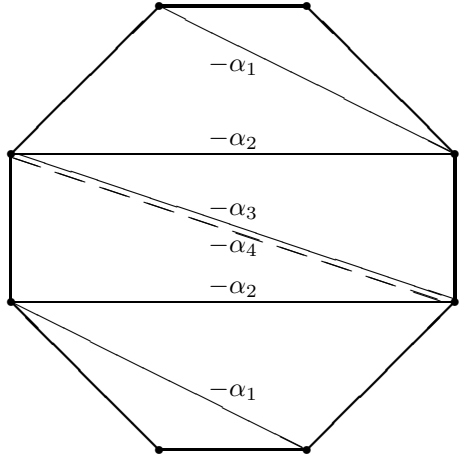

Let $\Phi$ be a root system of type $D_n$ 
(allowing for $n=3$, in which case $\Phi$ is of type~$A_3$). 
According to \cite[Section~3.5]{yga}, the almost positive roots 
(hence the cluster variables) for the type $D_n$ 
have a natural surjection onto those for~$B_{n-1}\,$. 
This surjection is one-to-one over the roots 
corresponding to pairs of diagonals of 
$\Poly_{2n}\,$, and two-to-one over those corresponding to diameters. 
Thus, the roots in $\Phi_{\geq -1}$ are represented 
by $\Theta$-orbits on the set of diagonals in a regular $2n$-gon, 
in which each diameter can be of one of two different ``colors''; 
we denote the two different kinds of diameters by 
$[a, \overline a]$ and $\widetilde{[a, \overline a]}$. 
The negative simple roots form a ``type~$D$ snake'' shown in 
 Figure~\ref{fig:d-octagon}. 
Two $\Theta$-orbits represent compatible roots if and 
only if the diagonals they involve do not cross each other; 
here we use the following convention: 
\begin{equation} 
\label{eq:diam-same-color} 
\text{diameters of the same color do not cross each other.} 
\end{equation} 
More generally, for $\alpha, \beta \in \Phi_{\geq -1}\,$, the 
compatibility degree 
$(\alpha \| \beta)$ is equal to the number of $\Theta$-orbits in the 
set of crossing points between the diagonals representing~$\alpha$ 
and~$\beta$ (again, with the convention (\ref{eq:diam-same-color})). 
Each positive root $\beta=\sum_{i} b_i \alpha_i$ 
is then represented by the unique 
$\Theta$-orbit such that the diagonals representing~$\beta$ 
cross the diagonals representing $-\alpha_i$ at $b_i$ pairs of 
centrally symmetric points (counting an intersection of two diameters 
of different color and location as one such pair). 
 
Accordingly, the cluster variables for $D_n$ can be denoted as 
$x_\alpha\,$, for all diagonals $\alpha$ in $\Poly_{2n}$ (with the convention 
$x_\alpha = x_{\Theta (\alpha)}$), plus $n$ extra variables 
$\tilde x_\beta$ for all diameters~$\tilde \beta$. 
 
With the help of \cite[Lemma~4.6]{cfz}, we obtain the 
following analogue of Proposition~\ref{pr:exchanges-C}.

\begin{proposition} 
\label{pr:exchanges-D} 
The exchange relations in a cluster algebra of type $D_n$ 
have the following form: 
\begin{equation} 
\label{eq:exchange-relation-D-1} 
x_{ac} x_{bd} = 
p^+_{ac,bd}\, 
x_{ab}\, x_{cd} 
+ 
p^-_{ac,bd}\, 
x_{ad}\, x_{bc} 
\end{equation} 
whenever $a,b,c,d,\overline a$ are 
in counter-clockwise order; 
\begin{equation} 
\label{eq:exchange-relation-D-2} 
x_{ac} x_{a \overline b} = 
p^+_{ac,a \overline b}\, 
x_{ab}\, x_{a \overline c} 
+ 
p^-_{ac,a\overline b}\, 
x_{a \overline a}  \, \tilde x_{a \overline a} \, x_{bc} 
\end{equation} 
whenever $a,b,c,\overline a$ are 
in counter-clockwise order; 
\begin{equation} 
\label{eq:exchange-relation-D-3} 
x_{a\overline a} \tilde x_{b \overline b} = 
p^+_{a\overline a,b \overline b}\, x_{ab} 
+ 
p^-_{a\overline a,b\overline b}\, 
x_{a \overline b} 
\end{equation} 
whenever $a,b,\overline a$ are 
in counter-clockwise order; 
\begin{equation} 
\label{eq:exchange-relation-D-4} 
x_{a \overline a} x_{b \overline c} = 
p^+_{a \overline a,b \overline c}\, 
x_{ab}\, x_{c \overline c} 
+ 
p^-_{a \overline a,b \overline c}\, 
x_{a \overline c}  \,  x_{b \overline b} 
\end{equation} 
and 
\begin{equation} 
\label{eq:exchange-relation-D-5} 
\tilde x_{a \overline a} x_{b \overline c} = 
\tilde p^+_{a \overline a,b \overline c}\, 
x_{ab}\, \tilde x_{c \overline c} 
+ 
\tilde p^-_{a \overline a,b \overline c}\, 
x_{a \overline c}  \, \tilde x_{b \overline b} 
\end{equation} 
whenever $a,b,c,\overline a$ are 
in counter-clockwise order. 
\end{proposition} 
 
We define a special coefficient system and the corresponding 
cluster algebra $\AA_\circ$ of type $D_n$ in 
precisely the same way as for the types $B_n$ and $C_n$ above. 
We conclude this paper with a geometric realization of this 
algebra similar to the one given in 
Example~\ref{ex:realization-Bn}. 
 
\begin{example} 
[\emph{Geometric realization for $\AA_\circ$ in type~$D_n$}] 
\label{ex:realization-Dn} 
Consider the same variety $X_{n+2}$ as in 
Example~\ref{ex:realization-Bn}. 
Let $X$ be the divisor in $X_{n+2}$ given by the equation 
$\Delta_{n+1,n+2}= 0$; thus, we have 
$$\CC[X] = \CC[X_{n+2}]/\langle \Delta_{n+1,n+2} \rangle \ .$$ 
(Geometrically, $X$ is the affine cone over 
the Schubert divisor in the Grassmannian ${\rm Gr}_{2,n+2}$.) 
Let $\ZZ[X]$ denote the $\ZZ$-form of $\CC[X]$ generated by all Pl\"ucker 
coordinates. 
 
By analogy with (\ref{eq:Delta-functions-Bn}), we introduce 
the following family of functions from $\ZZ[X]$: 
\begin{equation} 
\label{eq:Delta-functions-Dn} 
\begin{array}{rcll} 
[a, \overline a]  &\mapsto& \Delta_{a \overline a} = \Delta_{a,n+1}, 
&\qquad (1 \leq a \leq n),\\[.05in] 
\widetilde{[a, \overline a]}  &\mapsto& 
\tilde \Delta_{a \overline a} = \Delta_{a,n+2} 
&\qquad (1 \leq a \leq n),\\[.05in] 
\{[a,b], [\overline a, \overline b]\} &\mapsto& \Delta_{ab} 
 &\qquad (1 \leq a < b \leq n),\\[.05in] 
\{[a,\overline b], [\overline a, b]\} &\mapsto& 
\Delta_{a \overline b} = 
\Delta_{a,n+1} \Delta_{b,n+2} - \Delta_{ab} 
&\qquad 
(1 \leq a < b \leq n). 
\end{array} 
\end{equation} 
(Note that in $\ZZ[X]$, there is a relation 
$\Delta_{a,n+1} \Delta_{b,n+2} = \Delta_{a,n+2} \Delta_{b,n+1}$ 
since we have $\Delta_{n+1,n+2} = 0$.)

\begin{proposition} 
\label{pr:geom-realization-Dn} 
The correspondence sending each cluster variable 
and each coefficient for the cluster algebra $\AA_\circ$ 
of type $D_n$ to the element in {\rm (\ref{eq:Delta-functions-Dn})} 
with the same label extends uniquely to an algebra isomorphism between 
$\AA_\circ$ and~$\ZZ[X]$. 
\end{proposition} 
 
\proof 
The proof is completely analogous 
to that of Proposition~\ref{pr:geom-realization-Bn}. 
Details are left to the reader. 
\endproof 
 
\end{example}

\section*{Acknowledgments} 
 
We thank Tara Holm and Egon Schulte for bibliographical and 
terminological guidance, 
and Brian Conrad and Peter Magyar for useful discussions. 
We gratefully acknowledge the support 
provided by the Isaac Newton Institute in Cambridge, UK, in Spring 
2001. 
We thank Bill Fulton for 
generously supporting Andrei Zelevinsky's visit to 
the University of Michigan at Ann Arbor in May 2002.

\end{document}